\documentclass[11pt]{article}
\usepackage{amsfonts}
\usepackage{amssymb}
\usepackage{amsmath}
\usepackage{color}
=0pt
\def\sqr#1#2{{\vcenter{\vbox{\hrule height.#2pt
              \hbox{\vrule width.#2pt height#1pt \kern#1pt \vrule
width.#2pt}
              \hrule height.#2pt}}}}
\def\dbC{{\mathbb{C}}}
\def\dbD{{\mathbb{D}}}
\def\dbE{{\mathbb{E}}}
\def\dbF{{\mathbb{F}}}

\def\dbL{{\mathbb{L}}}

\def\dbN{{\mathbb{N}}}

\def\dbP{{\mathbb{P}}}

\def\dbR{{\mathbb{R}}}

\def\a{\alpha}

\def\g{\gamma}

\def\e{\varepsilon}

\def\si{\sigma}
\def\t{\tau}
\def\f{\varphi}
\def\th{\theta}
\def\om{\omega}

\def\3n{\negthinspace \negthinspace \negthinspace }
\def\2n{\negthinspace \negthinspace }
\def\1n{\negthinspace }
\def\ns{\noalign{\smallskip} }

\def\ds{\displaystyle}
\def\G{\Gamma}
\def\D{\Delta}
\def\Th{\Theta}
\def\L{\Lambda}
\def\Si{\Sigma}

\def\Om{\Omega}
\def\om{\omega}
\def\cA{{\cal A}}
\def\cB{{\cal B}}
\def\cC{{\cal C}}
\def\cD{{\cal D}}

\def\cF{{\cal F}}

\def\cH{{\cal H}}

\def\cJ{{\cal J}}
\def\cK{{\cal K}}
\def\cL{{\cal L}}

\def\cO{{\cal O}}

\def\cR{{\cal R}}
\def\cS{{\cal S}}

\def\cU{{\cal U}}
\def\cV{{\cal V}}

\def\cX{{\cal X}}
\def\cY{{\cal Y}}

\def\cl{{\cal l}}
\def\mE{{\mathbb{E}}}

\def\no{\noindent}

\def\ss{\smallskip}
\def\ms{\medskip}
\def\bs{\bigskip}
\def\q{\quad}
\def\qq{\qquad}
\def\hb{\hbox}

\def\lan{\mathop{\big\langle}}
\def\ran{\mathop{\big\rangle}}

\def\pa{\partial}

\def\wt{\widetilde}
\def\cd{\cdot}
\def\cds{\cdots}

\def\dim{\hbox{\rm dim$\,$}}

\def\ae{\hbox{\rm a.e.{ }}}
\def\as{\hbox{\rm a.s.{ }}}

\def\span{\hbox{\rm span$\,$}}

\def\cl{\overline}

\def\({\Big (}
\def\){\Big )}
\def\[{\Big[}
\def\]{\Big]}

\def\={\buildrel \triangle \over =}
\def\wh{\widehat}

\def\resp{{\it resp. }}

\def\be{\begin{equation}}
\def\bel{\begin{equation}\label}
\def\ee{\end{equation}}
\def\bea{\begin{eqnarray}}
\def\eea{\end{eqnarray}}
\def\bt{\begin{theorem}}
\def\et{\end{theorem}}
\def\bc{\begin{corollary}}
\def\ec{\end{corollary}}
\def\bl{\begin{lemma}}
\def\el{\end{lemma}}
\def\bp{\begin{proposition}}
\def\ep{\end{proposition}}
\def\br{\begin{remark}}
\def\er{\end{remark}}
\def\ba{\begin{array}}
\def\ea{\end{array}}
\def\bd{\begin{definition}}
\def\ed{\end{definition}}

\newtheorem{lemma}{Lemma}[section]
\newtheorem{remark}{Remark}[section]
\newtheorem{example}{Example}[section]
\newtheorem{theorem}{Theorem}[section]
\newtheorem{corollary}{Corollary}[section]

\newtheorem{definition}{Definition}[section]
\newtheorem{proposition}{Proposition}[section]

\allowdisplaybreaks
 \makeatletter
   
   \@addtoreset{equation}{section}
\makeatother

\begin{document}

\title{\bf Optimal Feedback for Stochastic Linear Quadratic Control and Backward Stochastic Riccati
Equations in Infinite Dimensions}

\author{Qi L\"u\thanks{School of
Mathematics, Sichuan University, Chengdu,
610064, China.  {\small\it E-mail:} {\small\tt
lu@scu.edu.cn}.\ms} \q and \q Xu Zhang\thanks{
School of Mathematics, Sichuan University,
Chengdu, 610064, China.  {\small\it E-mail:}
{\small\tt zhang\_xu@scu.edu.cn}.} }

\date{}

\maketitle

\begin{abstract}
It is a longstanding unsolved problem to
characterize the optimal feedbacks for general
SLQs (i.e., stochastic linear quadratic control
problems) with random coefficients in infinite dimensions; while the same problem but in finite dimensions was just addressed in a recent work \cite{LWZ}. This paper
is devoted to giving a solution to this
problem under some assumptions which can be verified for several interesting concrete models. More precisely, under these  assumptions, we establish the
equivalence between the existence of optimal
feedback operator for infinite dimensional SLQs
and the solvability of the corresponding
operator-valued, backward stochastic Riccati
equations. A key contribution of this work is to introduce a suitable notion of solutions (i.e., transposition solutions to the aforementioned Riccati equations), which plays
a crucial role in both the statement and the
proof of our main result.

\end{abstract}

\bs

\no{\bf 2010 Mathematics Subject
Classification}. Primary 60H15; Secondary 93E20,
60H25, 49J30.

\bs

\no{\bf Key Words}. Backward stochastic Riccati
equation, stochastic linear quadratic problem, infinite dimensions,
optimal feedback operator,  transposition
solution.

\section{Introduction}\label{s1}

Linear quadratic control problems (LQs for
short) are extensively studied in Control
Theory. It is an extremely important class of
optimal control problems because it can model
many problems in applications and, more
importantly, many nonlinear control problems can
be reasonably approximated by LQs.  It is
well-known that, one of the three milestones in
modern (finite dimensional) optimal control
theory is Kalman's theory for LQs
(\cite{Kalman}, see also \cite{LM, Wonham1} for
some further development).

In Control Theory, one of the fundamental issues
is to find feedback controls, which are
particularly important in practical
applications. Indeed, the main advantage of
feedback controls is that they keep the
corresponding control strategy to be robust with
respect to  (small) perturbation/disturbance,
which is usually unavoidable in real world.
Unfortunately, it is actually very difficult to
find feedback controls for many control
problems. So far, the most successful attempt in
this respect is that for LQs, in particular for
that in the deterministic and finite dimensional
setting, to be recalled as follows.

\subsection{Review on deterministic LQs in finite
dimensions and matrix valued Riccati equations}

Let $T>0$, $n,m\in \dbN$. Denote by
$\cS(\dbR^n)$ the set of all $n\times n$
symmetric matrices. For any  matric $M\in
\dbR^{n\times m}$,  denote by $M^\top$ and $\cR
(M)$ its transpose and range, respecitvely. More
notations (may be used below) will be given in
Chapter 2.

For any given $(s,\eta)\in[0,T]\times\dbR^n$, we
begin with the following control system:
\begin{equation}\label{11.1-eq1}
\left\{
\begin{array}{ll}\ds
\dot
x(t)=\mathrm{A}(t)x(t)+\mathrm{B}(t)u(t),\qq
t\in[s,T],\cr
\ns\ds x(s)=\eta,
\end{array}
\right.
\end{equation}
where $\mathrm{A}(\cd)\in
L^\infty(0,T;\dbR^{n\times n})$ and
$\mathrm{B}(\cd) \in L^\infty(0,T;\dbR^{n\times
m})$, and $u(\cd)\in L^2(s,T;\dbR^m)$ is the
control variable. The cost functional takes the
form:
\begin{equation}\label{11.1-eq2}
\begin{array}{ll}\ds
J(s,\eta;u(\cd))\3n&\ds = \frac{1}{2}\int_s^T
\big(\langle \mathrm{Q}(t) x(t),
x(t)\rangle_{\dbR^n} + \langle \mathrm{R}(t)
u(t),u(t)\rangle_{\dbR^m}
\big)dt\\
\ns&\ds\q+\frac{1}{2}\langle
\mathrm{G}x(T),x(T)\rangle_{\dbR^n},
\end{array}
\end{equation}
with $\mathrm{Q}(\cd)\in
L^\infty(0,T;\cS(\dbR^n))$, $\mathrm{R}(\cd)\in
L^\infty(0,T;\cS(\dbR^m))$ and
$\mathrm{G}\in\cS(\dbR^n)$. Let us consider the
following optimal control problem:

\ms

\no\bf Problem (LQ). \rm For each
$(s,\eta)\in[0,T]\times\dbR^n$, find (if
possible) a $\bar  u(\cd)\in L^2(s,T;\dbR^m)$,
called an {\it optimal control}, such that
\begin{equation}\label{11.1-eq3}
J(s,\eta;\bar  u(\cd))=\inf_{u(\cd)\in
L^2(s,T;\dbR^m)}J(s,\eta;u(\cd)).
\end{equation}
If the above is possible, then {\bf Problem
(LQ)}  is called {\it solvable}. The
corresponding state $\bar  x(\cdot)$ is called
an {\it optimal state}. If the $\bar u(\cd)$
which fulfills \eqref{11.1-eq3} is unique, then
{\bf Problem (LQ)}  is called {\it uniquely
solvable}.

\ms

Throughout this section, we assume that
$\mathrm{R}(\cdot)>\!\!>0$, i.e.,
$\mathrm{R}(t)-c I_m>0$ for some constant $c> 0$
and for a.e. $t\in [0,T]$, where $I_m$ stands
for the identity matrix in $\dbR^{m\times m}$.
One can show that {\bf Problem (LQ)} admits one
and only one optimal control $\bar  u(\cd)$,
which  can be characterized as \bel{2zx.q2--3}
 \bar  u(\cdot)=-\mathrm{R}(\cdot)^{-1}\mathrm{B}(\cdot)^\top\psi(\cdot),
 \ee
where $\psi(\cdot)$ solves \bel{2zx.q21}
 \left\{\ba{ll}
 \dot \psi(t)=-\mathrm{A}(t)^\top\psi(t)-\mathrm{Q}(t)\bar  x(t),\q
t\in[s,T],\\\ns\ds
 \psi(T)=-\mathrm{G}\bar  x(T).
 \ea\right.
 \ee

Note however that the optimal control $\bar
u(\cdot)$ given by (\ref{2zx.q2--3}) is NOT of
feedback form. In order to find the feedback
form of $\bar  u(\cdot)$, formally, assume that
\bel{2zx.q2013} \psi(t)=\mathrm{P}(t)\bar
x(t),\q t\in[s,T],
 \ee
for some $\dbR^{n\times n}$-valued function
$\mathrm{P}(\cdot)$ to be given later. Then,
combining (\ref{2zx.q2--3}) and
(\ref{2zx.q2013}), one obtains the feedback
control:
 \bel{2zx.q2-3}
 \bar  u(t)=-\mathrm{R}(t)^{-1}B(t)^\top
 \mathrm{P}(t)\bar  x(t), \ \mbox{ for a.e. }t\in [s,T].
 \ee
How to obtain the above $\mathrm{P}(\cdot)$?
Differentiating (\ref{2zx.q2013}), we find that
(e.g. \cite{YL})
 $$
 \ba{ll}
 \big(-\mathrm{A}(t)^\top \mathrm{P}(t)-\mathrm{Q}(t)\big)\bar  x(t)\\
 \ns\ds=-\mathrm{A}(t)^\top\psi(t)+\mathrm{Q}(t)\bar  x(t)\\[2mm]
 =\dot \psi(t)=\dot {\mathrm{P}(t)}\bar  x(t)+\mathrm{P}(t)\dot {{\bar  x}}(t)\\[2mm]
 =\dot {\mathrm{P}}(t)\bar  x(t)+\mathrm{P}(t)\big(\mathrm{A}(t)\bar  x(t)+\mathrm{B}(t)\bar  u(t)\big)\\[2mm]
 =\big(\dot {\mathrm{P}}(t)+\mathrm{P}(t)\mathrm{A}(t)-\mathrm{P}(t)\mathrm{B}(t)\mathrm{R}(t)^{-1}\mathrm{B}(t)^\top \mathrm{P}(t)\big)\bar  x(t).
 \ea
 $$
Clearly, it suffices to choose
$\mathrm{P}(\cdot)$ solving the following matrix
valued Riccati equation:
\begin{equation}\label{11.1-eq34}
\left\{
\begin{array}{ll}\ds
\dot{\mathrm{P}}(t)=-\big(\mathrm{P}(t)
\mathrm{A}(t)+\mathrm{A}(t)^\top \mathrm{P}(t)
+\mathrm{Q}(t)\\[1mm]
\ns\ds\qq\q\q - \mathrm{P}(t) \mathrm{B}(t)
\mathrm{R}(t)^{-1} \mathrm{B}(t)^\top
\mathrm{P}(t)\big),\q \ae t\in[0,T],\\[1mm]
\ns\ds \mathrm{P}(T)=\mathrm{G}.
\end{array}
\right.
\end{equation}

Due to many applications in the projective
differential geometry of curves and the calculus
of variations, the study of Riccati (type)
equations may date back to the very early period
of modern mathematics. Some particular cases
were studied more than three hundred years ago
by J. Bernoulli (1654--1705) and J. Riccati
(1676--1754). Other important contributors in
this respect include D. Bernoulli, L.~Euler,
A.-M.~Legendre, J.~d'Alembert and so on. In the
early stage, Riccati equations were in a narrow
sense, i.e., first-order ordinary differential
equations with quadratic unknowns. It is one of
the simplest type of nonlinear differential
equations which may have no explicit solutions.
Later on, the term of Riccati equation is also
used to refer to matrix or operator equations
with analogous quadratic unknowns. These
equations appear in many different mathematical
fields, such as boundary value problems,
scattering theory, spectral factorization of
operators, singular perturbation theory,
differential geometry, Hamiltonian system, the
theory of Bessel functions, functional calculus,
etc (e.g. \cite{Gromoll,Henry,Hermann,Tretter}).
Moreover, these equations are also serves as
fundamental tools in some other natural science
fields, such as quantum chemistry (e.g.
\cite{Fraga,Schuch}).

\ss

Let us recall that the simplest Riccati equation
takes the following form: \bel{eq115} \left\{
\ba{ll}
\dot y=y^2,\q t\in\dbR\\
\ns\ds y(0)=y_0(\in \dbR). \ea\right. \ee The
solution to (\ref{eq115}) is given by
 \bel{eq1151}
 y(t)=\frac{1}{y_0^{-1}-t},
 \ee
which blows up at $t=y_0^{-1}$. Nevertheless, by
means of the special structure of
(\ref{11.1-eq34}) (because it comes from {\bf
Problem (LQ)}), Kalman proved the following
fundamental result (\cite{Kalman}):

\begin{theorem}\label{sec1-th1}
{\bf Problem (LQ)} is uniquely solvable if and
only if the corresponding Riccati equation
(\ref{11.1-eq34}) admits a unique solution
$\mathrm{P}(\cdot)\in C([0,T];\cS(\dbR^n))$.
\end{theorem}

Note that the standard {\bf Problem (LQ)}, i.e.
the case when $\mathrm{R}(\cdot)>\!\!>0$,
$\mathrm{Q}(\cdot)\ge 0$ and $G\ge 0$, is always
uniquely solvable, and therefore, by Theorem
\ref{sec1-th1}, the  equation (\ref{11.1-eq34})
admits a unique global solution
$\mathrm{P}(\cdot)$ on $[0,T]$.

After the seminal work in \cite{Kalman}, the
matrix-valued Riccati equations were broadly
applied to solving control problems, such as
linear optimal control and filtering problems
with quadratic cost functionals, feedback
stabilization problems, linear dynamic games
with quadratic cost functionals, etc. We refer
to \cite{AFIJ, Paraev, Reid, Wonham1, Zelikin}
and the rich references therein for more
details.

Stimulated by Kalman's work, LQs are studied
extensively for more general control systems,
such as those governed by partial differential
equations (PDEs for short), by stochastic
differential equations (SDEs for short) and
stochastic partial differential equations (SPDEs
for short).

As far as we know, \cite{Lions} is the earliest
monograph addressing systematically to LQs for
PDEs. Compared with {\bf Problem (LQ)},  the
infinite dimensional setting brings about many
new difficulties (even for the deterministic
case). To overcome them, people introduce
several new tools, say (deterministic)
operator-valued differential Riccati equations
(see \cite{Curtain, LT, LY, Lions} for the
details). Nevertheless, the main concern of our
present work is to study LQs in stochastic
setting, especially those governed by SPDEs.

\subsection{Review on SLQs in finite dimensions}

Let $(\Om,\cF,\mathbf{F},\dbP)$ be a complete
filtered probability space with
$\mathbf{F}=\{\cF_t\}_{t\in[0,T]}$ (satisfying
the usual conditions), on which a
$1$-dimensional standard Brownian motion
$\{W(t)\}_{t\in[0,T]}$ is defined such that
$\mathbf{F}$ is the natural filtration generated
by $W(\cd)$ (augmented by all $\dbP$-null sets).

For any $(s,\eta)\in [0,T)\times
L^2_{\cF_s}(\Om;\dbR^n)$, consider the following
controlled linear stochastic (ordinary)
differential equation, i.e., stochastic
evolution equation (SEE for short) in finite
dimensions:
\begin{equation}\label{f5.2-eq1}
\left\{\begin{array}{ll}\ds dx(t)=\big(A x(t) +
B u(t)\big)dt + \big(C x(t)+D
u(t)\big)dW(t) & \mbox{ in }(s,T],\\
\ns\ds x(s)=\eta,
\end{array}
\right.
\end{equation}
with the quadratic cost functional:
\begin{equation}\label{f5.2-eq2}
\begin{array}{ll}\ds \cJ(s,\eta;u(\cd))
\3n&\ds=\frac{1}{2}\mE\Big[ \int_s^T
\big(\big\langle Q x(t),x(t)\big\rangle_{\dbR^n}
+\big\langle R
u(t),u(t)\big\rangle_{\dbR^m}\big)dt\\
\ns&\ds\qq\qq\q + \langle
Gx(T),x(T)\rangle_{\dbR^n}\Big],
\end{array}
\end{equation}
where the coefficients
$$
\begin{array}{ll}\ds
A\in L^\infty_\dbF(\Om;L^1(0,T;\dbR^{n\times
n})),\q B\in
L^\infty_\dbF(\Om;L^2(0,T;\dbR^{n\times m})),\\
\ns\ds C\in
L^\infty_\dbF(\Om;L^2(0,T;\dbR^{n\times n})), \q
D\in L^\infty_\dbF(0,T;\dbR^{n\times m}),\\
\ns\ds Q\in
L^\infty_\dbF(\Omega;L^1(0,T;\cS(\dbR^n))), \q R
\in L^\infty_\dbF(0,T;\cS(\dbR^m)),\q G\in
L^\infty_{\cF_T}(\Om;\cS(\dbR^n)).
\end{array}
$$
To simplify the notations, here and henceforth
the sample point $\omega(\in\Omega)$ and/or the
time variable $t(\in [0,T])$ in the coefficients
are often suppressed, in the case that no
confusion would occur.

In \eqref{f5.2-eq1}, $ u(\cd)\in
L^2_\dbF(s,T;\dbR^m)$ is the control variable,
$x(\cd)(=x(\cd;s,\eta,u(\cd)))\in
C_\dbF([s,T];L^2(\Om;$ $\dbR^n))$ is the state
variable.

\ss

Let us consider the following  SLQ in finite
dimensions:

\ms

\no\bf Problem (FSLQ): \rm For each $(s,\eta)\in
[0,T)\times L^2_{\cF_s}(\Om;\dbR^n)$, find (if
possible) a control $\bar  u(\cd)\in
L^2_\dbF(s,T;\dbR^m)$, called an {\it optimal
control}, such that
\begin{equation}\label{f5.2-eq3}
\cJ\big(s,\eta;\bar  u(\cd)\big)=\inf_{u(\cd)\in
L^2_\dbF(s,T;\dbR^m)}\cJ\big(s,\eta;u(\cd)\big).
\end{equation}
If the above is possible, then {\bf Problem
(FSLQ)} is  called {\it solvable}. If the $\bar
u(\cd)$ which fulfills \eqref{f5.2-eq3} is
unique, then {\bf Problem (FSLQ)}  is called
{\it uniquely solvable}.

\ms

SLQs in finite dimensions have been extensively
studied in the literature (See \cite{AMZ,
Athens, Bismut1, Bismut2, CLZ1, Tang1, Wonham2,
YZ} and the rich references therein). Similar to
the deterministic setting, Riccati equations
(and their variants) are introduced as
fundamental tools for constructing feedback
controls. Nevertheless, for stochastic problems
one usually has to consider backward stochastic
Riccati equations (BSREs for short). For {\bf
Problem (FSLQ)}, the desired (matrix-valued)
BSRE takes the following form:
\begin{equation}\label{f5.5-eq6}
 \left\{
\begin{array}{ll}\ds
dP =-\big( PA +
A^\top P + \Lambda C + C^\top \Lambda\\[1mm]
\ns\ds\qq\;\;\q + C^\top PC  + Q - L^\top
K^{\dag} L
\big)dt+ \Lambda dW(t) \q&\mbox{in }[0,T],\\
\ns\ds P(T)=G,
\end{array}
\right.
\end{equation}
where $K= R+D^\top PD$, $ L= B^\top P+D^\top
(PC+\Lambda)$, and $K^{\dag}$ denotes the
Moore-Penrose pseudo-inverse of $K$.

To the authors' best knowledge, \cite{Wonham2}
is the first work which employed Riccati
equations to study {\bf Problem (FSLQ)}. After
\cite{Wonham2}, Riccati equations were
systematically applied to study SLQs (e.g.
\cite{Athens, Bensoussan1, Bismut2, Davis, YZ}),
and the well-posedness of such equations was
studied in some literatures (See \cite{Tang1,
YZ} and the references therein). In the early
works in this respect (e.g., \cite{CLZ1,
Wonham2, YZ}),  the coefficients $A$, $B$, $C$,
$D$, $Q$, $R$, $G$ appeared in the control
system \eqref{f5.2-eq1} and the cost functional
\eqref{f5.2-eq2} were assumed to be
deterministic matrices. For this case, the
corresponding BSRE \eqref{f5.5-eq6} is
deterministic (i.e., $\Lambda\equiv0$ in
\eqref{f5.5-eq6}), as well. On the other hand,
\cite{Bismut1} is the first work addressed to
the study of SLQs (in finite dimensions) with
random coefficients. In \cite{Bismut1, Bismut2},
the equation \eqref{f5.5-eq6} was formally
derived. However, at that time only some special
and simple cases could be solved. Later,
\cite{Peng} proved the well-posedness of this
equation under the condition that $D=0$. This
condition was dropped in \cite{Tang1}. Now SLQs
in finite dimensions
 have been
studied in many literature, for which we refer
the readers to \cite{AMZ, Athens, Bismut1,
Bismut2, CLZ1, CLZ2, HZ, SY, SLY, Tang1,
Wonham2, YZ} and the rich references cited
therein.

Naturally, one hopes to employ solutions of
(\ref{f5.5-eq6}) to construct feedback controls
for {\bf Problem (FSLQ)}. Indeed, under some
mild assumptions, by Theorem \ref{sec1-th1}, one
can always find the desired feedback control
through the corresponding Riccati equation
whenever a deterministic LQ, i.e. {\bf Problem
(LQ)} is solvable.  However, as pointed out in
\cite{LWZ}, significantly different from its
deterministic counterpart, the problem of
feedback controls is much less well-understood
for {\bf Problem (FSLQ)}. Indeed, many new
difficulties are met:
\begin{itemize}
  \item The corresponding Riccati equation (\ref{f5.5-eq6}) is a
  backward stochastic differential equation
  with a quadratic nonlinear term. Compared with the equation \eqref{11.1-eq34}, generally speaking, the
  well-posedness of (\ref{f5.5-eq6}) is highly nontrivial.

  \item Although the BSRE \eqref{f5.5-eq6} is solvable whenever so is {\bf Problem (FSLQ)} (\cite{Tang1}), as noted in \cite[Remark 1.2]{LWZ}, the corresponding solution to \eqref{f5.5-eq6}
  may not be
  regular enough to serve as the design of feedback controls.
\end{itemize}

Due to the last difficulty mentioned above, it
is quite natural to ask such a question:

\ms

\no{\bf Question (Q):} {\it Is it possible to
link the existence of optimal feedback controls
(rather than the solvability) for SLQs directly
to the solvability of the corresponding BSRE?}

\ms

Clearly, from the viewpoint of applications, it
is more desirable to study the existence of
feedback controls for SLQs than the solvability
of the same problems. An affirmative answer to
the above question in finite dimensions was
given in \cite{LWZ}. More precisely, let us
recall first the notion of optimal feedback
operator for {\bf Problem (FSLQ)}:

\begin{definition}\label{5.7-def1.1}
A stochastic process $\Th(\cd)\in
L^\infty_\dbF(\Om;L^2(0,T;\dbR^{m\times n}))$ is
called an {\it optimal feedback  operator} for
{\bf Problem (FSLQ)} if
\begin{equation*}
\begin{array}{ll}\ds \cJ(s,\eta;\Th(\cd)\bar  x(\cd))\leq
\cJ(s,\eta;u(\cd)),\\[2mm]
\ns\ds\qq \forall\; (s,\eta)\in [0,T)\times
L^2_{\cF_s}(\Om;\dbR^n),\;\; u(\cd)\in
L^2_\dbF(s,T;\dbR^m),
\end{array}
\end{equation*}
where $\bar  x(\cd)=\bar  x(\cd\,;s,\eta,
\Th(\cd)\bar  x(\cd))$ solves the equation
(\ref{f5.2-eq1}) with $u(t)$ replaced by $\Th
(t)\bar x(t)$.
\end{definition}

The choice of the optimal feedback operator
space $L^\infty_\dbF(\Om;L^2(0,T;\dbR^{m\times
n}))$ (in Definition \ref{5.7-def1.1}) is sharp.
Indeed, this is a consequence of the following
sharp well-posedness result (See \cite[Chapter
V, Section 3]{Protter} for its proof):
\begin{lemma}\label{fflm2}
Let $\cA\in
L^\infty_\dbF(\Om;L^1(0,T;\dbR^{n\times n}))$,
$\cB\in L^\infty_\dbF(\Om;L^2(0,T;\dbR^{n\times
n}))$, $f\in L^1_\dbF(s,T;\dbR^n)$ and $g\in
L^2_\dbF(s,T;\dbR^n)$. Then, for any $\eta\in
L^2_{\cF_s}(\Om;\dbR^n)$, the stochastic
differential equation
$$
\left\{
\begin{array}{ll}\ds
dx = (\cA x  + f)dt + (\cB x+g)dW(t) &\mbox{ in }[s,T],\\
\ns\ds x(s)=\eta
\end{array}
\right.
$$
admits one and only one $\dbF$-adapted solution
$x(\cd)\in L^2_\dbF(\Om; C([s,T];\dbR^n))$.
\end{lemma}

In \cite{LWZ}, the following equivalence
(between the existence of optimal feedback
operator for {\bf Problem (FSLQ)} and the
solvability of the corresponding BSRE
\eqref{f5.5-eq6} in a suitable sense) was
proved.

\begin{theorem}\label{5.70000-th1}
\textbf{Problem (FSLQ)} admits an optimal
feedback operator $\Th(\cd)\in
L^\infty_\dbF(\Om;L^2(0,T;$ $\dbR^{m\times n}))$
if and only if the BSRE \eqref{f5.5-eq6} admits
a solution
$$\big(P(\cd),\L(\cd)\big) \in
L^{\infty}_{\dbF}(\Om;C([0,T];\cS(\dbR^n)))
\times L^p_{\dbF}(\Omega;L^2(0,T;\cS(\dbR^n)))$$
(for all $p\geq 1$)
such that
$$
 \cR(K(t,\omega))\supset\cR(L(t,\omega)) \q\hbox{and}\q K(t,\omega)\geq 0,\qq \ae (t,\omega)\in
 [0,T]\times\Omega,
$$
and
$$
K(\cd)^{\dag}L(\cd)\in
L^\infty_\dbF(\Om;L^2(0,T;\dbR^{m\times n})).
$$
In this case, the optimal feedback operator
$\Th(\cd)$ is given as
$$
\Th(\cd)=-K(\cd)^{\dag}L(\cd) + \big(I_m -
K(\cd)^{\dag}K(\cd)\big)\th,
$$
where $\th\in
L^\infty_\dbF(\Om;L^2(0,T;\dbR^{m\times n}))$ is
arbitrarily given. Furthermore,
$$
\inf_{u\in
L^2_\dbF(s,T;\dbR^m)}\cJ(s,\eta;u)=\frac{1}{2}\,\dbE\langle
P(s)\eta,\eta\rangle_{\dbR^n}.
$$
\end{theorem}

The following example (which is a modification
of \cite[Example 6.2]{LWZ}), shows that a
solvable {\bf Problem (FSLQ)} does not
necessarily  admit any optimal feedback
operator:
\begin{example}\label{counterexample-1}
Consider a solvable {\bf Problem (FSLQ)} with
the data:
\bel{201611195}
\begin{cases}\ds
m=n=1, \\
\ns\ds A=B=C=Q=0,\ \ D=1,\\
\ns\ds R=\frac{1}{4}>0,\ \
G=Y(T)^{-1}-\frac{1}{4}>0,
\end{cases}
\ee
where $Y(\cd)$ will be given later. We shall
show that for a suitable chosen $Y(\cd)$, the
resulting {\bf Problem (FSLQ)} does not  admit
optimal feedback operators.

Define two (one-dimensional) stochastic
processes $M(\cd)$ and $\zeta(\cd)$ and a
stopping time $\t$ as follows:
\begin{equation}\label{11.12-eq1}
\begin{cases}
\ds M(t)\triangleq
\int_0^t\frac{1}{\sqrt{T-s}}dW(s),\qq
t\in[0,T),\\[3mm]
\ns\ds \t\triangleq \inf\big\{t\in[0,T)\;\big|\; |M(t)|>1\big\},\\[3mm]
\ns\ds \zeta(t)\triangleq
\frac{\pi}{2\sqrt{2}\sqrt{T-t}}\chi_{[0,\t]}(t),\qq
t\in[0,T).
\end{cases}
\end{equation}
Here we agree that $\inf\emptyset=T$, also
$\chi_{[0,\t]}(\cdot)$ stands for the
characteristic function of $[0,\t]$. Clearly, by
the definition of $\tau$, it follows that, for
any $t\in [0,T]$,
\bel{boundedness-example}\ba{ll}
\ns\ds
\Big|\int_0^t\zeta(s)dW(s)\Big|=\frac{\pi}{2\sqrt{2}}\Big|\int_0^{\min
(t,\tau)}\frac{1}{ \sqrt{T-s}} dW(s)\Big|\leq
\frac{\pi}{2\sqrt{2}}.
\ea\ee
Further, It was shown in \cite[Lemma
A.1]{FreidosReis} that
\bel{aim-1}\ba{ll}
\ns\ds
\dbE\Big[\exp\Big(\int_0^T|\zeta(t)|^2dt\Big)\Big]=\infty.
\ea\ee

Consider the following backward stochastic
differential equation:
$$
Y(t)=\int_0^T\zeta(s)dW(s)+\frac{\pi}{2\sqrt{2}}+1-\int_t^TZ(s)dW(s),\q
t\in[0,T].
$$
This equation admits a unique solution $(Y,Z)$
as follows
$$
\begin{cases}\ds
Y(t)=\int_0^t\zeta(s)dW(s)+\frac{\pi}{2\sqrt{2}}+1,\\
\ns\ds
 Z(t)=\zeta(t),
\end{cases}
 \qquad t\in[0,T].
$$
From \eqref{11.12-eq1}--\eqref{aim-1}, it is
easy to see that
\begin{equation}\label{11.12-eq2}
\begin{cases}\ds
1\leq Y(\cd)\leq \frac{\pi}{\sqrt{2}}+1, \q
Z(\cd) \in
L^p_{\dbF}(\Omega;L^2(0,T;\dbR))\hbox{ for any }p\ge 1,\\[2mm]
\ns\ds Z(\cd) \notin
L^{\infty}_{\dbF}(\Omega;L^2(0,T;\dbR)).
\end{cases}
\end{equation}

The corresponding BSRE is specialized as
\begin{equation}\label{Riccati-equation-example00}
\left\{\ba{ll}
dP=(R+P)^{-1}\L^2dt+\L dW(t) \q \mbox{in }[0,T],\\[3mm]
P(T)=G, \ea\right.
\end{equation}
and $\Th(\cd)=-(R+P(\cd))^{-1}\L(\cd)$. By
applying It\^o's formula to $Y(\cd)^{-1}$, it is
easy to show that
$$ (P(\cd), \L(\cd))\triangleq ( Y(\cd)^{-1}-R, -Y(\cd)^{-2}Z(\cd))
$$
is the unique solution to
(\ref{Riccati-equation-example00}).

Now, by the contradiction argument, we suppose
the {\bf Problem (FSLQ)} under consideration
admitted an optimal feedback operator
$\Th(\cd)\in L^\infty_\dbF(\Om;L^2(0,T;\dbR))$.
Then, from Theorem \ref{5.70000-th1}, this
feedback control operator would be given
explicitly by
$$\Th(\cd)=-(R+P(\cd))^{-1}\L(\cd)=-Y(\cd)^{-1}Z(\cd).$$ By
\eqref{11.12-eq2}, we see that $\Th(\cd)$ does
not belong to
$L^{\infty}_{\dbF}(\Omega;L^2(0,T;\dbR))$, which
is a contradiction.
%
\end{example}

\subsection{SLQs in infinite dimensions}

From now on, we focus on SLQs in infinite
dimensions, the main concern of this work. In
order that our results can be applied to  as
many different SPDEs as possible, we shall
formulate the problem for SEEs.

Let $H$ and $U$ be separable Hilbert spaces, and
$A$ be an unbounded linear operator (with domain
$D(A)\subset H$), which generates a
$C_0$-semigroup $\{e^{At}\}_{t\geq 0}$. Denote
by $A^*$ the adjoint operator of $A$. More
assumptions (may be used below) will be given in
Chapter 2.

\ss

For any $(s,\eta)\in [0,T)\times
L^2_{\cF_s}(\Om;H)$, we consider the following
controlled linear SEE:
\begin{equation}\label{5.2-eq1}
\left\{\2n\begin{array}{ll}\ds
dx(t)=\big[(A+A_1) x(t) + B u(t)\big]dt + \big(C
x(t)+D
u(t)\big)dW(t) \q \mbox{in }(s,T],\\
\ns\ds x(s)=\eta,
\end{array}
\right.
\end{equation}
with the quadratic cost functional
\begin{equation}\label{5.2-eq2}
\begin{array}{ll}\ds \cJ(s,\eta;u(\cd))
\3n&\ds=\frac{1}{2}\mE\Big[ \int_s^T
\big(\big\langle Q x(t),x(t)\big\rangle_H
+\big\langle R
u(t),u(t)\big\rangle_U\big)dt\\
\ns&\ds\qq\qq\q + \langle
Gx(T),x(T)\rangle_H\Big].
\end{array}
\end{equation}
Here the coefficients $A_1$, $B$, $C$, $D$, $Q$
and  $R$ are suitable operator-valued stochastic
processes, and $G$ is a suitable operator-valued
random variable. In \eqref{5.2-eq1}, $ u(\cd)\in
L^2_\dbF(s,T;U)$ is the control variable,
$x(\cd)(=x(\cd;s,\eta,u(\cd)))\in
C_\dbF([s,T];L^2(\Om;H))$ is the state variable.

Let us consider the following optimal control
problem:

\ms

\no\bf Problem (SLQ): \rm For each $(s,\eta)\in
[0,T)\times L^2_{\cF_s}(\Om;H)$, find (if
possible) a control $\bar  u(\cd)\in
L^2_\dbF(s,T;U)$ such that
\begin{equation}\label{5.2-eq3}
\cJ\big(s,\eta;\bar  u(\cd)\big)=\inf_{u(\cd)\in
L^2_\dbF(s,T;U)}\cJ\big(s,\eta;u(\cd)\big).
\end{equation}
If the above is possible, then {\bf Problem
(SLQ)} is called {\it solvable}. Any $\bar
u(\cdot)$ satisfying (\ref{5.2-eq3}) is called
an {\it optimal control}. If the $\bar u(\cd)$
which fulfills \eqref{5.2-eq3} is unique, then
{\bf Problem (SLQ)}  is called {\it uniquely
solvable}. The corresponding state $\bar
x(\cdot)$ is called an {\it optimal state}, and
$\big(\cl x(\cdot),\bar u(\cdot)\big)$ is called
an {\it optimal pair}.

\begin{remark}
In this work, in order to present the key idea
in the simplest way,  we assume that $W(\cd)$ is
a $1$-dimensional standard Brownian motion. One
can also deal with the case that $W(\cd)$ is a
cylindrical (or other vector-valued) Brownian
motion by the method developed in this paper.
\end{remark}

Similarly to (\ref{f5.5-eq6}), we introduce the
following operator-valued (or more precisely,
$\cL(H)$-valued) BSRE for our {\bf Problem
(SLQ)}:
\begin{equation}\label{5.5-eq6}
 \left\{
\begin{array}{ll}\ds
dP =-\big[ P(A+A_1) +
(A+A_1)^* P + \Lambda C + C^* \Lambda\\[1mm]
\ns\ds\qq\;\;\q + C^* PC  + Q - L^* K^{-1} L
\big]dt+ \Lambda dW(t) \q&\mbox{in }[0,T),\\
\ns\ds P(T)=G,
\end{array}
\right.
\end{equation}
where
\begin{equation}\label{9.7-eq10}
K\equiv R+D^*PD>0, \qq L= B^* P+D^*
(PC+\Lambda).
\end{equation}
Besides the difficulties mentioned in the last
section for the case of finite dimensions, there
exists an essentially new one in the study of
\eqref{5.5-eq6} when $\dim H=\infty$, without
further assumption on the data $A_1$, $B$, $C$,
$D$, $Q$, $R$ and $G$. Indeed, in the infinite
dimensional setting, although $\cL(H)$ is still
a Banach space, it is neither reflexive
(needless to say to be a Hilbert space) nor
separable even if $H$ itself is separable. As
far as we know, in the previous literatures
there exists no such a stochastic
integration/evolution equation theory in general
Banach spaces that can be employed to treat the
well-posedness of \eqref{5.5-eq6}, especially to
handle the (stochastic integral) term ``$\Lambda
dW(t)$" effectively. For example, the existing
results on stochastic integration/evolution
equations in UMD Banach spaces (e.g.
\cite{vanNeerven1, vanNeerven2}) do  not fit the
present case because, if a Banach space is UMD,
then it is reflexive.

\ss

Because of the  aforementioned difficulty, there
exist only a quite limited number of works
 dwelling on some special cases of SLQs in
infinite dimensions (e.g., \cite{Ahmed, GT1,
GT2, Hafizoglu, Ichikawa, Lu, Tessitore}). We
list below some of these typical works:

\begin{itemize}
  \item In \cite{Ichikawa, Tessitore},  \textbf{Problem
(SLQ)} was studied under a key assumption that
the diffusion term in \eqref{5.2-eq1} is
$Cx(t)dW_1(t)+ D u(t)dW_2(t)$, where $W_1(\cd)$
and $W_2(\cd)$ are mutually independent Brownian
motions. This assumption plays a crucial role in
these papers. Indeed, under such an assumption,
the corresponding Riccati equation takes the
form:
\begin{equation}\label{6.19-eq3}
\left\{
\begin{array}{ll}\ds
dP=-\big[ P(A+A_1) + (A+A_1)^* P + C^* PC +
Q\\[1mm]
\ns\ds\qq\q\;\; - PB K^{-1}
B^*P \big]dt &\mbox{ in }[0,T),\\
\ns\ds P(T)=G.
\end{array}
\right.
\end{equation}
The equation \eqref{6.19-eq3} is a random
operator-valued Riccati equation (rather than
operator-valued BSRE), whose well-posedness is
not hard to be obtained.

\ss

 \item When the diffusion term in
\eqref{5.2-eq1} is of the form $\si dW(t)$ with
$\si$ a suitable $\mathbf{F}$-adapted $H$-valued
process, \cite{Ahmed} studied \textbf{Problem
(SLQ)} and found the optimal feedback control by
solving a random operator-valued Riccati
equation (similar to \eqref{6.19-eq3}) and a
backward stochastic evolution equation (BSEE for
short).

\ss

\item In \cite{GT1}, {\bf Problem (SLQ)} was considered
for the case that $R=I$, the identity operator
on $U$, and $D=0$ (the latter means that there
is no control in the diffusion term in
\eqref{5.2-eq1}).   In this case, the equation
\eqref{5.5-eq6} is specialized as
\begin{equation}\label{6.19-eq4}
\left\{
\begin{array}{ll}\ds
dP =-\big[P(A+A_1) + (A+A_1)^* P
+ \Lambda C + C^* \Lambda \\[1mm]
\ns\ds\qq\q\;+ C^* PC + Q - PB B^*P \big]dt +
\Lambda dW(t) &\mbox{in
}[0,T),\\
\ns\ds P(T)=G.
\end{array}
\right.
\end{equation}
Although \eqref{6.19-eq4} looks much simpler
than \eqref{5.5-eq6}, it is also an
operator-valued BSEE (because the ``bad" term
``$\Lambda dW(t)$" is still in
\eqref{6.19-eq4}). Nevertheless, \cite{GT1}
considered a kind of generalized solution to
\eqref{6.19-eq4}, which was a ``weak limit" of
solutions to some suitable finite dimensional
approximations of \eqref{6.19-eq4}. It is shown
in \cite{GT1} that the finite dimensional
approximations $P_n$ of $P$ are convergent in
some weak sense, and via which $P$ may be
obtained as a suitable generalized solution to
\eqref{6.19-eq4} although nothing can be said
about $\Lambda$. This is enough for this special
case that $D=0$. Indeed, the corresponding
optimal feedback operator in \eqref{5.10} (in
the next chapter) is then specialized as
 $$
 \Th(\cd)=-K(\cd)^{-1}B(\cd)^* P(\cd),
 $$
which is independent of $\Lambda$.

\ss

\item In \cite{GT2}, the well-posedness of
\eqref{6.19-eq4} was further studied when $A$ is
a self-adjoint operator on $H$ and there exists
an orthonormal basis $\{e_j\}_{j=1}^\infty$ in
$H$ and an increasing sequence of positive
numbers $\{\mu_j\}_{j=1}^\infty$ so that $ A e_j
= -\mu_j e_j$  for $ j\in\dbN$ and
$\ds\sum_{j=1}^\infty \mu_j^{-r}<\infty$ for
some $ r\in (\frac{1}{4}, \frac{1}{2})$.
Clearly, this assumption is not satisfied by
many controlled stochastic partial differential
equations, such as stochastic wave equations,
stochastic Schr\"odinger equations, stochastic
KdV equations, stochastic beam equations, etc.
It is even not fulfilled by the classical
$m$-dimensional stochastic heat equation for
$m\geq 2$. The well-posedness result (in
\cite{GT2}) for \eqref{6.19-eq4} was then
applied to {\bf Problem (SLQ)} in the case that
$R=I$ and $D=0$.
\end{itemize}

In this paper, we do not impose the assumptions
that $R=I$ and $D=0$. We remark that, dropping
the condition $D=0$ will lead to another
essential difficulty in the study of
\eqref{5.5-eq6}. Indeed, if $D=0$, then the
nonlinear term $L^* K^{-1} L$ is specialized as
$PB B^*P$, which enjoys a ``good" sign in the
energy estimate, and therefore, it is not very
hard to  obtain the well-posedness of the
linearized equation of \eqref{5.5-eq6}, at least
for some special situation. Once the
well-posedness of this linearized equation is
established,  the well-posedness of
\eqref{5.5-eq6} follows from a fixed point
argument (e.g. \cite{GT1, GT2}). However, when
$D\neq 0$, the situation is then quite subtle.
Indeed, it is not an easy task even to derive
the well-posedness of the linearized version of
\eqref{5.5-eq6} because actually one does not
have the desired energy estimate any more.

\ss

Generally speaking, in order to study the
difficult (deterministic/stochastic) nonlinear
partial differential equations, people need to
introduce suitable new concept of solutions,
such as viscosity solutions for Hamilton-Jacobi
equations (\cite{Crandall-Lions}) and fully
nonlinear second-order equations
(\cite{Lions-P}), and renormalized solutions for
the KPZ equation (\cite{Hairer1}). In this work,
we introduce another type of solution, i.e.,
transposition solution to the operator-valued
BSRE \eqref{5.5-eq6}. The concept of
transposition solution (or relaxed transposition
solution) to operator-valued, backward
stochastic (linear) Lyapunov equations was first
introduced in our previous work (\cite{LZ1}),
and was employed to study first and second order
necessary conditions for stochastic optimal
controls in infinite dimensions
(\cite{Frankowska,FZ1, LZZ, LZ1, LZ}). The core
of our stochastic transposition method is to
introduce two vector-valued forward stochastic
(test) evolution equations and via which
solutions to operator-valued backward stochastic
evolution equations can be interpreted in the
sense of transposition (See \cite{LZ} for more
analysis). As we shall see later, transposition
solution is also a suitable notion of solution
to the nonlinear equation \eqref{5.5-eq6}.

\ss

On the other hand, similarly to the setting of
finite dimensions, it is easy to construct a
counterexample (say by modifying suitably
Example \ref{counterexample-1}) in which a {\bf
Problem (SLQ)} (in infinite dimensions) is
solvable, but it does not  admit any optimal
feedback operator (See Definition
\ref{5.7-def1}). Because of this, {\bf Question
(Q)} that we posed in the last section makes
sense in infinite dimensions as well. Actually,
we shall show that, under some assumptions, the
existence of optimal feedback operator for {\bf
Problem (SLQ)} is equivalent to the solvability
of the corresponding operator-valued BSRE
\eqref{5.5-eq6} in the sense of transposition
solution (given in Definition \ref{4.8-def2} in
the next chapter).

\ss

Unsurprisingly, although (deterministic)
operator-valued differential Riccati equations
were introduced to study control problems, they
have been used elsewhere, for instance, the
study of some quantum systems (e.g., \cite{FLO,
Prugovecki}). Because of this, we believe that,
our operator-valued  BSRE (\ref{5.5-eq6}),
though appeared as a tool to study the optimal
feedback controls for {\bf Problem (SLQ)},
should have some independent interest and may be
applied in other places.

\ss

Since 1960s, simulated by a large number of
physical and biological problems (e.g.,
\cite{Carmona-Rozovskii,Kotelenez}), stochastic
partial differential equations or more generally
stochastic evolution equations in infinite
dimensions, have attracted the attention of many
researchers (e.g., \cite{Chow, Prato, Holden}
and the rich references therein). Quite
interestingly, the study of these equations
poses some challenging mathematical problems
(e.g., \cite{Carmona-Rozovskii}). Indeed, for
several infinite dimensional stochastic
equations, whose deterministic counterparts are
very simple and well studied, people has to
develop powerful new tools to study them (e.g.,
\cite{Balazs, Dalang, E, Gubinelli, Hairer1,
Hairer2, Hairer3}). Clearly, (\ref{5.5-eq6}) is
a quadratically nonlinear, operator-valued BSEE,
and therefore, one cannot expect to handle it
easily.

\ss

Because of the very  difficulties of both {\bf
Problem (SLQ)} and the equation (\ref{5.5-eq6}),
the main results in this paper are much less
than satisfactory. In our opinion, this paper
should be a starting point for further studies
in this respect in the future. It seems that
more delicate and powerful tools (especially
that from ``hard analysis") should be
introduced, and most likely the problems should
be analyzed for concrete models one by one.
Nevertheless, we believe that the transposition
solution notion (for (\ref{5.5-eq6})) introduced
in this work should be a basis for these
studies.

\ss

The rest of this work is organized as follows:
Chapter \ref{sec mr} is devoted to stating the
main results of this paper. In Capter \ref{pre},
we give some preliminary results which will be
used later. Chapters \ref{proof} and
\ref{proof-n} are addressed to the proof of our
main result. In  Chapter \ref{sect6}, we shall
prove  the existence of the transposition
solution to (\ref{5.5-eq6}) under some
assumptions on the coefficients. Finally, in
Chapter \ref{sect7}, we shall provide some
concrete illuminating examples.


\section{Statement of the main results}\label{sec
mr}


We begin with some notations to be used
throughout this paper.

Denote by $\dbF$ the progressive $\si$-field (in
$[0,T]\times\Omega$) with respect to
$\mathbf{F}$, by $\dbE f$ the (mathematical)
expectation of an integrable random variable
$f:(\Omega,\cF,\dbP)\to \dbC$, and by $\cC$ a
generic positive constant, which may be
different from one place to another.


Let $\cX$ be a Banach space. For any $t\in[0,T]$
and $p\in [1,\infty)$, denote by
$L_{\cF_t}^p(\Om;\cX)$ the Banach space of all
$\cF_t$-measurable random variables $\xi:\Om\to
\cX$ such that $\mathbb{E}|\xi|_\cX^p < \infty$,
with the canonical  norm. Denote by
$L^{p}_{\dbF}(\Om;C([t,T];\cX))$ the Banach
space of all $\cX$-valued $\mathbf{F}$-adapted
continuous processes $\phi(\cdot)$, with the
norm $$
|\phi(\cd)|_{L^{p}_{\dbF}(\Om;C([t,T];\cX))}
\triangleq
\[\mE\sup_{s\in
[t,T]}|\phi(s)|_\cX^p\]^{1/p}. $$ Similarly, one
can define
$L^{p}_{\dbF}(\Om;C([\tau_1,\tau_2];\cX))$ for
two stopping times $\tau_1$ and $\tau_2$ with
$\tau_1\leq \tau_2$, $\dbP$-a.s. Also, denote by
$C_{\dbF}([t,T];L^{p}(\Om;\cX))$ the Banach
space of all $\cX$-valued $\mathbf{F}$-adapted
processes $\phi(\cdot)$ such that
$\phi(\cdot):[t,T] \to L^{p}_{\cF_T}(\Om;\cX)$
is continuous, with the norm $$
|\phi(\cd)|_{C_{\dbF}([t,T];L^{p}(\Om;\cX))}
\triangleq \sup_{s\in
[t,T]}\left[\mE|\phi(s)|_\cX^p\right]^{1/p}. $$
Fix any $p_1,p_2,p_3,p_4\in[1,\infty]$. Put
$$
\ba{ll}
\ds L^{p_1}_\dbF(\Om;L^{p_2}(t,T;\cX))
=\Big\{\f:(t,T)\times\Om\to
\cX\;\Big|\;\f(\cd)\mbox{
is $\mathbf{F}$-adapted and }\\
\ns\ds\hspace{7.25cm}\dbE\(\int_t^T|\f(s)|_\cX^{p_2}ds\)^{\frac{p_1}{p_2}}<\infty\Big\},\\
\ns\ds
 L^{p_2}_\dbF(t,T;L^{p_1}(\Om;\cX)) =\Big\{\f:(t,T)\times\Om\to
\cX\; \Big|\;\f(\cd)\mbox{ is
$\mathbf{F}$-adapted and
}\\
\ns\ds\hspace{7.25cm}\int_t^T\(\dbE|\f(s)|_X^{p_1}\)^{\frac{p_2}
{p_1}}ds<\infty\Big\}.
 \ea
 $$
(When any one of $p_j$ ($j=1,2,3,4$) is equal to
$\infty$, it is needed to make the usual
modifications in the above definitions of
$L^{p_1}_\dbF(\Om;L^{p_2}(t,T;\cX))$ and
$L^{p_2}_\dbF(0,T;L^{p_1}(\Om;$ $\cX))$).
Clearly, both
$L^{p_1}_\dbF(\Om;L^{p_2}(t,T;\cX))$ and
$L^{p_2}_\dbF(t,T;L^{p_1}(\Om;$ $\cX))$ are
Banach spaces with the canonical norms. If
$p_1=p_2$, we simply write the above spaces as
$L^{p_1}_\dbF(t,T;\cX)$.

Let $\cY$ be another Banach space. Denote by
$\cL(\cX; \cY)$ the Banach space of all bounded
linear operators from $\cX$ to $\cY$, with the
usual operator norm (When $\cY=\cX$, we simply
write $\cL(\cX)$ instead of $\cL(\cX; \cY)$).
Let $\cH$ and $\cU$ be Hilbert spaces. For any
$M\in \cL(\cU; \cH)$, denote by $M^*(\in
\cL(\cH; \cU))$ the dual (operator) of $M$.
Also, denote by $\cS(\cH)$ the set of all
self-adjoint, bounded linear operators on $\cH$,
i.e., $\cS(\cH)=\{M\in \cL(\cH)\;|\;M=M^*\}$. An
$M\in \cS(\cH)$ is called nonnegative (\resp
positive), written as $M\ge 0$ (\resp $M> 0$) if
$(Mh,h)_\cH\ge0$ (\resp $(Mh,h)_\cH\ge
c|h|_\cH^2$ for some constant $c>0$) for any
$h\in\cH$.

Suppose $\cX_j$ and $\cY_j$ ($j=1,2$) are Banach
spaces satisfying $\cX_1\subset\cX\subset\cX_2$
and $\cY_1\subset\cY\subset\cY_2$. If
$M\in\cL(\cX; \cY)$ can be extended as an
operator $\wt M\in\cL(\cX_2; \cY_2)$, then, to
simplify the notations, (formally) we also write
$M\in\cL(\cX_2; \cY_2)$. Similarly, if
$M|_{\cX_1}\in\cL(\cX_1; \cY_1)$, then, we write
$M\in \cL(\cX_1; \cY_1)$.

\ss

Further, we put
$$
\;\,\begin{array}{ll}\ds
\cL_{pd}\big(L^{p_1}_{\dbF}(0,T;L^{p_2}(\Om;\cX));\;L^{p_3}_{\dbF}(0,T;L^{p_4}(\Om;\cY))\big)
\\[1mm]
\triangleq \big\{\cL\in
\cL\big(L^{p_1}_{\dbF}(0,T;L^{p_2}(\Om;\cX));L^{p_3}_{\dbF}(0,T;L^{p_4}(\Om;\cY))\,\big|
\mbox{ for a.e. } \\[1mm]
\ns\ds\q (t,\omega)\in (0,T)\times\Omega, \mbox{
there exists } L(t,\omega)\in\cL (\cX;\cY)
\mbox{ verifying } \\[1mm]
\ns\ds\q\mbox{ that }  \big(\cL
f(\cd)\big)(t,\omega)=L (t,\omega)f(t,\omega),\;
\forall\, f(\cd)\in
L^{p_1}_{\dbF}(0,T;L^{p_2}(\Om;\cX))\}
\end{array}
$$
and
$$
\begin{array}{ll}\ds
\cL_{pd}\big(\cX;\;L^{p_3}_{\dbF}(0,T;$
$L^{p_4}(\Om;\cY))\big)
\\[1mm]
\triangleq \{\cL\in
\cL\big(\cX;L^{p_3}_{\dbF}(0,T;L^{p_4}(\Om;\cY))|
\mbox{for a.e. } (t,\omega)\in
(0,T)\times\Omega,\mbox{ there exists }\\[1mm]
\ns\ds\q L(t,\omega)\in\cL (\cX;\cY) \mbox{
verifying that }  \big(\cL x\big)(t,\omega)=L
(t,\omega)x, \; \forall\, x\!\in\! \cX\}.
\end{array}
$$
To simplify the notations, in what follows we
shall identify $\cL$ with $L(\cd,\cd)$.
Similarly, one can define the spaces
$\cL_{pd}\big(L^{p_2}(\Om;\cX);\;L^{p_3}_{\dbF}(0,T;L^{p_4}(\Om;\cY))\big)$
and $\cL_{pd}\big(L^{p_2}(\Om;\cX);$
$L^{p_4}(\Om;\cY)\big)$, etc.

\ss

In this paper, for any operator-valued
process/random variable $M$, we denote by $M^*$
its pointwise dual. For example, if $M\in
L^{r_1}_\dbF(0,T; L^{r_2}(\Om;$ $\cL(H)))$, then
$M^*\in L^{r_1}_\dbF(0,T; L^{r_2}(\Om;
\cL(H)))$, and $$| M|_{L^{r_1}_\dbF(0,T;
L^{r_2}(\Om; \cL(H)))}=| M^*|_{L^{r_1}_\dbF(0,T;
L^{r_2}(\Om; \cL(H)))}.$$

Put
\begin{equation}\label{10.10-eq30}
\begin{array}{lll}\ds
\Upsilon_p(\cX;\cY) \ds \triangleq
\big\{L(\cd,\cd)\in
\cL_{pd}(L^2_{\dbF}(\Om;L^\infty(0,T;\cX));L_{\dbF}^2(\Om;L^p(0,T;\cY))) |\\[1mm]
\ns \ds \qq\qq\qq\qq\qq\qq
|L(\cd,\cd)|_{\cL(\cX;\cY)}\in
L^\infty_\dbF(\Om;L^p(0,T))\big\}.
\end{array}
\end{equation}
We shall simply denote $\Upsilon_p(\cX;\cX)$ by
$\Upsilon_p(\cX)$.

\begin{remark}
For any $J(\cd,\cd)\in \Upsilon_p(\cX;\cY)$, one
may not have $J(\cd,\cd)\in
L^\infty_{\dbF}(\Om;L^p(0,T;\cL(\cX;\cY)))$.
Nevertheless, as we shall see later, in some
sense $\Upsilon_p(\cX;\cY)$ is a nice
``replacement" of the space
$L^\infty_{\dbF}(\Om;L^p(0,T;\cL(\cX;\cY)))$.
\end{remark}

Similar to Definition \ref{5.7-def1.1}, we
introduce the following notion:
\begin{definition}\label{5.7-def1}
An operator $\Th(\cd)\in \Upsilon_2(H;U)$ is
called an {\it optimal feedback  operator} for
{\bf Problem (SLQ)} if \vspace{0.1cm}
\begin{equation}\label{5.7-eq2}
\begin{array}{ll}\ds \cJ(s,\eta;\Th(\cd)\bar x(\cd))\leq
\cJ(s,\eta;u(\cd)),\\[2mm]
\ns\ds\q \forall\; (s,\eta)\in [0,T)\times
L^2_{\cF_s}(\Om;H),\;\; u(\cd)\in
L^2_\dbF(s,T;U),
\end{array}
\end{equation}
where $\bar x(\cd)=\bar x(\cd\,;s,\eta,
\Th(\cd)\bar x(\cd))$ solves the following
equation:\vspace{0.1cm}
\begin{equation}\label{5.2-eq1.1}
\left\{\begin{array}{ll}\ds d\bar
x(t)=\big[(A+A_1)\bar x(t) + B\Th \bar
x(t)\big]dt\\[2mm]
\ns\ds\qq\qq + \big(C \bar x(t)+D \Th \bar
x(t)\big)dW(t) &\mbox{ in
}(s,T],\\
\ns\ds \bar x(s)=\eta.
\end{array}
\right.
\end{equation}
\end{definition}
\begin{remark}
As a natural generalization of Definition
\ref{5.7-def1.1}, the space $\Upsilon_2(H;U)$ in
Definition \ref{5.7-def1} should be replaced by
$L^\infty_{\dbF}(\Om;L^2(0,T;\cL(H;U)))$.
However, we believe that $\Upsilon_2(H;U)$ is
more suitable in the infinite dimensional
setting. See Theorems \ref{5.7-th1.1} and
\ref{5.7-th1} for more details.
\end{remark}
\begin{remark}
In Definition \ref{5.7-def1}, the operator
$\Th(\cd)$ is required to be independent of
$(s,\eta)\in [0,T)\times L^2_{\cF_s}(\Om;H)$.
For a fixed pair $(s,\eta)\in[0,T)\times
L^2_{\cF_s}(\Om;H)$, the inequality
\eqref{5.7-eq2} implies that the control $\bar
u(\cd)\equiv \Th(\cd)\bar x(\cd)\in
L^2_\dbF(s,T;U)$
is optimal for {\bf Problem (SLQ)}. Therefore,
for {\bf Problem (SLQ)}, the existence of an
optimal feedback  operator on $[0,T]$ implies
the existence of optimal controls for any pair
$(s,\eta)\in [0,T)\times L^2_{\cF_s}(\Om;H)$.
\end{remark}

Let us introduce the following assumptions:

\vspace{0.32cm}

({\bf AS1}) {\it The coefficients satisfy that
$$
\begin{cases}\ds
A_1(\cd)\in
L^\infty_\dbF(\Om;L^1(0,T;\cL(H))),\\
\ns\ds C(\cd)\in
L^\infty_\dbF(\Om;L^2(0,T;\cL(H))),\\
\ns\ds B(\cd)\in L^\infty_\dbF(\Om;L^2(0,T;
\cL(U;H))),\\
\ns\ds D(\cd)\in L^\infty_\dbF(0,T;\cL(U;H)),\\
\ns\ds Q(\cd)\in L^\infty_\dbF(0,T; \cS(H)),\\
\ns\ds R(\cd) \in L^\infty_\dbF(0,T;\cS(U)),\\
\ns\ds G\in L^\infty_{\cF_T}(\Om;$ $\cS(H)),\q
G\geq 0, \q R>0,\q Q\geq0.
\end{cases}
$$
}

\vspace{0.32cm}

({\bf AS2}) {\it The eigenvectors
$\{e_j\}_{j=1}^\infty$ of $A$ such that
$|e_j|_H=1$ for all $j\in\dbN$ constitute an
orthonormal basis of $H$.}

\vspace{0.32cm}

By ({\bf AS1}), it is easy to see that, for any
$(s,\eta)\in[0,T]\times L^2_{\cF_s}(\Om;H)$,
there exists a unique optimal control for {\bf
Problem (SLQ)}.

\vspace{0.1cm}

Let $\{\mu_j\}_{j=1}^\infty$ (corresponding to
$\{e_j\}_{j=1}^\infty$) be the eigenvalues of
$A$. Let $\{\lambda_j\}_{j=1}^\infty$ be an
arbitrarily given real number sequence
satisfying that $\lambda_j>0$ for all $j\in\dbN$
and $\sum_{j=1}^\infty \lambda_j^2<\infty$.
Define a norm $|\cd|_V$ on $H$ as follows:
$$
|h|_V = \sqrt{\sum_{j=1}^\infty \lambda_j^2
|e_j|_{D(A)}^{-2}h_j^2},\qq \forall\,
h=\sum_{j=1}^\infty h_j e_j \in H.
$$
Denote by $V$ the completion of $H$ with respect
to the norm $|\cd|_V$. Clearly, $V$ is a Hilbert
space, $V\subset H$ and
$\{\lambda_j^{-1}|e_j|_{D(A)}e_j\}_{j=1}^\infty$
is an orthonormal basis of $V$. Denote by
$\cV_H$ the set of all such kind of Hilbert
spaces $V$.

\vspace{0.1cm}

Denote by $\cL_2(H;V)$ the set of all
Hilbert-Schmidt operators from $H$ to $V$. It is
well known that $\cL_2(H;V)$ is a Hilbert space
itself. Denote by $V'$ the dual space of $V$
with respect to the pivot space $H\equiv H'$.

\vspace{0.1cm}

We also need the following technical conditions:

\vspace{0.32cm}

({\bf AS3}) {\it There is a $V\in\cV_H$ such
that $A_1,C,Q\in L^\infty_\dbF(0,T;\cL(V))$,
$G\in L^\infty_{\cF_T}(\Om;$ $\cL(V))$ and
$A_1,C,Q\in L^\infty_\dbF(0,T;\cL(V'))$, $G\in
L^\infty_{\cF_T}(\Om;\cL(V'))$.}

\vspace{0.32cm}

({\bf AS4}) {\it Let $\{\f_j\}_{j=1}^\infty$ be
an orthonormal basis of\; $U$. There is a $\wt
U\subset U$ such that $\wt U$ is dense in $U$,
$\{\f_j\}_{j=1}^\infty\subset \wt U$, $R\in
L^\infty_\dbF(0,T;\cL(\wt U))$ and $B,D\in
L^\infty_\dbF(0,T;\cL(\wt U;V'))$, where $V$ is
given in Assumption ({\bf AS3}).}

\vspace{0.2cm}

Consider the following two (forward) SEEs:
\begin{equation}\label{op-fsystem1}
\left\{
\begin{array}{ll}
\ds dx_1 = \big[(A+A_1) x_1 + u_1\big]d\tau + \big(C x_1 + v_1\big)dW(\tau) &\mbox{ in } (t,T],\\[1mm]
\ns\ds x_1(t)=\xi_1
\end{array}
\right.
\end{equation}
and
\begin{equation}\label{op-fsystem2}
\left\{
\begin{array}{ll}
\ds dx_2 =\big[ (A+A_1) x_2 + u_2\big]d\tau +\big( C x_2 + v_2 \big)dW(\tau) &\mbox{ in } (t,T],\\[1mm]
\ns\ds x_2(t)=\xi_2.
\end{array}
\right.
\end{equation}
Here $t\in [0,T)$, $\xi_1,\xi_2$ are suitable
random variables and $u_1,u_2,v_1,v_2$ are
suitable stochastic processes.

\vspace{0.1cm}

Put\vspace{0.1cm}
$$
\ba{ll} \ds C_{\dbF,w}([0,T];
L^{\infty}(\Om;\cL(H)))\\
\ns \ds\triangleq\Big\{P\in \Upsilon_2(H)\,
\big|\,P(t,\omega)\in \cS(H),\;\mbox{a.e. }
(t,\omega)\!\in\! [0,T]\!\times\!\Omega,\;
|P(\cd)|_{\cL(H)}\!\in\!
L^\infty_\dbF(0,T),\\
\ns \ds\qq\qq \qq\q\,\; \mbox{and
}P(\cd)\zeta\in C_\dbF([0,T];L^\infty(\Om;H)),\
\forall\;\zeta\in H\Big\} \ea
$$
and
$$
\begin{array}{ll}\ds
L^2_{\dbF,w}(0,T;\cL(H))
\\
\ns\ds\triangleq\Big\{\Lambda\in
\cL_{pd}(L^2_{\dbF}(\Om;L^\infty(0,T;H));L_{\dbF}^2(0,T;H))\bigcap
L^2_{\dbF}(0,T;\cL_2(V';V)) \;
| \\
\ns\ds\qq
D^*\Lambda\in\Upsilon_2(H;U),\;\Lambda(t,\omega)\in
\cS(H),\;\mbox{a.e. } (t,\omega)\in
[0,T]\times\Omega \Big\}.
\end{array}
$$

Now, we are in a position to introduce the
notion of {\it transposition solution} to
\eqref{5.5-eq6}:
\begin{definition}\label{4.8-def2}
We call $\big(P(\cd),\Lambda(\cd)\big)\in
C_{\dbF,w}([0,T]; L^{\infty}(\Om;\cL(H))) \times
L^2_{\dbF,w}(0,T;$ $\cL(H))$ a transposition
solution to \eqref{5.5-eq6} if the following
three conditions hold:

\ms

{\rm 1)}  $K(t,\om)\big(\equiv R(t,\om) +
D(t,\om)^*P(t,\om)D(t,\om)\big)\in \cS(U)$,
$K(t,\om) > 0$ and its left inverse
$K(t,\om)^{-1}$ is a densely defined closed
operator for a.e. $(t,\om)\in [0,T]\times\Om$;

\ms

{\rm 2)}  For any $t\in [0,T]$, $\xi_1,\xi_2\in
L^4_{\cF_t}(\Om;V')$, $u_1(\cd), u_2(\cd)\in
L^4_\dbF(\Om;L^2(t,T;V'))$ and $v_1(\cd),
v_2(\cd)$ $\in L^4_\dbF(\Om;L^2(t,T;V'))$, it
holds that
\begin{eqnarray}\label{6.18-eq1}
&&  \mE\langle Gx_{1}(T),x_{2}(T)\rangle_{H}
+\mE \int_t^T \big\langle Q(\tau) x_{1}(\tau),
x_{2}(\tau) \big\rangle_{H}d\tau\nonumber\\
&&\q- \mE \int_t^T \big\langle K(\tau)^{-1}
L(\tau) x_{1}(\tau), L(\tau)x_{2}(\tau)
\big\rangle_{H}d\tau\nonumber
\\
&& = \mE\big\langle P(t) \xi_{1},\xi_{2}
\big\rangle_{H} + \mE \int_t^T \big\langle
P(\tau)u_{1}(\tau),
x_{2}(\tau)\big\rangle_{H}d\tau  \\
&&  \q + \mE \int_t^T \big\langle
P(\tau)x_{1}(\tau),
u_{2}(\tau)\big\rangle_{H}d\tau + \mE
\int_t^T\big\langle P(\tau)C(\tau)x_{1}(\tau),
v_{2}(\tau)\big\rangle_{H}d\tau \nonumber\\
&&\q + \mE
\int_t^T \big\langle  P(\tau)v_{1}(\tau), C(\tau)x_{2}(\tau)+v_{2}(\tau)\big\rangle_{H}d\tau\nonumber\\
&& \q + \mE \int_t^T \big\langle v_{1}(\tau),
\Lambda(\tau)x_2(\tau)\big\rangle_{ V', V}d\tau+
\mE \int_t^T \big\langle \Lambda(\tau)x_1(\tau),
v_{2}(\tau) \big\rangle_{ V,V'}d\tau,\nonumber
\end{eqnarray}
where $x_1(\cd)$ and $x_2(\cd)$ solve
\eqref{op-fsystem1} and \eqref{op-fsystem2},
respectively\;\footnote{By Corollary \ref{lm15}
(in Chapter \ref{pre}), one has $x_1(\cd),
x_2(\cd)\in L^4_\dbF(\Om;C([0,T];V'))$.}; and

\ms

{\rm 3)}  For any $t\in [0,T]$, $\xi_1,\xi_2\in
L^2_{\cF_t}(\Om;H)$, $u_1(\cd), u_2(\cd)\in
L^2_\dbF(t,T;H)$  and $v_1(\cd),$ $v_2(\cd)\in
L^2_\dbF(t,T;U)$, it holds that
\begin{eqnarray}\label{10.10-eq10}
&& \mE\langle Gx_{1}(T),x_{2}(T)\rangle_{H} +\mE
\int_t^T \big\langle Q(\tau) x_{1}(\tau),
x_{2}(\tau) \big\rangle_{H}d\tau\nonumber\\
&&\q - \mE \int_t^T \big\langle K(\tau)^{-1}
L(\tau) x_{1}(\tau), L(\tau)x_{2}(\tau)
\big\rangle_{H}d\tau\nonumber
\\
&&  = \mE\big\langle P(t) \xi_{1},\xi_{2}
\big\rangle_{H} + \mE \int_t^T \big\langle
P(\tau)u_{1}(\tau),
x_{2}(\tau)\big\rangle_{H}d\tau  \nonumber\\
&&  \q + \mE \int_t^T \big\langle
P(\tau)x_{1}(\tau),
u_{2}(\tau)\big\rangle_{H}d\tau + \mE \int_t^T
\big\langle P(\tau)C(\tau)x_{1}(\tau), D(\tau)
v_{2}(\tau)\big\rangle_{H}d\tau\\
&&\q + \mE
\int_t^T \big\langle  P(\tau)D(\tau) v_{1}(\tau), C(\tau)x_{2}(\tau)+D (\tau)v_{2}(\tau)\big\rangle_{H}d\tau\nonumber\\
&&\q + \mE \int_t^T \big\langle v_{1}(\tau),
D(\tau)^*\Lambda(\tau)x_2(\tau)\big\rangle_Ud\tau+
\mE \int_t^T \big\langle
D(\tau)^*\Lambda(\tau)x_1(\tau), v_{2}(\tau)
\big\rangle_Ud\tau.\nonumber
\end{eqnarray}
Here, $x_1(\cd)$ and $x_2(\cd)$ solve
\eqref{op-fsystem1} and \eqref{op-fsystem2} with
$v_1$ and $v_2$ replaced by $Dv_1$ and $Dv_2$,
respectively.
\end{definition}
\ms

The main results of this paper, which reveal the
relationship between the existence of optimal
feedback operator  for {\bf Problem (SLQ)} and
the well-posedness of \eqref{5.5-eq6} in the
sense of transposition solution, are stated as
follow:

\begin{theorem}\label{5.7-th1.1}
Let (AS1) hold. If the Riccati equation
\eqref{5.5-eq6} admits a transposition solution
$\big(P(\cd),\Lambda(\cd)\big)\in
C_{\dbF,w}([0,T]; L^{\infty}(\Om;\cL(H))) \times
L^2_{\dbF,w}(0,T;\cL(H))$ such that
\begin{equation}\label{5.7-eq5}
\begin{array}{ll}\ds
K(\cd)^{-1}\big[B(\cd)^* P(\cd) +D(\cd)^*
P(\cd)C(\cd) + D(\cd)^*\Lambda(\cd)\big] \in
\Upsilon_2(H;U)\cap \Upsilon_2(V';\wt U),
\end{array}
\end{equation}
then {\bf Problem (SLQ)} admits an optimal
feedback operator $\Th(\cd)\in
\Upsilon_2(H;U)\cap \Upsilon_2(V';\wt U)$.
Furthermore, the optimal feedback operator
$\Th(\cd)$ is given by
\begin{equation}\label{5.10}
\begin{array}{ll}
\ns\ds\Th(\cd)=-K(\cd)^{-1}[B(\cd)^* P(\cd)
+D(\cd)^* P(\cd)C(\cd) + D(\cd)^*\Lambda(\cd)],
\end{array}
\end{equation}
and
\begin{equation}\label{Value}
\inf_{u\in
L^2_\dbF(s,T;U)}\cJ(s,\eta;u)=\frac{1}{2}\,\dbE\langle
P(s)\eta,\eta\rangle_H.
\end{equation}
\end{theorem}
\ms
\begin{theorem}\label{5.7-th1}
Let (AS1)--(AS4) hold and $A$ generate a
$C_0$-group on $H$. If {\bf Problem (SLQ)}
admits an optimal feedback operator $\Th(\cd)\in
\Upsilon_2(H;U)\cap \Upsilon_2(V';\wt U)$, then
the Riccati equation \eqref{5.5-eq6} admits a
unique transposition solution
$\big(P(\cd),\Lambda(\cd)\big)$ $\in
C_{\dbF,w}([0,T]; L^{\infty}(\Om;\cL(H))) \times
L^2_{\dbF,w}(0,T;\cL(H))$ such that
\eqref{5.7-eq5} holds and  the optimal feedback
operator $\Th(\cd)$ is given by \eqref{5.10}.
Furthermore, \eqref{Value} holds.
\end{theorem}
\begin{remark}
Theorem \ref{5.7-th1.1} concludes that if the
Riccati equation \eqref{5.5-eq6} admits a
transposition solution satisfying
\eqref{5.7-eq5}, then {\bf Problem (SLQ)} admits
an optimal feedback operator. Theorem
\ref{5.7-th1} says that the converse of Theorem
\ref{5.7-th1.1} is also true under some
additional conditions.
\end{remark}

More remarks are in order.

\begin{remark}
In Theorem \ref{5.7-th1} (and recalling
Definition \ref{4.8-def2}), we only conclude
that $K(t,\om)$ has left inverse for a.e.
$(t,\om)\in (0,T)\times\Om$, and therefore
$K(t,\om)^{-1}$ may be unbounded. Nevertheless,
these results cannot be improved. Let us show
this by the following example.

\vspace{0.1cm}

Let $\cO\subset\dbR^k$ (for some $k\in\dbN$) be
a  bounded domain with a smooth boundary
$\pa\cO$.  Let $$H= H_0^1(\cO)\times L^2(\cO),
\q  U=L^2(\cO),  \q A=\left(
     \begin{array}{cc}
       0 & I \\
       \D & 0 \\
     \end{array}
   \right),
$$ where  $\D$ is the Laplacian on $\cO$ with
the usual homogeneous Dirichlet boundary
condition. Let \vspace{0.1cm} $$B=\left(
     \begin{array}{c}
       0 \\
       I \\
     \end{array}
   \right), \q
 C=\left(
     \begin{array}{c}
       I \\
       0 \\
     \end{array}
   \right), \q  D=0, \q Q=0, \q R=(-\D)^{-1}, \q
   G=0.$$
Then \eqref{5.2-eq1} is specialized
as\vspace{0.1cm}
\begin{equation}\label{7.8-eq1}
\left\{ \2n\begin{array}{ll}\ds dX =\big(A X  +
u \big)dt +  X dW(t) &\mbox{ in
}(s,T], \\
\ns\ds X(s)=\eta.
\end{array}
\right.
\end{equation}
The cost functional reads
\begin{equation}\label{7.8-eq2}
\cJ(s,\eta;u(\cd)) =\frac{1}{2}\mE \int_s^T
\big\langle
(-\D)^{-1}u(t),u(t)\big\rangle_{L^2(\cO)} dt.
\end{equation}
Clearly, for any $(s,\eta)\in [0,T)\times
L^2_{\cF_s}(\Om;L^2(\cO))$, there is a unique
optimal control $u\equiv 0$. For the present
case, it is easy to check that
$(P(\cd),\Lambda(\cd))=(0,0)$ is the unique
 transposition solution to \eqref{5.5-eq6}. However,
$K=(-\D)^{-1}$ is not surjective and $K^{-1}$ is
unbounded.
\end{remark}

\br
In Theorem \ref{5.7-th1}, we assume that $A$
generates a $C_0$-group on $H$. This assumption
is used to guarantee the well-posednss of the
(\ref{5.26-eq4}) in Chapter \ref{proof-n}. We
believe that  it should be a technical condition
but, so far we do not know how to drop it.
\er

\begin{remark}
In this paper, we have assumed that $R(t)>0$ for
a.e. $t\in [0,T]$. Following the proof of
\cite[Theorem 2.3]{Lu}, one can show that this
assumption can be relaxed to the condition that
the map $u$ to $\cJ(0,0;u)$ is uniformly convex.
However, as we have said, to present the key
idea in a simple way, we will not do this
technical generalization. On the other hand, it
is a really challengeable problem to drop the
positivity condition on $R$ or the uniform
convexity condition on the map $u$ to
$\cJ(0,0;u)$ in \cite[Theorem 2.3]{Lu}. This was
done when all the coefficients of {\bf Problem
(SLQ)} are deterministic (e.g. \cite{Lu}), where
the Riccati equation \eqref{5.5-eq6} becomes a
deterministic operator-valued evolution
equation. However, the method in \cite{Lu}
cannot be used to handle the stochastic problem.
\end{remark}

\br It is easy to see that, under Assumption
(AS3), if $\xi_j\in L^2_{\cF_t}(\Om;V')$ and
$u_j(\cd), v_j(\cd) \in L^2_\dbF(t,T;V')$, then
the solutions $x_j$ ($j=1,2$) to
\eqref{op-fsystem1}--\eqref{op-fsystem2} belong
to $ L^2_\dbF(\Om;$ $C([t,T];V'))$. This plays a
key role in Step 5 in the proof of Theorem
\ref{5.7-th1}. We believe that this assumption
can be dropped. However, we do not know how to
do it at this moment. \er

\br
In Theorem \ref{5.7-th1}, the most natural
choice of optimal feedback operator set should
be $ \Upsilon_2(H;U)$ rather than $
\Upsilon_2(H;U)\cap \Upsilon_2(V';\wt U)$.
Nevertheless, at this moment, in the proof of
Theorem \ref{5.7-th1} (see Chapter
\ref{proof-n}), we do need to suppose that
$\Th(\cd)\in \Upsilon_2(H;U)\cap
\Upsilon_2(V';\wt U)$. \er

\br
 It would be quite interesting to extend the main result in this paper
to linear quadratic stochastic differential
games in infinite dimensions whereas this
remains to be done.
\er

\section{Some preliminary results}\label{pre}


\ms

\subsection{Well-posedness for some SEEs and
BSEEs}

\ms

In this section, we present well-posedness
results for some SEEs and BSEEs, which will be
useful in the sequel.

\ms

First, for any $s\in[0,T)$,  consider the
following SEE:
\begin{equation}\label{6.20-eq1}
\left\{
\begin{array}{ll}\ds
dx = [(A+\cA) x  + f]dt + (\cB x+g)dW(t) &\mbox{ in }(s,T],\\
\ns\ds x(s)=\eta.
\end{array}
\right.
\end{equation}
Here $\cA \in \Upsilon_1(H)$, $\cB\in
\Upsilon_2(H)$, $\eta\in L^2_{\cF_s}(\Om;H)$,
$f\in L^2_\dbF(\Om;L^1(s,T;H))$  and $g\in
L^2_\dbF(s,T;H)$.

\vspace{0.1cm}

We have the following result.

\begin{lemma}\label{lm2}
The equation \eqref{6.20-eq1} admits a unique
mild solution $x(\cd)\in L^2_\dbF(\Om;$
$C([s,T];H))$ satisfying
\begin{equation}\label{lm2-eq1}
| x(\cd)|_{L^2_\dbF(\Om;C([s,T];H))}\leq
\cC\big(|\eta|_{L^2_{\cF_s}(\Om;H)} +
|f|_{L^2_\dbF(\Om;L^1(0,T;H))} +
|g|_{L^2_\dbF(s,T;H)}\big).
\end{equation}
\end{lemma}

{\it Proof}\,: We borrow some idea from
\cite[Chapter V, Section 3]{Protter}. Without
loss of generality, let us assume that $s=0$.
Write
$$N=\lceil\frac{1}{\e}\big(||\cA|_{\cL(H)}|^2_{L^\infty_\dbF(\Om;L^1(0,T))}
+
||\cB|_{\cL(H)}|^2_{L^\infty_\dbF(\Om;L^2(0,T))}\big)\rceil
+1,
 $$
where $\e>0$ is a constant to be determined
later. Define a sequence of stopping times
$\{\tau_{j,\e}\}_{j=1}^N$ as follows:
$$
\left\{\ba{ll}
\ds\tau_{1,\e}(\om)=\inf\Big\{t\in
[0,T]\,\Big|\,\(\int_0^t|\cA(s,\om)|_{\cL(H)}ds\)^2+\int_0^t|\cB(s,\om)|^2_{\cL(H)}ds=\e\Big\},\\
\ns\ds \tau_{k,\e}(\om)= \inf\Big\{t\in
[0,T]\;\Big|\;\(\int_{\tau_{k-1,\e}}^t|\cA(s,\om)|_{\cL(H)}ds\)^{2}\\
\ns\ds\hspace{4cm}+\int_{\tau_{k-1,\e}}^t|\cB(s,\om)|^2_{\cL(H)}ds=\e\Big\},\q
k= 2,\cdots,N. \ea\right.
$$
Here, we agree that $\inf\emptyset=T$.

\vspace{0.1cm}

Consider the following  SEE:
\begin{equation}\label{10.29-eq9}
\left\{
\begin{array}{ll}\ds
dx = \big(A x + \tilde f\big)dt + \tilde gdW(t) &\mbox{ in }(0,T],\\
\ns\ds x(0)=\eta,
\end{array}
\right.
\end{equation}
where $\tilde f\in L^2_\dbF(\Om;L^1(0,T;H))$ and
$\tilde g\in L^2_\dbF(0,T;H)$. Clearly,
\eqref{10.29-eq9} admits a unique mild solution
$x\in L^2_\dbF(\Om;C([0,T];H))$.  Define a map
$$\cJ:\;L^2_\dbF(\Om;C([0,\tau_{1,\e}];H))\to
L^2_\dbF(\Om;C([0,\tau_{1,\e}];H))$$ as follows:
$$
L^2_\dbF(\Om;C([0,\tau_{1,\e}];H))\ni\tilde
x\mapsto x=\cJ(\tilde x),
$$
where $x$ is the solution to \eqref{10.29-eq9}
with $\tilde f$ and $\tilde g$ replaced by $\cA
\tilde x + f$ and $\cB \tilde x + g$, \linebreak
respectively. We claim that $\cJ$ is
contractive. Indeed, for any $\tilde x_1,\tilde
x_2\in L^2_\dbF(\Om;$ $C([0,\tau_{1,\e}];H))$,
\begin{eqnarray}\label{10.29-eq11}
&& |\cJ(\tilde x_1) - \cJ(\tilde
x_2)|_{L^2_\dbF(\Om;C([0,\tau_{1,\e}];H))}^2\nonumber\\
&& \leq \mE\(\sup_{t\in
[0,\tau_{1,\e}]}\Big|\int_0^t e^{A(t-r)}\cA
(\tilde x_1-\tilde x_2)dr + \int_0^t
e^{A(t-r)}\cB (\tilde x_1-\tilde
x_2)dW(r)\Big|_H\)^2\nonumber\\
&& \leq 2\mE\(\sup_{t\in [0,\tau_{1,\e}]}\Big|
\int_0^t  e^{A(t-r)}\cA (\tilde x_1 - \tilde
x_2)dr\Big|_H\)^2\\
&&\q + 2\mE\sup_{t\in [0,\tau_{1,\e}]}\(\Big|
\int_0^t\! e^{A(t-r)}\cB (\tilde
x_1\!-\!\tilde x_2)dW(r)\Big|_H\)^2 \nonumber\\
&& \leq 2\(\sup_{t\in
[0,T]}|e^{At}|^2_{\cL(H)}\)\big||\cA|_{\cL(H)}\big|^2_{L^\infty_\dbF(\Om;L^1(0,\tau_{1,\e}))}|\tilde
x_1-\tilde x_2|^2_{L^2_\dbF(\Om;C([0,\tau_{1,\e}];H))}\nonumber\\
&& \q + 2\(\sup_{t\in
[0,T]}|e^{At}|^2_{\cL(H)}\)\big||\cB|_{\cL(H)}\big|^2_{L^\infty_\dbF(\Om;L^2(0,\tau_{1,\e}))}|\tilde
x_1-\tilde
x_2|^2_{L^2_\dbF(\Om;C([0,\tau_{1,\e}];H))}.\nonumber
\end{eqnarray}
Let us choose
$$
 \e = \frac{1}{16\sup_{t\in
[0,T]}|e^{At}|^2_{\cL(H)}}.
$$
Then, from \eqref{10.29-eq11}, we find that
\begin{equation}\label{10.29-eq12}
|\cJ(\tilde x_1) - \cJ(\tilde
x_1)|_{L^2_\dbF(\Om;C([0,\tau_{1,\e}];H))}^2
\leq \frac{1}{4}|\tilde x_1-\tilde
x_2|^2_{L^2_\dbF(\Om;C([0,\tau_{1,\e}];H))}.
\end{equation}
Hence, $\cJ$ is contractive, and it has a unique
fixed point $x\in
L^2_\dbF(\Om;C([0,\tau_{1,\e}];H))$, which
solves \eqref{6.20-eq1} (with $s=0$) in
$[0,\tau_{1,\e}]$ (in the sense of mild
solution). Inductively, we conclude that
\eqref{6.20-eq1} admits a mild solution $x$ in
$[\tau_{k-1,\e},\tau_{k,\e}]$ for $k=2,\cds,N$.
Furthermore,
\begin{eqnarray}\label{10.29-eq13}
&&\3n\3n
|x|_{L^2_\dbF(\Om;C([0,\tau_{1,\e}];H))}^2\nonumber\\
&&\3n\3n \leq \mE\sup_{t\in
[0,\tau_{1,\e}]}\(\Big|e^{At}\eta + \int_0^t
e^{A(t-r)}(\cA  x+f) dr +
\int_0^t e^{A(t-r)}(\cB x +g)dW(r) \Big|_H\)^2\nonumber\\
&&\3n\3n \leq 4\mE\sup_{t\in
[0,\tau_{1,\e}]}\(\Big|\int_0^t e^{A(t-r)}\cA x
dr\Big|_H\)^2 + 4\mE\sup_{t\in
[0,\tau_{1,\e}]}\(\Big| \int_0^t e^{A(t-r)}\cB
x dW(r)\Big|_H\)^2\nonumber \\
&&\3n\3n \q + 4 \mE\sup_{t\in [0,T]}
\(\Big|e^{At}\eta + \int_0^t e^{A(t-r)}f dr +
\int_0^t e^{A(t-r)}g dW (r)\Big|_H\)^2\\
&& \3n\3n\leq 4\(\sup_{t\in
[0,T]}|e^{At}|^2_{\cL(H)}\)\big||\cA|_{\cL(H)}\big|^2_{L^\infty_\dbF(\Om;L^1(0,\tau_{1,\e}))}|
x|^2_{L^2_\dbF(\Om;C([0,\tau_{1,\e}];H))}\nonumber\\
&&\3n\3n\q + 4\(\sup_{t\in
[0,T]}|e^{At}|^2_{\cL(H)}\)\big||\cB|_{\cL(H)}\big|^2_{L^\infty_\dbF(\Om;L^2(0,\tau_{1,\e}))}|
x|^2_{L^2_\dbF(\Om;C([0,\tau_{1,\e}];H))}\nonumber\\
&& \3n\3n\q + \cC\big(|\eta|_{H}^2 +
|f|_{L^2_\dbF(\Om;L^1(0,T;H))}^2 +
|g|_{L^2_\dbF(0,T;H)}^2\big).\nonumber
\end{eqnarray}
This, together with the choice of $\tau_{1,\e}$,
implies that
\begin{equation}\label{10.29-eq14}
|x|_{L^2_\dbF(\Om;C([0,\tau_{1,\e}];H))}^2\leq
\cC\big(|\eta|_{H}^2 +
|f|_{L^2_\dbF(\Om;L^1(0,T;H))}^2 +
|g|_{L^2_\dbF(0,T;H)}^2\big).
\end{equation}
Repeating the above argument, we obtain
\eqref{lm2-eq1}. The uniqueness of the solution
is obvious.

\br
From the proof of Lemma \ref{lm2}, it is easy to
see that our assumptions $\cA \in \Upsilon_1(H)$
and $\cB\in \Upsilon_2(H)$  are sharp for the
well-posedness of \eqref{6.20-eq1}. Hence, our
choice of the optimal feedback  operator
$\Th(\cd)\in \Upsilon_2(H;U)$ in Definition
\ref{5.7-def1} is also sharp. On the other hand,
it is easy to see that Lemma \ref{fflm2} is a
special case of Lemma \ref{lm2}.
\er

Next, we consider the following BSEE:
\begin{equation}\label{6.20-eq10}
\left\{
\begin{array}{ll}\ds
dy = -[(A+A_1)^*y+ \cD z + h]dt + zdW(t) &\mbox{ in }[0,T),\\
\ns\ds y(T)=\xi.
\end{array}
\right.
\end{equation}
Here $\xi\in L^2_{\cF_T}(\Om;H)$, $\cD\in
L^\infty_\dbF(0,T;\cL(H))$ and $h\in
L^2_\dbF(0,T;H)$. Let us recall the following
known result (e.g. \cite{Mahmudov1}).

\begin{lemma}\label{lm3}
The equation \eqref{6.20-eq10} admits a unique
mild solution $(y(\cd),z(\cd))\in
L^2_\dbF(\Om;C([0,T];$ $H))\times
L^2_\dbF(0,T;H)$, and
$$
|(y(\cd),z(\cd))|_{L^2_\dbF(\Om;C([0,T];H))\times
L^2_\dbF(0,T;H)}\leq
\cC\big(|\xi|_{L^2_{\cF_T}(\Om;H)}+|h|_{L^2_\dbF(0,T;H)}\big).
$$
\end{lemma}

Also, let us recall the following
Pontryagin-type maximum principle (\cite[Theorem
5.2]{LZ}).

\begin{lemma}\label{5.4-th1}
Let {\bf Problem (SLQ)} be solvable at
$(s,\eta)\in[0,T)\times L^2_{\cF_s}(\Om;H)$ with
$(\bar x(\cd),\bar u(\cd))$ $\in
L^2_\dbF(\Om;C([s,T];H))\times L^2_\dbF(s,T;U)$
being an optimal pair. Then there exists a pair
$ (\bar y(\cd), \bar z(\cd)) \in
L^2_\dbF(\Om;C([s,T];H))\times L^2_\dbF(s,T;H) $
satisfying the following BSEE:
$$
\left\{
\begin{array}{ll}\ds
d\bar y(t)=-\big[(A+A_1)^*\bar y(t)+C^*\bar
z(t)+Q\bar x(t)\big]dt +\bar z(t)dW(t)
&\mbox{\rm in }[s,T),\cr
\ns\ds\bar y(T)= G\bar x(T),
\end{array}
\right.
$$
and
$$
R\bar u+B^*\bar y+D^*\bar z =0,\q \mbox{ for
a.e. }  (t,\om)\in [s,T]\times\Om.
$$
\end{lemma}

As an immediate consequence of Lemmas \ref{lm2}
and \ref{5.4-th1}, we have the following result.
\begin{corollary}\label{5.7-prop1}
Let $\Th(\cd)$ be an optimal feedback operator
for {\bf Problem (SLQ)}. Then, for any
$(s,\eta)\in[0,T)\times L^2_{\cF_s}(\Om;H)$, the
following forward-backward stochastic evolution
equation:
$$
\left\{
\begin{array}{ll}
\ns\ds d\bar x(t) =[(A+A_1)+B\Th]\bar x(t) dt+
(C+D\Th)\bar x(t) dW(t)
&\mbox{\rm in } (s,T],\\[1mm]
\ns\ds d\bar y(t)=-\big[(A+A_1)^* \bar y(t)+C^*
\bar z(t)+ Q\bar x(t)\big]dt+\bar
z(t)dW(t)\quad&\mbox{\rm in }
[s,T),\\[1mm]
\ns\ds \bar x(s)=\eta,\q\bar y(T)=G\bar x(T)
\end{array}
\right.
$$
admits a unique mild solution
$$ (\bar x(\cd),\bar y(\cd),\bar z(\cd))\!\in\!
L^2_\dbF(\Om;C([s, T];\!H))\times
L^2_\dbF(\Om;C([s,T];\!H))\times
L^2_\dbF(s,T;\!H), $$
and
$$
R\Th \bar x+B^* \bar y +D^* \bar z  =0,\q \mbox{
for a.e. }  (t,\om)\in [s,T]\times\Om.
$$
\end{corollary}

\ms

\subsection{$\cL(H;V)$-valued SEEs and BSEEs}

\ms

For any $V\in\cV_H$, from the definition of $V$,
it is easy to see that
$\{e_j\}_{j=1}^\infty\subset V'$ is an
orthogonal basis of $V'$ and the norm on $V'$ is
given as follows:
$$
|\xi|_{V'} = \sqrt{\sum_{j=1}^\infty
|\xi_j|^2|e_j|_{D(A)}^2\lambda_j^{-2}},\qq
\forall\, \xi \in V',
$$
where $\xi_j=\langle \xi,e_j\rangle_H$.
Furthermore,
$\{\lambda_j|e_j|_{D(A)}^{-1}e_j\}_{j=1}^\infty$
is an orthonormal basis of $V'$.

\begin{lemma}\label{lm14}
Let $V\in\cV_H$ and (AS2) hold. If
$\{e^{At}\}_{t\in\dbR}$ is a $C_0$-group on $H$,
then it is a $C_0$-group on $V'$, and it can be
uniquely extended to a $C_0$-group (also denoted
by itself) on $V$.
\end{lemma}

{\it Proof}\,:  We only prove that
$\{e^{At}\}_{t\geq 0}$ is a $C_0$-group on $V'$.
The proof for the other conclusion is similar.

\vspace{0.1cm}

Let $\ds \xi=\sum_{j=1}^\infty \xi_j
\lambda_j|e_j|_{D(A)}^{-1}e_j\in V'$ with
$\{\xi_j\}_{j=1}^\infty\in\ell^2$. Then,
$\ds\tilde \xi=\sum_{j=1}^\infty \xi_j e_j\in H$
and $|\tilde \xi|_H = |\xi|_{V'}$. Clearly,
$$
e^{At}\xi=\sum_{j=1}^\infty \xi_j \lambda_j
|e_j|_{D(A)}^{-1} e^{\mu_j t}e_j.
 $$

For any $t_1,t_2\in \dbR$,
\begin{equation}\label{12.15-eq1}
\begin{array}{ll}\ds
e^{At_2}e^{At_1}\xi\3n&\ds=\sum_{j=1}^\infty
\xi_j \lambda_j |e_j|_{D(A)}^{-1}
e^{At_2}e^{At_1}e_j\\
\ns&\ds= \sum_{j=1}^\infty \xi_j \lambda_j
|e_j|_{D(A)}^{-1} e^{\mu_j(t_1+t_2)}e_j\\
\ns&\ds = e^{A(t_2+t_1)}\xi.
\end{array}
\end{equation}
This indicates that $\{e^{At}\}_{t\geq 0}$ is a
group on $V'$.

\vspace{0.1cm}

For any $t_2>t_1>0$,
\begin{equation*}\label{12.15-eq2}
\begin{array}{ll}\ds
\big|\big(e^{At_2} - e^{At_2}\big)\xi
\big|_{V'}\3n&\ds=\Big|\sum_{j=1}^\infty \xi_j
\lambda_j |e_j|_{D(A)}^{-1} \big(e^{At_2} -
e^{At_1}\big)e_j\Big|_{V'}\\
\ns&\ds = \Big|\sum_{j=1}^\infty \xi_j \lambda_j
|e_j|_{D(A)}^{-1}
\big(e^{\mu_jt_2}\!-\!e^{\mu_jt_1}\big)e_j\Big|_{V'}\\
\ns&\ds = \[\sum_{j=1}^\infty |\xi_j|^2
\big(e^{\mu_jt_2}-e^{\mu_jt_1}\big)^2\]^{\frac{1}{2}}\\
\ns&\ds = \big|\big(e^{At_2} -
e^{At_2}\big)\tilde \xi\big|_{H}.
\end{array}
\end{equation*}
This, together with that $\{e^{At}\}_{t\geq 0}$
is strongly continuous on $H$, implies that
$\{e^{At}\}_{t\geq 0}$ is strongly continuous on
$V'$.

\ms

With the aid of Lemma \ref{lm14}, as an
immediate consequence of the standard well-posed
result for SEEs, we have the following result.
\begin{corollary}\label{lm15}
Let (AS2)--(AS3) hold. Then, for any $\xi_1\in
L^4_{\cF_t}(\Om;V')$ and $u_1(\cd), v_1(\cd) \in
L^4_\dbF(\Om;$ $L^2(t,T;V'))$, the mild solution
$x_1(\cd)$ to \eqref{op-fsystem1} belongs to
$L^4_\dbF(\Om;C([0,T];V'))$.
\end{corollary}

For any $n\in\dbN$, denote by $\G_n$ the
projection operator from $H$  to
$H_n\triangleq\span_{1\leq j\leq n}\{e_j\}$. Let
\bel{20161214e1}
\begin{cases} \ds A_n = \G_n
A \G_n,
\qq \cA_n = \G_n \cA \G_n, \\
\ns\ds \cB_n = \G_n \cB \G_n,\qq \cD_n = \G_n
\cD \G_n,\\
\ns\ds G_n = \G_n G \G_n, \qq f_n = \G_n f,\\
\ns\ds g_n = \G_n g,\qq\;\, h_n = \G_n h.
\end{cases} \ee
It is easy to show that  for all $\zeta\in H$
and a.e. $(t,\om)\in [0,T]\times\Om$,
\begin{equation}\label{9.21-eq2}
\begin{cases} \ds
\lim_{n\to+\infty} e^{A_n t}\zeta = e^{A t}
\zeta \mbox{ in }H, ,\q\forall t\in [0,T],\\
\ns\ds
\lim_{n\to+\infty} \cA_n\zeta = \cA \zeta \mbox{ in }H,\\
\ns\ds \lim_{n\to+\infty} \cB_n\zeta = \cB \zeta\mbox{ in }H, \\
\ns\ds \lim_{n\to+\infty} \cD_n\zeta = \cD
\zeta\mbox{ in }H,
\end{cases}
\end{equation}
\begin{equation}\label{6.8-eq7}
\begin{array}{ll}\ds
\lim_{n\to+\infty} G_n\zeta = G \zeta \mbox{ in
}H,  \q \mbox{ for all }\zeta\in H \mbox{ and
a.e. } \om \in \Om,
\end{array}
\end{equation}
and
\begin{equation}\label{9.21-eq1.1}
\begin{cases}\ds
\lim_{n\to\infty} f_n = f\q\mbox{ in
}\; L^2_\dbF(0,T;H), \\
\ns\ds \lim_{n\to\infty} g_n = g \q\mbox{ in
}\; L^2_\dbF(0,T;H), \\
\ns\ds \lim_{n\to\infty} h_n = h  \q\mbox{ in
}\; L^2_\dbF(0,T;H).
\end{cases}
\end{equation}

For any $\xi\in D(A)$,
$$
\begin{array}{ll}\ds
\lim_{n\to\infty}|A_n \xi - A\xi|_H\\
\ns\ds =
\lim_{n\to\infty}|\G_n A\G_n \xi - A\xi|_H \\
\ns\ds \leq \lim_{n\to\infty}|\G_n( A\G_n \xi -
A \xi)|_H + \lim_{n\to\infty}|(\G_n -I) A \xi|_H \\
\ns \ds \leq \lim_{n\to\infty}|(A\G_n \xi - A
\xi)|_H + \lim_{n\to\infty}|(\G_n -I) A
\xi|_H \\
\ns \ds\leq
\lim_{n\to\infty}|A|_{\cL(D(A);H)}|(\G_n \xi -
\xi)|_{D(A)} + \lim_{n\to\infty}|(\G_n -I) A
\xi|_H =0.
\end{array}
$$
By the Trotter-Kato approximation theorem (e.g.
\cite[page 209]{EN}), we have that, for any
$\zeta\in H$,
\begin{equation}\label{9.21-eq1}
\lim_{n\to+\infty} e^{A_n t}\zeta = e^{At}\zeta
\mbox{ in } H, \q \mbox{uniformly for }\, t\in
[0,T].
\end{equation}
Similarly, we can prove that for any $\zeta\in
V$,
\begin{equation}\label{9.21-eq1.1xx}
\lim_{n\to+\infty} e^{A_n t}\zeta = e^{At}\zeta
\mbox{ in } V, \q \mbox{uniformly for }\, t\in
[0,T].
\end{equation}

We need the following result.

\begin{lemma}\label{lm13}
Let (AS2) hold, and $V\in \cV_H$. Then
$A\in\cL_2(H;V)$ and $ \ds\lim_{n\to\infty} |A_n
-A|_{\cL_2(H;V)}=0. $

\end{lemma}

{\it Proof}\,: Since $\{e_j\}_{j=1}^\infty$ are
eigenvectors of $A$, it holds that
$$
\begin{array}{ll}\ds
\sum_{j=1}^\infty|A e_j|_V^2
\3n&\ds=\sum_{j=1}^\infty|\mu_{j}e_j|_V^2 =
\sum_{j=1}^\infty \lambda_j^2
|\mu_j|^2|e_j|_{D(A)}^{-2} \leq
\sum_{j=1}^\infty \lambda_j^2.
\end{array}
$$
Hence, $A\in \cL_2(H;V)$.

\vspace{0.1cm}

Next,
$$
\begin{array}{ll}\ds
 \lim_{n\to\infty} |A_n - A|_{\cL_2(H;V)}^2
\3n&\ds= \lim_{n\to\infty} \sum_{j=1}^\infty
|(A_n - A)e_j|_V^2\\
\ns&\ds =  \lim_{n\to\infty} \sum_{j=n+1}^\infty
\lambda_j^2
|\mu_j|^2|e_j|_{D(A)}^{-2}\\
\ns&\ds\leq \lim_{n\to\infty}\sum_{j=n+1}^\infty
\lambda_j^2 =0.
\end{array}
$$

For any $s\in[0,T)$ and $n\in\dbN$, let us
consider the following two equations:
\begin{equation}\label{6.7-eq1.1}
\left\{
\begin{array}{ll}
\ds dx_n= \big[ (A_n+\cA_n) x_n  + f_n \big]dt+
(\cB_n x_n+g_n)
dW(t)& \mbox{ in }[s,T], \\
\ns\ds x_n(0)=\G_n\eta,
\end{array}
\right.
\end{equation}
and
\begin{equation}\label{6.7-eq1.2}
\left\{
\begin{array}{ll}
\ds  dy_n=-\big(A^*_n y_n+ \cD_n z_n+
h_n\big)dt+ z_ndW(t) &
\mbox{ in }[s,T],\\
\ns\ds  y_n(T)= \G_n \xi.
\end{array}
\right.
\end{equation}
By \eqref{9.21-eq2}, \eqref{9.21-eq1.1} and
\eqref{9.21-eq1}, it is easy to show the
following result (Hence we omit the details).
\begin{lemma}\label{lm7}
Let (AS2) hold. For any $\eta\in
L^2_{\cF_s}(\Om;H)$ and $\xi\in
L^2_{\cF_T}(\Om;H)$, it holds that
\begin{equation}\label{lm7-eq1}
\left\{ \ba{ll}\ds
\lim_{n\to\infty}\mE\big(\sup_{t\in
[s,T]}|x_n(t)-x(t)|_{H}^2\big) =0,\\
\ns \ds \lim_{n\to\infty}\mE\big(\sup_{t\in
[s,T]}|y_n(t)-y(t)|_{H}^2\big) =0,\\
\ns\ds \lim_{n\to\infty}
|z_n(\cd)-z(\cd)|_{L^2_\dbF(s,T;H)} =0. \ea
\right.
\end{equation}
\end{lemma}
\begin{remark}
The equalities \ref{lm7-eq1} hold without the
assumption (AS2). In this paper, we only need
them under the assumption (AS2). Hence, we do
not deal with the general case for avoiding
technical complexity.
\end{remark}
{\it Proof of Lemma \ref{lm7}}\,: Put
$$
M\triangleq\sup_{t\in[s,T]}|e^{A(t-s)}|_{\cL(H)}\big(|\cA|_{\cL(H)}+
|\cB|_{\cL(H)}\big)
$$
and let $t_0\in [s,T]$ satisfying that
\begin{equation}\label{6.1-eq9}
t_0\triangleq \max\Big\{t\in [s,T]:\max\{t-s,
(t-s)^2\}\leq \frac{1}{4M^2}\Big\}.
\end{equation}

From \eqref{6.20-eq1} and \eqref{6.7-eq1.1}, we
have that
\begin{eqnarray}\label{6.1-eq1}
&& \mE\big(\sup_{t\in
[s,t_0]}|x_n(t)-x(t)|_{H}^2\big)\nonumber\\
&& \leq \cC\mE\sup_{t\in [s,T]}\(\Big|
e^{A(t-s)}\eta
-e^{A_n(t-s)}\G_n\eta \Big|_H^2\nonumber\\
&&\qq + \Big| \int_s^t e^{A(r-s)}\cA x(r)dr
-\int_s^t
e^{A_n(r-s)}\cA_nx_n(r)dr \Big|_H^2\\
&&\qq + \Big| \int_s^t e^{A(r-s)}f(r)dr
-\int_s^t e^{A_n(r-s)}f_n(r)dr \Big|_H^2\nonumber\\
&&\qq  + \Big| \int_s^t e^{A(r-s)}\cA x(r)dW(r)
-\int_s^t e^{A_n(r-s)}\cA_nx_n(r)dW(r)
\Big|_H^2\nonumber\\
&&\qq + \Big| \int_s^t e^{A(r-s)}g(r)dr
-\int_s^t e^{A_n(r-s)}g_n(r)dW(r)
\Big|_H^2\).\nonumber
\end{eqnarray}
Let us estimate the terms in the right hand side
of \eqref{6.1-eq1} one by one. First,
\begin{equation}\label{6.1-eq2}
\begin{array}{ll}\ds
\sup_{t\in [s,t_0]} \big| e^{A(t-s)}\eta
-e^{A_n(t-s)}\G_n\eta \big|_H^2\\
\ns\ds \leq  2\sup_{t\in [s,t_0]} \big|
e^{A(t-s)}\eta -e^{A_n(t-s)} \eta\big|_H^2 +
2\sup_{t\in [s,t_0]}\big| e^{A_n(t-s)}\eta
-e^{A_n(t-s)}\G_n\eta \big|_H^2\\
\ns\ds \leq 2\sup_{t\in [s,t_0]} \big|
e^{A(t-s)}\eta -e^{A_n(t-s)} \eta\big|_H^2 +
\cC\mE\big|\eta-\G_n\eta \big|_H^2.
\end{array}
\end{equation}
Next, from the definition of $A_n$ and $\cA_n$,
we know that
\begin{equation}\label{6.1-eq8}
\big(e^{A(t-s)}-e^{A_n(t-s)}\big)\cA_n=0,\qq
\forall t\in [s,T].
\end{equation}
Thus, we have
\begin{equation}\label{6.1-eq4}
\begin{array}{ll}\ds
\mE \sup_{t\in [s,t_0]}\Big| \int_s^t
e^{A(r-s)}\cA x(r)dr -\int_s^t
e^{A_n(r-s)}\cA_nx_n(r)dr \Big|_H^2\\
\ns\ds \leq 2\mE \sup_{t\in [s,t_0]}\Big|
\int_s^t e^{A(r-s)}\cA x(r)dr
-\int_s^t e^{A_n(r-s)}\cA_nx(r)dr \Big|_H^2\\
\ns\ds\q +2\mE \sup_{t\in [s,t_0]}\Big| \int_s^t
e^{A_n(r-s)}\cA_nx(r)dr -\int_s^t
e^{A_n(r-s)}\cA_nx_n(r)dr \Big|_H^2\\
\ns\ds\leq 2\mE \int_s^{t_0}
\big|\big[\big(e^{A(r-s)}  -
e^{A_n(r-s)}\big)\cA +\big(e^{A(r-s)}
- e^{A_n(r-s)}\big)\cA_n\big]x(r)\Big|_H dr \\
\ns\ds\q + M^2(t_0-s)^2\sup_{t\in [s,t_0]}
|x(t)- x_n(t)|_H^2\\
\ns\ds \leq 2\mE \int_s^{t_0} \big|
\big(e^{A(r-s)}  - e^{A_n(r-s)}\big)\cA
x(r)\Big|_H dr \\
\ns\ds\q + M^2(t_0-s)^2\sup_{t\in [s,t_0]}
|x(t)- x_n(t)|_H^2
\end{array}
\end{equation}
Using \eqref{6.1-eq8} again, we get that
\begin{eqnarray}\label{6.1-eq7}
&& \mE \sup_{t\in [s,t_0]}\Big|\int_s^t
e^{A(r-s)}f(r)dr -\int_s^t e^{A_n(r-s)}f_n(r)dr
\Big|_H^2 \nonumber\\
&& \leq 2\mE \sup_{t\in [s,t_0]}\Big| \int_s^t
e^{A(r-s)}g(r)dW(r)
-\int_s^t e^{A_n(r-s)}g(r)dW(r) \Big|_H^2\nonumber\\
&&\q +2\mE \sup_{t\in [s,t_0]}\Big| \int_s^t
e^{A_n(r-s)}f(r)dr -\int_s^t
e^{A_n(r-s)}f_n(r)dr \Big|_H^2\\
&&\leq 2\mE  \int_s^{t_0} \big|\big(e^{A(r-s)} -
e^{A_n(r-s)} \big) f(r)\big|_H^2dr +
\cC\mE\int_s^{t_0}|f(r)- f_n(r)|_H^2dr.\nonumber
\end{eqnarray}

By Burkholder-Davis-Gundy inequality and
\eqref{6.1-eq8}, we have that
\begin{eqnarray}\label{6.1-eq5}
&& \mE \sup_{t\in [s,t_0]}\Big| \int_s^t
e^{A(r-s)}\cB x(r)dW(r) -\int_s^t
e^{A_n(r-s)}\cB_nx_n(r)dW(r)
\Big|_H^2 \nonumber\\
&& \leq 2\mE \sup_{t\in [s,t_0]}\Big| \int_s^t
e^{A(r-s)}\cB x(r)dW(r)
-\int_s^t e^{A_n(r-s)}\cB_nx(r)dW(r) \Big|_H^2\nonumber\\
&&\q +2\mE \sup_{t\in [s,t_0]}\Big| \int_s^t
e^{A_n(r-s)}\cB_nx(r)dW(r) -\int_s^t
e^{A_n(r-s)}\cB_nx_n(r)dW(r) \Big|_H^2\\
&&\leq 2\mE \int_s^{t_0} \big|\big(e^{A(r-s)}
- e^{A_n(r-s)} \big)\cB x(r)\big|_H^2dr \nonumber\\
&&\q + M^2(t_0-s) \sup_{t\in [s,t_0]} |x(t)-
x_n(t)|_H^2\nonumber
\end{eqnarray}
and
\begin{eqnarray}\label{6.1-eq6}
&& \mE \sup_{t\in [s,t_0]}\Big|\int_s^t
e^{A(r-s)}g(r)dr -\int_s^t
e^{A_n(r-s)}g_n(r)dW(r)
\Big|_H^2 \nonumber\\
&& \leq 2\mE \sup_{t\in [s,t_0]}\Big| \int_s^t
e^{A(r-s)}g(r)dW(r)
-\int_s^t e^{A_n(r-s)}g(r)dW(r) \Big|_H^2\nonumber\\
&&\q +2\mE \sup_{t\in [s,t_0]}\Big| \int_s^t
e^{A_n(r-s)}g(r)dW(r) -\int_s^t
e^{A_n(r-s)}g_n(r)dW(r) \Big|_H^2\\
&&\leq 2\mE  \int_s^{t_0} \big|\big(e^{A(r-s)} -
e^{A_n(r-s)} \big) g(r)\big|_H^2dr +
\cC\mE\int_s^{t_0}|g(t)- g_n(t)|_H^2.\nonumber
\end{eqnarray}

From \eqref{6.1-eq9} to \eqref{6.1-eq6}, we find
that
\begin{eqnarray}\label{6.1-eq1.1}
&& \mE\big(\sup_{t\in
[s,t_0]}|x_n(t)-x(t)|_{H}^2\big)\nonumber\\
&& \leq \cC\(\mE\sup_{t\in [s,T]} \big|
e^{A(t-s)}\eta -e^{A_n(t-s)} \eta\big|_H^2 +
\mE\big|\eta-\G_n\eta \big|_H^2\nonumber\\
&&\q  + \mE \int_s^{t_0} \big|\big(e^{A(r-s)}
- e^{A_n(r-s)}\big)\cA x(r)\Big|_H dr\\
&&\q  + \mE  \int_s^{t_0} \big|\big(e^{A(r-s)}-
e^{A_n(r-s)} \big) f(r)\big|_H^2dr +
 \mE\int_s^{t_0}|f(r)- f_n(r)|_H^2dr \nonumber\\
&&\q   + \mE  \int_s^{t_0} \big|\big(e^{A(r-s)}-
e^{A_n(r-s)} \big)\cB x(r)\big|_H^2dr\nonumber\\
&&\q  + \mE  \int_s^{t_0} \big|\big(e^{A(r-s)}-
e^{A_n(r-s)} \big) g(r)\big|_H^2dr +
 \mE\int_s^{t_0}|g(t)- g_n(t)|_H^2\).\nonumber
\end{eqnarray}

For all $n\in\dbN$,
$$
\sup_{t\in [s,T]} \big| e^{A(t-s)}\eta
-e^{A_n(t-s)} \eta \big|_H^2 \leq \cC|\eta|_H^2.
$$
This, together with Lebesgue's dominated
convergence theorem and \eqref{9.21-eq2},
implies that
\begin{equation}\label{6.1-eq3.1}
\begin{array}{ll}\ds
\lim_{n\to\infty}\mE\sup_{t\in [s,T]} \big|
e^{A(t-s)}\eta
-e^{A_n(t-s)}\eta \big|_H^2\\
\ns\ds =  \lim_{n\to\infty}\mE\sup_{t\in [s,T]}
\big| e^{A(t-s)}\eta -e^{A_n(t-s)}\eta
\big|_H^2=0.
\end{array}
\end{equation}
Similarly, we can prove that
\begin{equation}\label{6.1-eq3.2}
\lim_{n\to\infty}\mE\big|\eta-\G_n\eta
\big|_H^2=0.
\end{equation}
Next, noting that for all $n\in\dbN$,
$$
\big|\big(e^{A(r-s)}- e^{A_n(r-s)} \big) \cA
x(r)\big|_H^2  \leq \cC \big|x(r)\big|_H^2,\q
\mbox{ for a.e. }(t,\om)\in [s,t_0]\times\Om,
$$
it follows from Lebesgue's dominated convergence
theorem and \eqref{9.21-eq2} that
\begin{equation}\label{6.1-eq3.3}
\begin{array}{ll}\ds
\lim_{n\to\infty}\mE\int_s^{t_0}
\big|\big(e^{A(r-s)}- e^{A_n(r-s)} \big)
f(r)\big|_H^2dr\\
\ns\ds = \lim_{n\to\infty}\mE\int_s^{t_0}
\big|\big(e^{A(r-s)}- e^{A_n(r-s)} \big) \cA
x(r)\big|_H^2dr
\end{array}
\end{equation}
Similar to the above arguments, we can show that
all the terms in the right hand side of
\eqref{6.1-eq1.1} tend to zero as $n$ tends to
$\infty$. Consequently, we obtain that
$$
\lim_{n\to\infty}\mE\big(\sup_{t\in
[s,t_0]}|x_n(t)-x(t)|_{H}^2\big) =0.
$$

If $t_0=T$, then we complete our proof.
Otherwise, let
\begin{equation}\label{6.1-eq9.1}
t_1\triangleq \max\Big\{t\in
[t_0,T]:\max\{t-t_0, (t-t_0)^2\}\leq
\frac{1}{4M^2}\Big\}.
\end{equation}
Repeating the above argument, we get that
$$
\lim_{n\to\infty}\mE\big(\sup_{t\in
[s,t_1]}|x_n(t)-x(t)|_{H}^2\big) =0.
$$
By an induction argument, we can obtain that
$$
\lim_{n\to\infty}\mE\big(\sup_{t\in
[s,T]}|x_n(t)-x(t)|_{H}^2\big) =0.
$$

Consider the following stochastic differential
equation:
\begin{equation}\label{op-fsystem2.1}
\left\{
\begin{array}{ll}
\ds dx_{1,n} = \big[(A+A_1) x_{1,n}+ \G_n u_1\big]d\tau + \big(\G_n C x_{1,n}  + \G_n v_1 \big)dW(\tau) &\mbox{ in } (t,T],\\
\ns\ds x_{1,n}(t)=\G_n\xi_1.
\end{array}
\right.
\end{equation}
Similar to Lemma \ref{lm7}, we can easily
establish the following result.
\begin{lemma}\label{lm16}
Assume that (AS2)--(AS3) holds.  Then, for any
$\xi_1\in L^4_{\cF_t}(\Om;V')$ and $u_1(\cd),
v_1(\cd)$ $\in L^4_\dbF(\Om;L^2(t,T;V'))$, the
solution $x_{1,n}(\cd)\in
L^4_\dbF(\Om;C([0,T];V'))$ to
\eqref{op-fsystem2.1} satisfies
\begin{equation}\label{xxx9.21-eq3}
\lim_{n\to\infty}
|x_{1,n}(\cdot)-x_1(\cdot)|_{L^4_\dbF(\Om;C([t,T];V'))}
=0,
\end{equation}
where $x_1(\cd)$ is the solution to
\eqref{op-fsystem1}.
\end{lemma}

For a.e. $\tau\in [s,T]$, let us define six
operators $\Phi$, $\Phi_n$, $\Psi$, $\Psi_n$,
$\Xi$ and $\Xi_n$ as follows:
$$
\begin{array}{ll}\ds \left\{
\begin{array}{ll}\ds
\Phi: L^2_{\cF_s}(\Om;H)\to
L^2_\dbF(\Om;C([s,T];H)),\\
\ns\ds (\Phi\eta)(\tau) = x(\tau)
\end{array}
\right. \\
\ns\ds \left\{
\begin{array}{ll}\ds
\Phi_n: L^2_{\cF_s}(\Om;H)\to
L^2_\dbF(\Om;C([s,T];H)),\\
\ns\ds (\Phi_n\eta)(\tau) = x_n(\tau),
\end{array}
\right.
\end{array}
$$
$$
\begin{array}{ll}
\left\{
\begin{array}{ll}\ds
\Psi: L^2_{\cF_s}(\Om;H)\to
L^2_\dbF(\Om;C([s,T];H)),\\
\ns\ds (\Psi\eta)(\tau) = y(\tau)
\end{array}
\right. \\
\ns\ds \left\{
\begin{array}{ll}\ds
\Psi_n: L^2_{\cF_s}(\Om;H)\to
L^2_\dbF(\Om;C([s,T];H)),\\
\ns\ds (\Psi_n\eta)(\tau) = y_n(\tau),
\end{array}
\right.
\end{array}
$$
and
$$
\begin{array}{ll}
\left\{
\begin{array}{ll}\ds
\Xi: L^2_{\cF_s}(\Om;H)\to
L^2_\dbF(s,T;H),\\
\ns\ds (\Xi\eta)(\tau) = z(\tau)
\end{array}
\right. \\
\ns\ds \left\{
\begin{array}{ll}\ds
\Xi_n: L^2_{\cF_s}(\Om;H)\to
L^2_\dbF(s,T;H),\\
\ns\ds (\Xi_n\eta)(\tau) = z_n(\tau).
\end{array}
\right.
\end{array}
$$
Here $\eta\in L^2_{\cF_s}(\Om;H)$, $x(\cd)$
(\resp $x_n(\cd)$) is the solution to
\eqref{6.20-eq1} (\resp \eqref{6.7-eq1.1}) with
$f=g=0$, $(y(\cd),z(\cd))$ (\resp
$(y_n(\cd),z_n(\cd))$) is the solution to
\eqref{6.20-eq10} (\resp \eqref{6.7-eq1.2}) with
$h$ and $\xi$ replaced by $\cK x$ for some
$\cK\in L^\infty_\dbF(0,T;\cL(H))$ and $Gx(T)$
(\resp $h_n$ and $\xi$ replaced by $\cK_n x_n$
with $\cK_n = \G_n \cK \G_n$ and $Gx_n(T)$),
respectively.

\vspace{0.1cm}

Denote by $I_{HV}$ the embedding operator from
$H$ to $V$. We have the following result.

\begin{lemma}\label{lm8}
If $\cA\in \Upsilon_1(V)$, $\cB\in
\Upsilon_2(V)$, and $\cD, \cK\in
L^\infty_\dbF(0,T;\cL(V))$, then,
\begin{equation}\label{lm8-eq1}
\lim_{n\to\infty}|I_{HV}\Phi_n -
I_{HV}\Phi|_{L^4_\dbF(\Om;C([s,T];\cL_2(H;V)))}=0,
\end{equation}
\begin{equation}\label{lm8-eq2}
\begin{cases}
 \ds \lim_{n\to\infty}|I_{HV}\Psi_n -
I_{HV}\Psi|_{L^4_\dbF(\Om;C([s,T];\cL_2(H;V)))}=0,\\
\ns\ds \lim_{n\to\infty}|I_{HV}\Xi_n -
I_{HV}\Xi|_{L^4_\dbF(\Om;L^2(s,T;\cL_2(H;V)))}=0.
\end{cases}
\end{equation}
\end{lemma}

{\it Proof}\,: We first prove \eqref{lm8-eq1}.
 It is easy to show that, for any
$\varrho\in L_\dbF^\infty(0,T;$ $L^2(\Om;V))$,
\begin{equation}\label{9.25-eq6}
\lim_{n\to\infty}|\cA_n\varrho-\cA\varrho|_{L^4_\dbF(\Om;L^1(0,T;V))}=0
\end{equation}
and
\begin{equation}\label{9.25-eq6.1}
\lim_{n\to\infty}|\cB_n\varrho-\cB\varrho|_{L^4_\dbF(\Om;L^2(0,T;V))}=0.
\end{equation}

From the definitions of $\Phi$ and $\Phi_n$, we
see that, for any $(s,\eta)\in[0,T)\times
L^2_{\cF_s}(\Om;H)$, $t\in[s,T]$ and
$\dbP\mbox{-a.s.}$,
$$
\begin{array}{ll}\ds
\Phi(t)\eta\3n&\ds = e^{A(t-s)}\eta + \int_s^t
e^{A(t-r)}\cA(r) \Phi(r)\eta dr \\
\ns&\ds\q+ \int_s^t e^{A(t-r)}\cB(r)\Phi(r)\eta
dW(r)\;\, \mbox{ in } H
\end{array}
$$
and
$$
\begin{array}{ll}\ds
\Phi_n(t)\eta \3n&\ds =  e^{A_n(t-s)}\G_n\eta
+\int_s^t e^{A_n(t-r)}\cA_n(r) \Phi_n(r)\eta
dr\\
\ns&\ds\q + \int_s^t e^{A_n(t-r)}
\cB_n(r)\Phi_n(r)\eta dW(r)\;\, \mbox{ in } H.
\end{array}
$$
Thus, by
$$e^{A_n(t-s)}\G_n\eta=e^{A_n(t-s)}\eta,$$  it
holds that
\begin{equation}\label{10.10-eq38}
\begin{array}{ll}\ds
I_{HV}\Phi(t)\eta\3n&\ds = I_{HV}e^{A(t-s)}\eta
+
\int_s^t I_{HV}e^{A(t-r)} \cA(r) \Phi(r)\eta dr \\
\ns&\ds \q+ \int_s^t I_{HV}e^{A(t-r)}
\cB(r)\Phi(r)\eta dW(r)\;\, \mbox{ in } V
\end{array}
\end{equation}
and
$$
\begin{array}{ll}\ds
I_{HV}\Phi_n(t)\eta\3n&\ds = I_{HV}
e^{A_n(t-s)}\eta +\int_s^t
I_{HV}e^{A_n(t-r)} \cA_n(r) \Phi_n(r)\eta dr \\
\ns&\ds \q + \int_s^tI_{HV} e^{A_n(t-r)}
\cB_n(r)\Phi_n(r)\eta dW(r)\;\, \mbox{ in } V.
\end{array}
$$

Noting that $I_{HV}e_j=e_j$ for all $j\in\dbN$,
if $O\in \cL(H)$ can be extended to a bounded
linear operator on $V$, then
$$
\begin{array}{ll}\ds
\big|I_{HV} O -O
I_{HV}\big|_{\cL_2(H;V)}^2\3n&\ds =
\sum_{j=1}^\infty\big|I_{HV} O e_j - O
I_{HV}e_j\big|^2_V\\
\ns&\ds= \sum_{j=1}^\infty\big|O e_j - O
e_j\big|^2_V=0.
\end{array}
$$
Consequently,  $I_{HV} O =O I_{HV}$. This,
together with \eqref{10.10-eq38}, implies that
$$
\begin{array}{ll}\ds
I_{HV}\Phi(t)\eta\3n&\ds = I_{HV}e^{A(t-s)}\eta
+ \int_s^t e^{A(t-r)} \cA(r) I_{HV}\Phi(r)\eta
dr\\
\ns&\ds\q + \int_s^t e^{A(t-r)}
\cB(r)I_{HV}\Phi(r)\eta dW(r)\;\, \mbox{ in } V.
\end{array}
$$
Since $\cL_2(H;V)$ is a Hilbert space, for any
$t\in[s,T]$,
\begin{equation}\label{10.10-eq40}
\begin{array}{ll}\ds
I_{HV}\Phi(t) \3n&\ds = I_{HV}e^{A(t-s)}  +
\int_s^t e^{A(t-r)} \cA(r) I_{HV}\Phi(r)  dr \\
\ns&\ds \q+ \int_s^t e^{A(t-r)}
\cB(r)I_{HV}\Phi(r)  dW(r)\;\, \mbox{ in }
\cL_2(H;V),\;\dbP\mbox{-a.s.}
\end{array}
\end{equation}
Similarly, we can prove that for any
$t\in[s,T]$,
\begin{equation}\label{10.10-eq39}
\begin{array}{ll}\ds
I_{HV}\Phi_n(t)\3n&\ds = I_{HV} e^{A_n(t-s)}
+\int_s^t
e^{A_n(t-r)} \cA_n(r) I_{HV}\Phi_n(r) dr \\
\ns&\ds \q + \int_s^t e^{A_n(t-r)}
\cB_n(r)I_{HV}\Phi_n(r) dW(r)\;\, \mbox{ in }
\cL_2(H;V),\;\dbP\mbox{-a.s.}
\end{array}
\end{equation}
In what follows, to simplify the notations, we
omit the operator $I_{HV}$ if there is no
confusion.

\vspace{0.1cm}

It follows from \eqref{10.10-eq40} and
\eqref{10.10-eq39} that for any stopping time
$\tau_0$ with $\tau_0(\om)\in (s,T]$,
$\dbP$-a.s.,
\begin{equation}\label{9.25-eq7}
\begin{array}{ll}\ds
\mE\sup_{r\in [s,\tau_0]}\big|\Phi(r) - \Phi_n(r)\big|_{\cL_2(H;V)}^4 \\
\ns\ds \leq \cC\mE\sup_{r\in
[s,\tau_0]}\[\big|e^{A(r-s)}-e^{A_n(r-s)}\big|_{\cL_2(H;V)}^4\\
\ns\ds \q+\Big|\int_s^r \big(e^{A(r-\tau)}
\cA(\tau)\Phi(\tau) - e^{A_n(r-\tau)}
\cA_n(\tau)
\Phi_n(\tau)\big)d\tau\Big|_{\cL_2(H;V)}^4\\
\ns\ds \q + \Big|\int_s^r
\big(e^{A(r-\tau)}\cB(\tau)\Phi(\tau) -
e^{A_n(r-\tau)}\cB_n(\tau)\Phi_n(\tau) \big)
dW(\tau)\Big|_{\cL_2(H;V)}^4\].
\end{array}
\end{equation}
Hence, by Burkholder-Davis-Gundy inequality, and
$|\cA_n(\cd) |_{\cL(V)}\le |\cA(\cd) |_{\cL(V)}$
and $|\cB_n(\cd) |_{\cL(V)}\le |\cB(\cd)
|_{\cL(V)}$, we deduce that
\begin{eqnarray}\label{9.25-eq7.1}
&& \mE\sup_{r\in [s,\tau_0]}\big|\Phi(r) -
\Phi_n(r)\big|_{\cL_2(H;V)}^4 \nonumber\\
&&  \leq \cC\[\mE\(\sup_{r\in
[s,\tau_0]}\big|e^{A(r-s)}-e^{A_n(r-s)}\big|_{\cL_2(H;V)}^4\)
\nonumber\\ && \q +\mE\(\int_s^{\tau_0}
\Big|\big(e^{A(r-\tau)}
\cA(\tau)-e^{A_n(r-\tau)}
\cA_n(\tau)\big)\Phi(\tau)\nonumber\\
&& \qq\qq\q  + e^{A_n(r-\tau)} \cA_n(\tau)
\big(\Phi(\tau)-\Phi_n(\tau)\big)\Big|_{\cL_2(H;V)}
d\tau\)^4 \nonumber\\
&&  \q + \mE\( \int_s^{\tau_0}
\Big|\big(e^{A(r-\tau)}\cB(\tau) -
e^{A_n (r-\tau)}\cB_n(\tau)\big) \Phi(\tau)\\
&& \qq\qq\q + e^{A_n (r-\tau)}\cB_n (\tau)
\big(\Phi(\tau) - \Phi_n (\tau)\big)
\Big|_{\cL_2(H;V)}^2 d\tau \)^2\] \nonumber\\
&&   \leq  \cC
\[\mE\( \sup_{r\in
[s,\tau_0]}
|e^{A(r-s)}-e^{A_n(r-s)}|_{\cL_2(H;V)}^4\)\nonumber\\
&& \qq   + \mE\Big( \int_0^T
\Big|\big(e^{A(r-\tau)} \cA(\tau) -
e^{A_n(r-\tau)} \cA_n(\tau) \big)
\Phi(\tau)\Big|_{\cL_2(H;V)} d\tau\)^4\nonumber
 \\ &&   \qq  + \big(\big|
|\cA(\cd)|_{\cL(V)}\big|_{L_\dbF^\infty(\Omega;L^1(s;\tau_0))}^4
+ \big|
|\cB(\cd)|_{\cL(V)}\big|_{L_\dbF^\infty(\Omega;L^2(s;\tau_0))}^4\big)\nonumber\\
&& \qq\q\times\mE\big(\sup_{r\in
[s,\tau_0]}|\Phi(r) -
\Phi_n(r)|_{\cL_2(H;V)}^4\big) \nonumber\\
&&  \q  + \mE\(\int_0^T
\big|\big(e^{A(r-\tau)}\cB(\tau) -
e^{A_n(r-\tau)}\cB_n(\tau)\big)\Phi(\tau)\big|_{\cL_2(H;V)}^2
d\tau\)^2\].\nonumber
\end{eqnarray}

From \eqref{10.10-eq30}, and noting that $\cA\in
\Upsilon_1(V)$ and $\cB\in \Upsilon_2(V)$,  we
conclude that there is a stopping time
$\tau_0\in (s,T]$, $\dbP$-a.s., such that
$$
\begin{array}{ll}\ds
\cC\(\Big|\int_s^{\tau_0}\big|\cA(\tau)
\big|_{\cL(V)}d\tau\Big|_{L^\infty(\Omega)}^4 +
\Big|\int_s^{\tau_0}\big|\cB(\tau)
\big|_{\cL(V)}^2d\tau\Big|_{L^\infty(\Omega)}^2\)\leq\frac{1}{2}.
\end{array}
$$
For such kind of $\tau_0$, it follows from
\eqref{9.25-eq7} that
\begin{equation}\label{9.25-eq7.1xx}
\begin{array}{ll}\ds
\mE\big(\sup_{r\in [0,\tau_0]}|\Phi(r) - \Phi_n(r)|_{\cL_2(H;V)}^4\big)  \\
\ns\ds\leq  \cC \[\mE\(\sup_{r\in [s,\tau_0]}
|e^{A(r-s)} - e^{A_n(r-s)}|_{\cL_2(H;V)}^4\)\\
\ns\ds\qq +\mE\Big| \int_0^T \big(e^{A(r-\tau)}
\cA(\tau) - e^{A_n(r-\tau)} \cA_n(\tau) \big)
\Phi(\tau)d\tau\Big|_{\cL_2(H;V)}^4  \\
\ns\ds \qq + \mE\(\int_0^T
\big|\big(e^{A(r-\tau)}\cB(\tau) -
e^{A_n(r-\tau)}\cB_n(\tau)\big)\Phi(\tau)\big|_{\cL_2(H;V)}^2
d\tau\)^2\].
\end{array}
\end{equation}
Now,  \eqref{9.21-eq1} implies that
$$
\lim_{n\to\infty}\mE\(\sup_{r\in
[s,\tau_0]}|e^{A(r-s)}-e^{A_n(r-s)}|_{\cL_2(H;V)}^4\)
=0.
$$
From \eqref{9.21-eq1} and \eqref{9.25-eq6}, we
get that
$$
\lim_{n\to\infty}\mE\Big|\int_0^T
\big(e^{A(r-\tau)} \cA(\tau) - e^{A_n(r-\tau)}
\cA_n(\tau) \big)
\Phi(\tau)d\tau\Big|_{\cL_2(H;V)}^4=0
$$
and
$$
\lim_{n\to\infty}\mE\(\int_0^T
\big|\big(e^{A(r-\tau)}\cB(\tau) -
e^{A_n(r-\tau)}\cB_n(\tau)\big)\Phi(\tau)\big|_{\cL_2(H;V)}^2
d\tau\)^2=0.
$$
These, together with \eqref{9.25-eq7.1xx}, imply
that
\begin{equation}\label{9.25-eq8}
\lim_{n\to\infty} |\Phi_n -\Phi
|_{L^4_\dbF(\Om;L^\infty(s;\tau_0;\cL_2(H;V)))}=0.
\end{equation}
Repeating the above argument gives
\eqref{lm8-eq1}.

\vspace{0.2cm}

Next, we prove \eqref{lm8-eq2}.  It is easy to
see that, for any $\varrho\in
L_\dbF^\infty(0,T;L^2(\Om;V))$,
\begin{equation}\label{9.25-eq6xx}
\begin{cases}\ds
\lim_{n\to\infty}|\cK_n\varrho-\cK\varrho|_{L^\infty_\dbF(0,T;L^2(\Om;V))}=0,\\
\ns\ds
\lim_{n\to\infty}|\cD_n\varrho-\cD\varrho|_{L^\infty_\dbF(0,T;L^2(\Om;V))}=0.
\end{cases}
\end{equation}
Similar to the proof of \eqref{10.10-eq40}, we
obtain that for any $t\in[s,T]$,
\begin{equation}\label{10.10-eq19}
\begin{array}{ll}\ds
\Psi(t) \3n&\ds= e^{A(T-t)}G\Phi(T) + \int_t^T
e^{A(r-t)}\big(\cK(r) \Psi(r) + \cD(r)\Xi(r)\big) dr \\
\ns&\ds \q + \int_t^T e^{A(t-r)} \Xi(r)
dW(r)\;\, \mbox{ in }\; \cL_2(H;V),\q
\dbP\mbox{-a.s.}
\end{array}
\end{equation}
and
\begin{equation}\label{10.10-eq20}
\begin{array}{ll}\ds
\Psi_n(t) \3n&\ds= e^{A_n(T-t)}G_n\Phi_n(T)
+\int_t^T e^{A_n(r-t)}\big(\cK_n(r) \Psi_n(r) +
\cD_n(r)\Xi_n(r)\big)
dr\\
\ns&\ds\q + \int_t^T e^{A_n(r-t)}  \Xi_n(r)
dW(r)\q \mbox{ in }\;
\cL_2(H;V),\q\dbP\mbox{-a.s.}
\end{array}
\end{equation}
Since $\cL_2(H;V)$ is a Hilbert space, by
\eqref{10.10-eq19}--\eqref{10.10-eq20}, it is
easy to see that $(\Psi,\Xi)$ and
$(\Psi_n,\Xi_n)$ are respectively weak solutions
of the following $\cL_2(H;V)$-valued BSEEs
\begin{equation}\label{10.10-eq12}
\begin{cases}\ds
d\Psi = -(A\Psi+\cK\Psi+\cD\Xi)dt + \Xi dW(t)
&\mbox{ in } [s,T),\\
\ns\ds \Psi(T)=G\Phi(T)
\end{cases}
\end{equation}
and
\begin{equation}\label{10.10-eq12xx}
\begin{cases}\ds
d\Psi_n = -(A_n\Psi_n+\cK_n\Psi_n+\cD_n\Xi_n)dt
+ \Xi_n dW(t)
&\mbox{ in } [s,T),\\
\ns\ds \Psi_n(T)=G_n\Phi_n(T).
\end{cases}
\end{equation}
Then, for any $t\in (s,T]$, by It\^o's formula
and noting that $(A-A_n)\Psi_n=0$ (by our
assumption (\textbf{AS2})),
\begin{equation}\label{10.10-eq13}
\begin{array}{ll}\ds
|\Psi(t) - \Psi_n(t)|_{\cL_2(H;V)}^2 + \int_t^T |\Xi(r)-\Xi_n(r)|_{\cL_2(H;V)}^2dr  \\
\ns\ds =|G\Phi(T) - G_n\Phi_n(T)|_{\cL_2(H;V)}^2
+ 2\int_t^T \big\langle A(\Psi -\Psi_n), \Psi
-\Psi_n \big\rangle_{\cL_2(H;V)}dW(\tau)
 \\
\ns\ds\q + 2\int_t^T\! \big[\big\langle
(\cK-\cK_n)\Psi, \Psi-\Psi_n
\big\rangle_{\cL_2(H;V)}\!+ \big\langle
\cK_n(\Psi-\Psi_n),
\Psi-\Psi_n \big\rangle_{\cL_2(H;V)}\big]dW(\tau) \\
\ns\ds\q   + 2\int_t^T \big[ \big\langle
(\cD-\cD_n)\Xi, \Psi-\Psi_n
\big\rangle_{\cL_2(H;V)}+ \big\langle
\cD_n(\Xi-\Xi_n),
\Psi-\Psi_n \big\rangle_{\cL_2(H;V)}\big]dW(\tau) \\
\ns\ds\q   - 2\int_t^T\big\langle \Xi-\Xi_n ,
\Psi-\Psi_n \big\rangle_{\cL_2(H;V)} dW(\tau).
\end{array}
\end{equation}
Since $A$ generates a $C_0$-group on $V$, we
have that for any $\varrho\in \cL_2(H;V)$,
$$
\begin{array}{ll}\ds
\langle A \varrho,
\varrho\rangle_{\cL_2(H;V)}\3n&\ds =
\sum_{k=1}^\infty \langle A (\varrho e_k),
\varrho e_k\rangle_{V}\leq \cC\sum_{k=1}^\infty
|\varrho e_k|_{V}^2 \\[2mm]
\ns&\ds= \cC |\varrho|_{\cL_2(H;V)}^2.
\end{array}
$$
Thus,
\begin{equation}\label{10.10-eq15}
\int_t^T \big\langle A(\Psi-\Psi_n),
\Psi-\Psi_n\big\rangle_{\cL_2(H;V)}d\tau \leq
\cC\int_t^T |\Psi-\Psi_n|^2_{\cL_2(H;V)}d\tau.
\end{equation}
Clearly,
\begin{eqnarray}\label{10.10-eq16}
&& \int_t^T \big\langle \cK_n(\Psi-\Psi_n),
\Psi-\Psi_n \big\rangle_{\cL_2(H;V)} d\tau\nonumber\\
&&\q +\int_t^T \big\langle \cD_n(\Xi-\Xi_n),
\Psi-\Psi_n \big\rangle_{\cL_2(H;V)} d\tau \nonumber\\
&& \leq |\cK_n|_{L^\infty_\dbF(0,T;V)} \int_t^T
\big|\Psi-\Psi_n \big|_{\cL_2(H;V)}^2 d\tau \\
&&\q + |\cD_n|_{L^\infty_\dbF(0,T;V)} \int_t^T
\big|\Psi-\Psi_n
\big|_{\cL_2(H;V)}\big|\Xi-\Xi_n
\big|_{\cL_2(H;V)}
d\tau\nonumber\\
&&   \leq \cC\int_t^T \big|\Psi-\Psi_n
\big|_{\cL_2(H;V)}^2 d\tau + \frac{1}{2}\int_t^T
\big|\Xi-\Xi_n \big|_{\cL_2(H;V)}^2
d\tau.\nonumber
\end{eqnarray}
From \eqref{10.10-eq13}, \eqref{10.10-eq15} and
\eqref{10.10-eq16}, we find that
\begin{equation}\label{10.10-eq14}
\begin{array}{ll}\ds
\mE|\Psi(t) - \Psi_n(t)|_{\cL_2(H;V)}^4 + \mE\(\int_t^T |\Xi(r)-\Xi_n(r)|_{\cL_2(H;V)}^2dr\)^2 \\[1mm]
\ns\ds \leq   \cC\mE\[|G\Phi(T) -
G_n\Phi_n(T)|_{\cL_2(H;V)}^4  +\(\int_t^T |\Psi-\Psi_n|^2_{\cL_2(H;V)}d\tau\)^2 \\[1mm]
\ns\ds \q+ \(\int_t^T
\big|(\cK-\cK_n)\Psi\big|^2_{\cL_2(H;V)}
d\tau\)^2 +\int_t^T
\big|(\cD-\cD_n)\Xi\big|_{\cL_2(H;V)} d\tau\].
\end{array}
\end{equation}
Since
$$
\begin{array}{ll}\ds
\lim_{n\to\infty}\mE\(\int_s^T
\big|(\cK-\cK_n)\Psi\big|^2_{\cL_2(H;V)}
d\tau\)^2\\[1mm]
\ns\ds +\lim_{n\to\infty}\mE\(\int_s^T
\big|(\cD-\cD_n)\Xi\big|^2_{\cL_2(H;V)}
d\tau\)^2=0,
\end{array}
$$
the estimate \eqref{10.10-eq14} implies that
\begin{equation}\label{10.10-eq17}
\begin{array}{ll}\ds
\lim_{n\to\infty}\[\sup_{t\in[s,T]}\mE|\Psi(t) -
\Psi_n(t)|_{\cL_2(H;V)}^4\\[1mm]
\ns\ds\qq\q + \mE\(\int_s^T
|\Xi(r)-\Xi_n(r)|_{\cL_2(H;V)}^2dr\)^2\] =0.
\end{array}
\end{equation}
This gives the second equality in
\eqref{lm8-eq2}.

It remains to prove the first equality in
\eqref{lm8-eq2}. From \eqref{10.10-eq19} and
\eqref{10.10-eq20}, it follows that, for any
$s_0\in [s,T)$,
\begin{equation}\label{10.10-eq18}
\begin{array}{ll}\ds
\mE\sup_{r\in [s_0,T]}|\Psi(r) - \Psi_n(r)|_{\cL_2(H;V)}^4 \\[1mm]
\ns\ds \leq \cC\mE\sup_{r\in
[s_0,T]}\[|e^{A(T-r)}G\Phi(T)-e^{A_n(T-r)}G_n\Phi_n(T)|_{\cL_2(H;V)}^4\\[1mm]
\ns\ds \q+\Big|\int_r^T \big(e^{A(\tau-r)}
\cK(\tau)\Psi(\tau) - e^{A_n(\tau-r)}
\cK_n(\tau)
\Psi_n(\tau)\big)d\tau\Big|_{\cL_2(H;V)}^4\\[1mm]
\ns\ds \q + \Big|\int_r^T
\big(e^{A(\tau-r)}\cD(\tau)\Xi(\tau) -
e^{A_n(\tau-r)}\cD_n(\tau)\Xi_n(\tau)
\big) \tau \Big|_{\cL_2(H;V)}^4\\[1mm]
\ns\ds \q + \Big|\int_r^T  \big(\Xi(\tau) -
 \Xi_n(\tau) \big) dW(\tau)\Big|_{\cL_2(H;V)}^4\].
\end{array}
\end{equation}
Therefore, by Burkholder-Davis-Gundy inequality,
similar to \eqref{9.25-eq7.1}, we deduce that
\begin{equation}\label{10.10-eq18.1}
\begin{array}{ll}\ds
\mE\sup_{r\in [s_0,T]}|\Psi(r) - \Psi_n(r)|_{\cL_2(H;V)}^4 \\[1mm]
\ns\ds \leq \cC\[\mE\sup_{r\in
[s,T]}\Big\{|e^{A(T-r)}G\Phi(T)-e^{A_n(T-r)}G_n\Phi_n(T)|_{\cL_2(H;V)}^4\\[1mm]
\ns\ds \qq+\mE\(\int_0^T \Big|\big(e^{A(r-\tau)}
\cK(\tau) - e^{A_n(r-\tau)} \cK_n(\tau) \big)
\Psi(\tau) \Big|_{\cL_2(H;V)}d\tau\)^4 \\[1mm]
\ns\ds \qq + (T-s_0)^3 \big|\cK
\big|^4_{L^\infty_\dbF(0,T;\cL(V))}\mE\sup_{r\in [s_0,T]}|\Psi(r) - \Psi_n(r)|_{\cL_2(H;V)}^4\\[1mm]
\ns\ds \qq + \mE\(\int_0^T
\big|\big(e^{A(r-\tau)}\cD(\tau) -
e^{A_n(r-\tau)}\cD_n(\tau)\big)\Xi(\tau)\big|_{\cL_2(H;V)}^2
d\tau\)^2 \\[1mm]
\ns\ds\qq + \big(\big|\cD
\big|^4_{L^\infty_\dbF(0,T;\cL(V))}+1\big) |\Xi
- \Xi_n|_{L^2_\dbF(0,T;\cL_2(H;V))}^4\big)\].
\end{array}
\end{equation}

Let
$$s_0=T-(2\cC|\cK|^4_{L^\infty_\dbF(0,T;\cL(V))})^{-\frac{1}{3}}.$$
Then $$\cC(T-s_0)^3
|\cK|^4_{L^\infty_\dbF(0,T;\cL(V))}=\frac{1}{2}.$$
By \eqref{9.21-eq1.1}, \eqref{lm8-eq1} and
\eqref{10.10-eq18.1},  we conclude that that $$
\ds\lim_{n\to\infty} |\Psi_n -\Psi
|_{L^4_\dbF(\Om;L^\infty(s_0,T;\cL_2(H;V)))}=0.$$
Repeating this argument gives the first equality
in \eqref{lm8-eq2}.

\ms

\subsection{An auxilliary controllability result
}

\ms

In this section, we give a controllability
result concerning the trajectories of solutions
to \eqref{op-fsystem2}, which plays an important
role in the proof of the uniqueness of the
transposition solution to \eqref{5.5-eq6}.

\begin{lemma}\label{lm11}
The set
$$
\begin{array}{ll}\ds
\Xi\triangleq\big\{x_2(\cd)\;\big|\;
x_2(\cd)\mbox{ solves } \eqref{op-fsystem2}
\mbox{ with }t=0,\;\xi_2=0,\; v_2=0 \\
\ns\ds\hspace{2cm}\mbox{ and }u_2\in
L^4_{\dbF}(\Om;L^2(0,T;H)) \big\}
\end{array}
$$
is dense in $L^2_{\dbF}(0,T;H)$.
\end{lemma}

{\it Proof}\,: The proof of Lemma \ref{lm11} is
very similar to an intermediate step in the
proof of \cite[Theorem 4.1]{LZ1} (See \cite[pp.
38--39]{LZ1}). We give it here for the sake of
completeness.

If Lemma \ref{lm11} was not true, then there
would be a nonzero $r\in L^{2}_{\dbF}(0,T;H)$
such that
\begin{equation}\label{r}
\mE\int_0^T \big\langle r,x_2 \big\rangle_H ds =
0,\q\mbox{ for any } x_2\in \Xi.
\end{equation}
Consider the following $H$-valued BSEE:
\begin{equation}\label{xz2}
\left\{
\begin{array}{ll} \ds
dy=-A^*ydt+\big(r-A_1^* y-C^*
Y\big)dt+YdW(t)&\hb{in } [0,T),\\[1mm] \ns\ds
y(T)=0.
\end{array}
\right.
\end{equation}
The equation (\ref{xz2}) admits a unique
solution
$$(y(\cdot), Y(\cdot)) \in
L^{2}_{\dbF}(\Om;C([0,T];H)) \times
L^2_{\dbF}(0,T;H).$$ Hence, for any
$\phi_1(\cdot)\in L^1_{\dbF}(0,T;L^4(\Om;H))$
and $\phi_2(\cdot)\in
L^2_{\dbF}(0,T;L^4(\Om;H))$, it holds that
\begin{equation}\label{zx3}
\begin{array}{ll}\ds
\q - \dbE\int_0^T \big\langle z(s),r(s)
-A_1^*y(s)-C^*Y(s)\big\rangle_Hds\\[1mm]
\ns\ds = \dbE\int_0^T \big\langle
\phi_1(s),y(s)\big\rangle_H ds + \dbE\int_0^T
\big\langle \phi_2(s),Y(s)\big\rangle_H ds,
\end{array}
\end{equation}
where $z(\cd)$ solves
 \bel{xz4}
 \left\{
 \ba{ll}
 \ds dz=(Az+\phi_1)dt+\phi_2dW(t) &\hb{ in } (0,T],\\
 \ns\ds
 z(0)=0.
 \ea
 \right.
 \ee
In particular, for any $x_2(\cd)$ solving
\eqref{op-fsystem2} with $t=0$, $\xi_2=0$,
$v_2=0$ and an arbitrarily given $u_2\in
L^4_{\dbF}(0,T;H)$, we choose $z=x_2$,
$\phi_1=A_1x_2+u_2$ and $\phi_2=Cx_2$. By
(\ref{zx3}), it follows that
\begin{equation}\label{zx-5}
\begin{array}{ll}\ds
- \dbE\int_0^T \big\langle
x_2(s),r(s)\big\rangle_Hds= \dbE\int_0^T
\big\langle u_2(s),y(s)\big\rangle_H
ds,\\
\ns\ds\hspace{4.2cm}\q\forall\; u_2\in
L^4_{\dbF}(0,T;H).
\end{array}
\end{equation}
By (\ref{zx-5}) and recalling \eqref{r}, we
conclude that $y(\cd)=0$. Hence, (\ref{zx3}) is
reduced to
 \begin{equation}\label{zx-6}
- \dbE\int_0^T \big\langle
z(s),r(s)-C^*Y(s)\big\rangle_Hds = \dbE\int_0^T
\big\langle \phi_2(s),Y(s)\big\rangle_H ds.
\end{equation}
Choosing $\phi_2(\cd)=0$ in (\ref{xz4}) and
(\ref{zx-6}), we obtain that
\begin{equation}\label{zx-7}
\begin{array}{ll}\ds
\dbE\int_0^T \Big\langle
\int_0^sS(s-\si)\phi_1(\si)d\si,r(s)-C^*Y(s)\Big\rangle_Hds
= 0,\\
\ns\ds \qq\qq\qq\qq\qq\qq\forall\;
\phi_1(\cdot)\in L^1_{\dbF}(0,T;L^4(\Om;H)).
\end{array}
\end{equation}
Hence,
\begin{equation}\label{zx-8}
\int_\si^TS(s-\si)\big[r(s)-C^*Y(s)\big]ds
=0,\qq\forall\;\si\in [0,T].
\end{equation}
Then, for any given $\lambda_0\in\rho(A)$ and
$\si\in [0,T]$, we have
\begin{equation}\label{zx-9}
\begin{array}{ll}\ds
\q\int_\si^TS(s-\si)(\lambda_0-A)^{-1}\big[r(s)-C^*Y(s)\big]ds\\
\ns\ds
=(\lambda_0-A)^{-1}\int_\si^TS(s-\si)\big[r(s)-C^*Y(s)\big]ds
=0.
\end{array}
\end{equation}
Differentiating the equality (\ref{zx-9}) with
respect to $\si$, and noting (\ref{zx-8}), we
see that
 $$
\ba{ll}\ds
(\lambda_0-A)^{-1}\big[r(\si)-C^*Y(\si)\big]\\[1mm]
\ns\ds=-\int_\si^TS(s-\si)A(\lambda_0-A)^{-1}\big[r(s)-C^*Y(s)\big]ds\\\ns
\ds=\int_\si^T
S(s-\si)\big[r(s)-C^*Y(s)\big]ds\\\ns\ds
 \ds\q-\lambda_0\int_\si^TS(s-\si)(\lambda_0-A)^{-1}\big[r(s)-C^*Y(s)\big]ds\\[1mm]
 \ns\ds
 \ds=0,\qq\forall\;\si\in [0,T].
 \ea
 $$
Therefore,
 \bel{zx-10}
r(\cd)=C(\cd)^*Y(\cd).
 \ee
By (\ref{zx-10}), the equation (\ref{xz2}) is
reduced to
 \bel{xoz2}
 \left\{
 \ba{ll}
 \ds dy=-A^*ydt-A_1 ^*ydt+YdW(t),\qq\hb{in }
 [0,T),\\[1mm]
 \ns\ds
 y(T)=0.
 \ea
 \right.
 \ee
It is clear that the unique solution of
(\ref{xoz2}) is $(y(\cd),Y(\cd))=(0,0)$. Hence,
by (\ref{zx-10}), we conclude that $r(\cd)=0$,
which is a contradiction. Therefore, $\Xi$ is
dense in $L^2_{\dbF}(0,T;H)$.

\ms

\subsection{Some preliminaries on Malliavinian
Calculus}

\ms

In this section, we recall some basic results
for Malliavinian Calculus (See \cite{Nualart}
for the material on this topic).

\vspace{0.1cm}

Let $\widehat H$ be a separable Hilbert space.
Let $\xi\in L^2_{\cF_T}(\Om;\widehat H)$ be
Malliavinian differentiable. Denote by
$\mathbf{D}_\tau \xi$ the Malliavinian
derivative of $\xi$ at $\tau\in [0,T]$. Put
$$
\mathbb{D}^{1,2}(\widehat
H)\triangleq\Big\{\xi\in
L^2_{\cF_T}(\Om;\widehat
H)\;\Big|\;\mE\(|\xi|_{\widehat
H}^2+\int_0^T|\mathbf{D}_\tau\xi|_{\widehat
H}^2d\tau\)<\infty\Big\},
$$
which is a Hilbert space with the canonical
norm.

\vspace{0.1cm}

Denote by $\dbL^a_{1,2}(\widehat H)$ the set of
all $\widehat H$-valued, $\mathbf{F}$-adapted
processes $f(\cd)$ such that
\begin{enumerate}
  \item \ For  a.e. $t \in [0, T]$, $f(t,\cd)\in
  \mathbb{D}^{1,2}(\widehat
H)$;
  \item  \  $(t, \om) \mapsto \mathbf{D}_\tau f(t,\om)\in  L^2(0, T;\widehat
H)$ admits a progressively measurable version;

  \item \ $\ds|f|_{\dbL_{1,2}^a(\widehat
H)}\triangleq\sqrt{\mE\(\int_0^T|f(t)|_{\widehat
H}^2dt + \int_0^T\int_0^T|\mathbf{D}_\tau
f(t)|_{\widehat H}^2d\tau dt\)}<\infty$.
\end{enumerate}

\ms

We shall need the following result:
\begin{lemma}\label{10.25-lm1}\cite[Proposition 3.1]{Dou}
Let $Z\in L^2_\dbF(0,T;\wh H)$ and $t\in [0,T)$.
If $\xi\triangleq\int_t^T Z(r)dW(r)\in
\dbD^{1,2}(\wh H)$, then $Z\in
\dbL^a_{1,2}(\widehat H)$ and
$d\tau\otimes\dbP$-a.e.
$$
\cD_\tau \xi =\left\{ \begin{array}{ll}\ds
\int_t^T \mathbf{D}_\tau Z(r)dW(r), &\mbox{ if }
\tau\leq
t,\\[1mm]
\ns\ds Z(\tau)+\int_\tau^T \mathbf{D}_\tau
Z(r)dW(r), &\mbox{ if } \tau> t.
\end{array}
\right.
$$
\end{lemma}

\ms

Next, we present a result concerning the
regularity of solution to BSEE.

Assume that $\wh A$ generates a $C_0$-semigroup
on $\wh H$. Consider the following equation:
\begin{equation}\label{system3}
\left\{
\begin{array}{ll}
\ds dz +  A^*z dt=F(t,z,Z) dt + Z
dW(t)&\mbox{ in }[0,T), \\[1mm]
\ns\ds z(T)=z_T.
\end{array}
\right.
\end{equation}
Here   $F:\Om\times[0,T]\times \wh H \times \wh
H\to \wh H$ satisfies the following conditions:

\vspace{0.2cm}

1) $F$ is continuously differentiable with
respect to the second and third variables and
$$
|F_z(t,z,Z)|_{\cL(\wh H)}\leq L_1(t),\q
|F_Z(t,z,Z)|_{\cL(\wh H)}\leq L_2(t),\q
\mbox{a.e. } t\in [0,T]
$$
for some $L_1(\cd),L_2(\cd)\in L^2(0,T)$;

\vspace{0.1cm}

2) for each $(\eta_1,\eta_2)\in \wh H\times \wh
H$, $F(\cd,\eta_1,\eta_2)\in
L^1_\dbF(0,T;L^2(\Om;\wh H))\cap
\dbL^a_{1,2}(\wh H)$;

\vspace{0.1cm}

3) for all $t\in [0,T]$ and
$(\eta_1,\eta_2,\eta_3,\eta_4)\in \wh H\times
\wh H\times \wh H\times \wh H$, and for a.e.
$\tau\in [0,T]$,\vspace{1mm}
$$
\big|\mathbf{D}_\tau
F(t,\eta_1,\eta_2)-\mathbf{D}_\tau
F(t,\eta_3,\eta_4)\big|_{\wh H}\leq L_3(t)
|\eta_1-\eta_2|_{\wh
H}+L_4(t)|\eta_3-\eta_4|_{\wh H},
$$
where $L_3\in L^1(0,T)$ and $L_4\in L^2(0,T)$.

Put
\begin{equation}\label{1.18-eq1}
\cX(0,T)\triangleq\big[L^2_\dbF(\Om;C([0,T]; \wh
H))\times L^2_\dbF(0,T;\wh H)\big]\cap
\big[\dbL^a_{1,2}(\wh H)\times \dbL^a_{1,2}( \wh
H)\big].
\end{equation}
\begin{proposition}\label{prop2}\cite[Proposition 3.2]{Dou}
Let $z_T\in \dbD^{1,2}(\wh H)$. Then $(z,Z)\in
\cX(0,T)$ satisfying\vspace{1mm}
\begin{equation}\label{prop2-eq1}
|(z,Z)|_{\cX(0,T)}\leq
\cC\big(|z_T|_{\dbD^{1,2}( \wh
H)}+|F(\cd,0,0)|_{\dbL^a_{1,2}(\wh H)}\big).
\end{equation}
Further,  a version of $(\mathbf{D}_\tau
z(\cd),\mathbf{D}_\tau Z(\cd))$
solves\vspace{1mm}
\begin{equation}\label{system4}
\!\left\{\2n
\begin{array}{ll}
\ds d\mathbf{D}_\tau z +  A ^*\mathbf{D}_\tau z
dt=\big[F_z(t,z,Z)\mathbf{D}_\tau z +
F_Z(t,z,Z)\mathbf{D}_\tau Z\\
\ns\ds\hspace{3.4cm} + \mathbf{D}_\tau
F(t,z,Z)\big] dt + \mathbf{D}_\tau Z
dW(t)\;\mbox{ in }[\tau,T), \\
\ns\ds \mathbf{D}_\tau z(T)=\mathbf{D}_\tau z_T.
\end{array}
\right.
\end{equation}
\end{proposition}
%


\section{Proof of the first main result}\label{proof}


In this chapter, we give the proof of the first
main result, i.e., Theorem \ref{5.7-th1.1}.

\vspace{0.2cm}

{\it Proof of Theorem \ref{5.7-th1.1}}\,: Let us
assume that the equation \eqref{5.5-eq6} admits
a transposition solution\vspace{0.1cm}
$$\big(P(\cd),\Lambda(\cd)\big)\in
C_{\dbF,w}([0,T]; L^{\infty}(\Om;\cL(H))) \times
L^2_{\dbF,w}(0,T;\cL(H))$$\vspace{0.1cm} such
that
$$
K(\cd)^{-1}\big[B(\cd)^* P(\cd) +D(\cd)^*
P(\cd)C(\cd) + D(\cd)^*\Lambda(\cd)\big] \in
\Upsilon_2(H;U)\cap \Upsilon_2(V';\wt U).
$$
Then,\vspace{0.1cm}
\begin{equation}\label{6.19-eq5}
\Th \triangleq-K^{-1}(B^* P+D^* PC+
D^*\Lambda)\in \Upsilon_2(H;U)\cap
\Upsilon_2(V';\wt U).
\end{equation}

\ms

For any $s\in [0,T)$, $\eta\in
L^2_{\cF_s}(\Om;H)$ and $u(\cd)\in
L^2_\dbF(s,T;U)$, let $x(\cd)\equiv
x(\cd\,;s,\eta,$ $u(\cd))$ be the corresponding
state process of the system \eqref{5.2-eq1}.
Choose
$$
\begin{cases}\ds
\xi_1=\xi_2=\eta,\\
\ns\ds u_1=u_2=Bu,\\
\ns\ds v_1=v_2=Du
\end{cases}
$$
in the equations
\eqref{op-fsystem1}--\eqref{op-fsystem2}.

\ms

From the definitions of $K$ and $L$ (see
\eqref{9.7-eq10}), and  \eqref{10.10-eq10} in
the definition of the transposition solution to
\eqref{5.5-eq6}, and the pointwise symmetry of
$K(\cd)$, we obtain that\vspace{0.1cm}
\begin{eqnarray}\label{6.8-eq19}
&& \mE\langle G x(T),x(T)\rangle_H + \mE
\int_s^T \big\langle Q(r) x(r), x(r)
\big\rangle_{H}dr\nonumber  \\
&& \q- \mE \int_s^T \big\langle \Th(r)^*
K(r)\Th(r) x(r), x(r) \big\rangle_{H}dr
\nonumber  \\[1mm]
&&  =\mE\big\langle P(s) \eta,\eta
\big\rangle_{H}+\mE\int_s^T\big\langle
P(r)B(r)u(r), x(r)\big\rangle_{H}dr  \\[1mm]
&& \q +\mE \int_s^T\big\langle P(r)x(r),
B(r)u(r)\big\rangle_{H}dr  +\mE\int_s^T
\big\langle P(r)C(r)x(r),
D(r)u(r)\big\rangle_{H}dr \nonumber\\[1mm]
&&  \q  + \mE
\int_s^T\big\langle  P(r)D(r)u(r), C(r)x(r)+D(r)u(r)\big\rangle_{H}dr\nonumber\\[1mm]
&&  \q + \mE \int_s^T \big\langle u(r),
D(r)^*\Lambda(r)x(r)\big\rangle_Udr+ \mE
\int_s^T \big\langle D(r)^*\Lambda(r)x(r), u(r)
\big\rangle_Udr.\nonumber \vspace{0.1cm}
\end{eqnarray}
By \eqref{6.8-eq19}, and recalling the
definition of the cost functional
$\cJ(s,\eta;u(\cd))$ in \eqref{5.2-eq2},  we
arrive at\vspace{0.1cm}
\begin{eqnarray*}
&&\3n\3n 2\cJ(s,\eta;u(\cd))\\[1mm]
&&\3n\3n= \dbE\Big[\int_s^T\Big(\big\langle
Qx(r),x(r)\big\rangle_H
+\big\langle Ru(r),u(r)\big\rangle_U \Big)dr+\big\langle Gx(T),x(T)\big\rangle_H\Big]\\[1mm]
&&\3n\3n=\mE\big\langle P(s) \eta,\eta
\big\rangle_{H}+\mE\int_s^T\big\langle PBu(r),
x(r)\big\rangle_{H}dr+\mE \int_s^T\big\langle
Px(r),
Bu(r)\big\rangle_{H}dr \\[1mm]
&&\3n\3n \q  +\mE\int_s^T \big\langle PCx(r),
Du(r)\big\rangle_{H}dr + \mE
\int_s^T\big\langle  PDu(r), Cx(r)+Du(r)\big\rangle_{H}dr\\
&&\3n\3n \q + \mE \int_s^T \big\langle u(r),
D^*\Lambda(r)x(r)\big\rangle_Udr+ \mE \int_s^T
\big\langle D^*\Lambda(r)x(r), u(r)
\big\rangle_Udr\\[1mm]
&& \3n\3n\q  + \mE \int_s^T \big\langle \Th^* K
\Th x(r), x(r) \big\rangle_{H}dr  + \mE\int_s^T
\big\langle Ru(r),u(r)\big\rangle_U dr\\[1mm]
&& \3n\3n= \dbE\[\big\langle
P(s)\eta,\eta\big\rangle_H + \int_s^T \(
\big\langle \Th^* K\Th x(r),x(r)\big\rangle_H
\!+\! 2\big\langle [B^* P + D^*
(PC\!+\!\Lambda)]x(r),
\\[1mm]
&&\3n\3n\hspace{4cm}u(r)\big\rangle_U
+\big\langle (R+D^*PD)
u(r),u(r)\big\rangle_U\)dr\]
\\[1mm]
&&\3n\3n= \dbE\[\big\langle
P(s)\eta,\eta\big\rangle_H + \int_s^T \(
\big\langle \Th^* K\Th x(r),x(r)\big\rangle_H +
2\big\langle Lx(r),u(r)\big\rangle_U
\\[1mm]
&&\3n\3n\hspace{4cm}+\big\langle K
u(r),u(r)\big\rangle_U\)dr\].\vspace{0.1cm}
\end{eqnarray*}
This, together with the definition of $\Th$ in
\eqref{6.19-eq5}, implies that \vspace{0.1cm}
\begin{equation}\label{5.31-eq1}
\begin{array}{ll}\ds
\cJ(s,\eta;u(\cd))  \\[1mm]
\ns\ds=\frac{1}{2}\dbE\Big[\big\langle
P(s)\eta,\eta\big\rangle_H+\int_s^T\big(\big\langle
K\Th x,\Th x\big\rangle_U-2\big\langle K\Th x,u\big\rangle_U \\[1mm]
\ns\ds\hspace{4.4cm}+\big\langle Ku,u\big\rangle_U\big)dr\Big]\\[1mm]
\ns\ds=\frac{1}{2}\dbE\Big(\big\langle
P(s)\eta,\eta\big\rangle_H+\int_s^T\big\langle
K(u-\Th x),u-\Th
x\big\rangle_Udr\Big).\vspace{0.1cm}
\end{array}
\end{equation}
By taking $u=\Th x$, that is, the control is
chosen as the feedback form given by $\Th$, we
find that\vspace{0.1cm}
$$
\cJ(s,\eta;(\Th \bar x)(\cd)) =\frac{1}{2}\dbE
\big\langle
P(s)\eta,\eta\big\rangle_H.\vspace{0.1cm}
$$
This, together with \eqref{5.31-eq1}, implies
that\vspace{0.1cm}
$$
\begin{array}{ll}\ds
\cJ(s,\eta;u(\cd))\\
\ns\ds=\cJ\big(s,\eta;\Th \bar x\big)
+\frac{1}{2}\dbE\int_s^T\big\langle K(u-\Th
x),u-\Th x\big\rangle_U dr.
\end{array}
$$
Consequently,
$$
\cJ(s,\eta;\Th \bar x)\leq
\cJ(s,\eta;u),\q\forall\, u(\cd)\in
L^2_\dbF(s,T;U).\vspace{0.1cm}
$$
Therefore, $\Th(\cd)$ is an optimal feedback
operator for {\bf Problem (SLQ)}, and
\eqref{Value} holds. This completes the proof of
Theorem \ref{5.7-th1.1}.


\section{Proof of the second main result}\label{proof-n}


This chapter is devoted to the proof of the
second main result, i.e., Theorem \ref{5.7-th1}.
We borrow some ideas from \cite{Reid1, Kalman,
Bismut2, AMZ, Tang1}.

\ms

{\it Proof of Theorem \ref{5.7-th1}}\,: Without
loss of generality, we assume that $s=0$, and
$H$ and $U$ are real separable Hilbert spaces.
The proof is rather long, and therefore we
divide it into several steps.

\medskip

{\bf Step 1}.  In this step, we introduce some
operators $X(\cd)$\;\footnote{In the sequel, we
shall interchangeably use $X(t)$, $X(t,\cd)$, or
even $X$ to denote the operator $X(\cd)$. The
same can be said for $Y(\cd)$, $Z(\cd)$ and $\wt
X(\cd)$.}, $Y(\cd)$, $Z(\cd)$ and $\wt X(\cd)$.

\vspace{0.1cm}

Let $\Th(\cd)\in \Upsilon_2(H;U)\cap
\Upsilon_2(V';\wt U)$ be an optimal feedback
operator of {\bf Problem (SLQ)} on $[0,T]$.
Then, by Corollary \ref{5.7-prop1}, for any
$\zeta\in H$, the following forward-backward
stochastic evolution equations
\begin{equation}\label{5.7-eq4.1}
\left\{
\begin{array}{ll}
\ds d\hat x(t)= \big[(A+A_1)+B\Th\big]\hat
x(t)dt+ (C+D\Th)\hat x(t)
dW(t)& \mbox{ in }(0,T],\\[1mm]
\ns\ds
dy(t)=-\big[(A+A_1)^* y(t)+C^* z(t)+Q\hat x(t)\big]dt+z(t)dW(t) & \mbox{ in }[0,T),\\[1mm]
\ns\ds \hat x(0)=\zeta,\q y(T)= G\hat x(T)
\end{array}
\right.
\end{equation}
admits a unique mild solution $$(\hat
x(\cd),y(\cd),z(\cd))\in
L^2_\dbF(\Om;C([0,T];H))\times
L^2_\dbF(\Om;C([0,T];H))\times L^2_\dbF(0,T;
H)$$ such that
\begin{equation}\label{5.7-eq3}
R\Th \hat x+B^* y+D^* z =0,\q\mbox{ a.e. }
(t,\om)\in (0,T)\times\Om.
\end{equation}

Further, consider the following SEE:
\begin{equation}\label{5.26-eq4}
\left\{
\begin{array}{ll}
\ns\ds d\tilde x(t) =
\big[-A-A_1-B\Th+\big(C+D\Th\big)^2
\big]^*\tilde x(t) dt\\[1mm]
\ns\ds\qq\qq - \big(C+D\Th\big)^*\tilde
x (t)dW(t) &  \mbox{ in }\;(0,T],\\[1mm]
\ns\ds \tilde x(0)=\zeta.
\end{array}
\right.
\end{equation}
Note that $A$ generates a $C_0$-group, and
hence, so does $-A^*$. By Lemma \ref{lm2}, the
equation \eqref{5.26-eq4} admits a unique mild
solution $\tilde x(\cd)\in
L^2_\dbF(\Om;C([0,T];H))$.

\vspace{0.1cm}

For each $n\in\dbN$, denote by $\wt \G_n$ the
projection operator from $U$ to
$U_n\triangleq\span_{1\leq j\leq n}\{\f_j\}$
(Recall that $\{\f_j\}_{j=1}^\infty$ is an
orthonormal basis of $U$). Write (Recall Chapter
\ref{pre} for $ \G_n$)
$$
\begin{array}{ll}\ds
A_{1,n}=\G_nA_1\G_n,\q B_n = \G_n B \wt \G_n,\\[1mm]
\ns\ds C_n
= \G_n C \G_n, \qq D_n = \G_n D \wt\G_n,\\[1mm]
\ns\ds  Q_n = \G_n Q \G_n,\;\;\;\q R_n = \wt\G_n
R \wt\G_n,\; \q \Th_n = \wt\G_n \Th \G_n.
\end{array}
$$
It is easy to show that
\begin{equation}\label{6.8-eq5}
\begin{cases} \ds
\lim_{n\to+\infty} A_{1,n}\zeta = A_1 \zeta \q\mbox{ in } H,  \\
\ns\ds
\lim_{n\to+\infty} C_n\zeta = C \zeta \q\mbox{ in } H,\\
\ns\ds \lim_{n\to+\infty} Q_n\zeta = Q \zeta
\q\mbox{ in } H,
\\
\ns\ds
\lim_{n\to+\infty} \Th_n\zeta = \Th \zeta \q\mbox{ in } H,\\
\ns\ds \hspace{0.7cm}\mbox{ for all }\zeta\in H
\mbox{ and a.e. } (t,\om)\in [0,T]\times\Om,
\end{cases}
\end{equation}
\begin{equation}\label{6.8-eq5.1}
\begin{cases} \ds
\lim_{n\to+\infty} B_n\varsigma = B \varsigma
\q\mbox{ in } H,
\\[1mm]
\ns\ds \lim_{n\to+\infty} D_n\varsigma = D
\varsigma \q\mbox{ in } H,\\[1mm]
\ns\ds \lim_{n\to+\infty} R_n\varsigma= R \varsigma \q\mbox{ in } U, \\[1mm]
\ns\ds \hspace{0.6cm}\mbox{ for all
}\varsigma\in U \mbox{ and a.e. } (t,\om)\in
[0,T]\times\Om.
\end{cases}
\end{equation}

Consider the following forward-backward
stochastic differential equation:
\begin{equation}\label{6.7-eq1}
\left\{
\begin{array}{ll}
\ds d\hat x_n=  (A_n+A_{1,n}+B_n\Th_n) \hat x_n
dt+ (C_n+D_n\Th_n) \hat x_n
dW(t)& \mbox{ in }[0,T],\\[1mm]
\ns\ds dy_n=-\big[(A_n+A_{1,n})^* y_n+ C^*_n
z_n+ Q_n\hat x_n\big]dt+ z_ndW(t) &
\mbox{ in }[0,T],\\[1mm]
\ns\ds \hat x_n(0)=\G_n\zeta, \q y_n(T)= G_n
\hat x_n(T)
\end{array}
\right.
\end{equation}
and the following stochastic differential
equation
\begin{equation}\label{6.7-eq2}
\left\{
\begin{array}{ll}
\ds d\tilde x_n= \big[ -A_n - A_{1,n} - B_n\Th_n
+\big(C_n + D_n\Th_n\big)^2 \big]^*\tilde x_n
dt\\[1mm]
\ns\ds\qq\q  - (C_n+D_n\Th_n)^*\tilde x_n
dW(t)& \mbox{in }[0,T],\\[1mm]
\ns\ds \tilde x_n(0)=\G_n\zeta,
\end{array}
\right.
\end{equation}
where $A_n$ and $G_n$ are given in
\eqref{20161214e1}. For each $t\in [0,T]$,
define three operators $X_{n,t}$, $Y_{n,t}$ and
$\wt X_{n,t}$ on $H_n$ as follows:
\begin{equation}\label{7.19-eq5}
\begin{cases}\ds
X_{n,t}\G_n\zeta\triangleq \hat
x_n(t;\G_n\zeta), \\[1mm]
\ns\ds Y_{n,t}\G_n\zeta\triangleq
y_n(t;\G_n\zeta), \\[1mm]
\ns\ds \wt X_{n,t}\G_n\zeta\triangleq \tilde
x_n(t;\G_n\zeta),
\end{cases}
\q\forall\,\zeta\in H.
\end{equation}
For a.e. $t\in [0,T]$, define an operator
$Z_{n,t}$ on $H_n$ by
\begin{equation}\label{7.19-eq5zz}
Z_{n,t}\G_n\zeta\triangleq  z_n(t;\G_n\zeta),
\qq \forall\,\zeta\in H.
\end{equation}
By the well-posedness results for the equations
\eqref{6.7-eq1} and \eqref{6.7-eq2}, and the
fact that both $A$ and $-A^*$  generate
$C_0$-semigroups on $H$ (because $A$ generates a
$C_0$-group on $H$), we see that
$$
\begin{cases} \ds
|X_{n,t}\G_n\zeta|_{L^2_{\cF_t}(\Om;H)}\leq
\cC|\zeta|_H,\\[1mm]
\ns\ds
|Y_{n,t}\G_n\zeta|_{L^2_{\cF_t}(\Om;H)}\leq
\cC|\zeta|_H,\\[1mm]
\ns\ds
|Z_{n,\cd}\G_n\zeta|_{L^2_{\dbF}(0,T;H)}\leq
\cC|\zeta|_H,\\[1mm]
\ns\ds |\wt
X_{n,t}\G_n\zeta|_{L^2_{\cF_t}(\Om;H)}\leq
\cC|\zeta|_H,
\end{cases}
$$
where the constant $\cC$ is independent of $n$.
This implies that
\begin{equation}\label{7.2016-eq2}
\begin{cases} \ds
|X_{n,t}\G_n|_{\cL(H;L^2_{\cF_t}(\Om;H))}\leq
\cC,\\[1mm]
\ns\ds
|Y_{n,t}\G_n|_{\cL(H;L^2_{\cF_t}(\Om;H))}\leq
\cC,\\[1mm]
\ns\ds
|Z_{n,\cd}\G_n|_{\cL(H;L^2_{\dbF}(0,T;H))}\leq
\cC,\\[1mm]
\ns\ds |\wt
X_{n,t}\G_n|_{\cL(H;L^2_{\cF_t}(\Om;H))}\leq
\cC.
\end{cases}
\end{equation}

Denote by $I_n$ the identity matrix on $\dbR^n$
(or, the identity map on $H_n$). Consider the
following equations:
\begin{equation}\label{7.19-eq1}
\left\{
\begin{array}{ll}
\ds dX_n=  (A_n+A_{1,n}+B_n\Th_n)X_n dt \\[1mm]
\ns\ds\qq\q  + (C_n+D_n\Th_n) X_n
dW(t)& \mbox{ in }[0,T],\\[1mm]
\ns\ds dY_n=-\big[(A_n+A_{1,n})^* Y_n+ C^*_n
Z_n+ Q_nX_n\big]dt \\[1mm]
\ns\ds\qq\q+ Z_ndW(t) &
\mbox{ in }[0,T],\\[1mm]
\ns\ds X_n(0)=I_{n}, \q Y_n(T)= G_nX_n(T)
\end{array}
\right.
\end{equation}
and
\begin{equation}\label{7.19-eq2}
\left\{
\begin{array}{ll}
\ds d\wt X_n = \big[ - A_n - A_{1,n} - B_n\Th_n
+ \big(C_n + D_n\Th_n\big)^2 \big]^*\wt X_n dt
\\[1mm]
\ns\ds\qq\q - (C_n + D_n\Th_n)^*\wt X_n
dW(t)&\mbox{ in }[0,T],\\[1mm]
\ns\ds \wt X_n(0)=I_n.
\end{array}
\right.
\end{equation}
Clearly, both \eqref{7.19-eq1} and
\eqref{7.19-eq2} can be viewed as $\dbR^{n\times
n}\equiv\dbR^{n^2}$-valued equations.

By Lemmas \ref{lm2}--\ref{lm3}, the equations
\eqref{7.19-eq1} and \eqref{7.19-eq2} admit
unique solutions
$$ (X_n,Y_n,Z_n)\in L^2_\dbF(\Om;C([0,T]; \dbR^{n\times n}))\times
L^2_\dbF(\Om;C([0,T];\dbR^{n\times n}))\times
L^2_\dbF(0,T;\dbR^{n\times n}) $$
and
$$ \wt X_n \in L^2_\dbF(\Om;C([0,T];\dbR^{n\times
n})),$$
respectively. It follows from
\eqref{7.19-eq5}--\eqref{7.19-eq2} that, for
a.e. $t\in[0,T]$,
\begin{equation}\label{7.19-eq6}
\begin{cases} \ds
X_{n,t}\G_n\zeta =  X_{n}(t)\G_n\zeta, \\[1mm]
\ns\ds Y_{n,t}\G_n\zeta =  Y_{n}(t)\G_n\zeta, \\[1mm]
\ns\ds Z_{n,t}\G_n\zeta =  Z_{n}(t)\G_n\zeta,
\\[1mm]
\ns\ds \wt X_{n,t}\G_n\zeta = \wt
X_{n}(t)\G_n\zeta,
\end{cases}
\qq \forall\,\zeta\in H.
\end{equation}
Thus, $X_{n,t}\G_n$, $Y_{n,t}\G_n$ and $\wt
X_{n,t}\G_n$ belong to $L^2_{\cF_t}(\Om;\cL(H))$
and $Z_{n,\cd}\G_n\in L^2_\dbF(0,T;$ $\cL(H))$.
It is easy to see that $$X_{n,t}\G_n,
Y_{n,t}\G_n, \wt X_{n,t}\G_n\in
\cL_{pd}(H;L^2_{\cF_t}(\Om;H)) \mbox{ for each }
t\in [0,T]$$ and $$Z_{n,\cd}\G_n\in
\cL_{pd}(H;L^2_{\dbF}(0,T;H)).$$

\ms

By \eqref{7.2016-eq2} and using  Theorems 5.2
and 5.3 in \cite{LZ1}, we deduce that, there
exist suitable subsequences
$\{X_{n_k,t}\}_{k=1}^\infty\subset
\{X_{n,t}\}_{n=1}^\infty$,
$\{Y_{n_k,t}\}_{k=1}^\infty$ $\subset
\{Y_{n,t}\}_{n=1}^\infty$,
$\{Z_{n_k,t}\}_{k=1}^\infty\subset
\{Z_{n,t}\}_{n=1}^\infty$ and $\{\wt
X_{n_k,t}\}_{k=1}^\infty\subset \{\wt
X_{n,t}\}_{n=1}^\infty$ (these sequences may
depend on $t$), and (pointwise defined)
operators\vspace{0.2mm}
$$
\begin{array}{ll}\ds
\begin{cases}\ds
X(t,\cd) \in \cL_{pd}(H;
L^2_{\cF_t}(\Om;H)),\\
\ns\ds  Y(t,\cd)  \in
\cL_{pd}(H; L^2_{\cF_t}(\Om;H)),\\
\ns\ds  \wt X(t,\cd) \in \cL_{pd}(H;
L^2_{\cF_t}(\Om;H))
\end{cases}\mbox{ for each }t\in
[0,T]
\end{array}
$$
and $$Z(\cd,\cd)\in \cL_{pd}(H;L^2_{\dbF}(0,T;
H))$$  such that
\begin{equation}\label{7.24-eq3}
\left\{
\begin{array}{ll}\ds
\lim_{k\to+\infty}X_{n_k,t}\G_{n_k}\zeta =
X(t,\cd)\zeta &\mbox{ weakly in }L^2_{\cF_t}(\Om;H),\\
\ns\ds \lim_{k\to+\infty}Y_{n_k,t}\G_{n_k}\zeta
=
Y(t,\cd)\zeta &\mbox{ weakly in }L^2_{\cF_t}(\Om;H),\\
\ns\ds\lim_{k\to+\infty}Z_{n_k,t}\G_{n_k}\zeta =
Z(\cd,\cd)\zeta &\mbox{ weakly in }L^2_{\dbF}(0,T;H),\\
\ns\ds \lim_{k\to+\infty}\wt
X_{n_k,t}\G_{n_k}\zeta = \wt X(t,\cd)\zeta
&\mbox{ weakly in }L^2_{\cF_t}(\Om;H),
\end{array}
\right.
\end{equation}
and that
\begin{equation}\label{7.24-eq6}
\begin{cases} \ds
|X(t,\cd)\zeta|_{L^2_{\cF_t}(\Om;H)}\leq
\cC|\zeta|_H,\\
\ns\ds |Y(t,\cd)\zeta|_{L^2_{\cF_t}(\Om;H)}\leq
\cC|\zeta|_H,\\
\ns\ds |Z(\cd,\cd)\zeta|_{L^2_{\dbF}(0,T;H)}\leq
\cC|\zeta|_H,\\
\ns\ds |\wt
X(t,\cd)\zeta|_{L^2_{\cF_t}(\Om;H)}\leq
\cC|\zeta|_H.
\end{cases}
\end{equation}

On the other hand, from the definitions of $\hat
x_n(\cd;\G_n\zeta)$ and $\tilde
x_n(\cd;\G_n\zeta)$, by Lemma \ref{lm7}, we have
that
\begin{equation}\label{6.7-eq3}
\left\{
\begin{array}{ll}\ds
\lim_{n\to+\infty} \hat x_n(\cd;\G_n\zeta) =
\hat x(\cd;\zeta) &\mbox{ in }
L^2_\dbF(\Om;C([0,T];H)),\\
\ns\ds \lim_{n\to+\infty} y_n(\cd;\G_n\zeta) =
y(\cd;\zeta) &\mbox{ in }
L^2_\dbF(\Om;C([0,T];H)), \\
\ns\ds \lim_{n\to+\infty} z_n(\cd;\G_n\zeta) =
z(\cd;\zeta) &\mbox{ in } L^2_\dbF(0,T;H),
\\
\ns\ds \lim_{n\to+\infty} \tilde
x_n(\cd;\G_n\zeta) =\tilde x(\cd;\zeta) &\mbox{
in } L^2_\dbF(\Om;C([0,T];H)).
\end{array}
\right.
\end{equation}
Hence, in view of \eqref{7.19-eq5}, we find that
\begin{equation}\label{7.24-eq4}
\left\{
\begin{array}{ll}\ds
\lim_{n\to+\infty}X_{n,t}\G_{n}\zeta =
\hat x(t;\zeta) &\mbox{ strongly in }L^2_{\cF_t}(\Om;H),\\[1mm]
\ns\ds \lim_{n\to+\infty}Y_{n,t}\G_{n}\zeta =
y(t;\zeta) &\mbox{ strongly in }L^2_{\cF_t}(\Om;H),\\[1mm]
\ns\ds\lim_{n\to+\infty}Z_{n,t}\G_{n}\zeta =
z(t;\zeta) &\mbox{ strongly in }L^2_{\dbF}(0,T;H),\\[1mm]
\ns\ds \lim_{n\to+\infty}\wt X_{n,t}\G_{n}\zeta
= \tilde x(t;\zeta) &\mbox{ strongly in
}L^2_{\cF_t}(\Om;H).
\end{array}
\right.
\end{equation}
According to \eqref{7.24-eq3} and
\eqref{7.24-eq4}, we obtain that
\begin{equation}\label{7.24-eq5}
\begin{cases}
X(t,\cd)\zeta=\hat  x(t;\zeta),\\
\ns\ds Y(t,\cd)\zeta= y(t;\zeta),\\
\ns\ds Z(t,\cd)\zeta = z(t;\zeta),\\
\ns\ds\wt X(t,\cd) \zeta = \tilde x(t;\zeta).
\end{cases}
\end{equation}
Also, from the equality \eqref{5.7-eq3} and
noting \eqref{7.24-eq5}, we find
that\vspace{0.2mm}
\begin{equation}\label{5.7-eq5.1}
R\Th X+B^* Y+D^* Z =0, \q \mbox{for a.e. }
(t,\om)\in[0,T]\times\Om.
\end{equation}

Combining \eqref{6.7-eq3} and \eqref{7.24-eq5},
we find that\vspace{0.2mm}
\begin{equation}\label{2.10-eq1}
\left\{
\begin{array}{ll}\ds
\lim_{n\to+\infty}X_{n,t}\G_{n}\zeta =
X(t,\cd)\zeta &\mbox{ strongly in }L^2_{\cF_t}(\Om;H),\\[1mm]
\ns\ds \lim_{n\to+\infty}Y_{n,t}\G_{n}\zeta =
Y(t,\cd)\zeta &\mbox{ strongly in }L^2_{\cF_t}(\Om;H),\\[1mm]
\ns\ds\lim_{n\to+\infty}Z_{n,t}\G_{n}\zeta =
Z(\cd,\cd)\zeta &\mbox{ strongly in }L^2_{\dbF}(0,T;H),\\[1mm]
\ns\ds \lim_{n\to+\infty}\wt X_{n,t}\G_{n}\zeta
= \wt X(t,\cd)\zeta &\mbox{ strongly in
}L^2_{\cF_t}(\Om;H).
\end{array}
\right.
\end{equation}
Moreover, from Lemma \ref{lm8}, it follows that
\begin{equation}\label{10.10-eq21}
\begin{cases}
X, Y, \wt X\in
L^4_\dbF(\Om;C([0,T];\cL_2(H;V))), \\[1mm]
\ns\ds Z\in L^4_\dbF(\Om;L^2(0,T;\cL_2(H;V)))
\end{cases}
\end{equation}
and
\begin{equation}\label{10.10-eq22}
\begin{cases}\ds
\lim_{n\to\infty}X_n = X &\mbox{ in
}L^4_\dbF(\Om;C([0,T];\cL_2(H;V))),\\
\ns\ds \lim_{n\to\infty}Y_n = Y &\mbox{ in
}L^4_\dbF(\Om;C([0,T];\cL_2(H;V))),\\
\ns\ds \lim_{n\to\infty}Z_n = Z &\mbox{ in
}L^4_\dbF(\Om;L^2(0,T;\cL_2(H;V))),\\
\ns\ds \lim_{n\to\infty}\wt X_n = \wt X &\mbox{
in }L^4_\dbF(\Om;C([0,T];\cL_2(H;V))).
\end{cases}
\end{equation}

\medskip

\ms

{\bf Step 2}. Denote by $I$ the identity
operator on $H$. In this step, we shall prove
that $X(\cd)\wt X(\cd)^*=I$, for a.e.
$(t,\om)\in(0,T)\times\Om$.

\vspace{0.1cm}

For any $\zeta,\rho\in H$ and $t\in [0,T]$, by
It\^o's formula, we have\vspace{0.1cm}
$$
\begin{array}{ll}\ds\big\langle
\hat x_n(t;\G_n\zeta),\tilde
x_n(t;\G_n\rho)\big\rangle_{H_n} -
\big\langle \G_n\zeta,\G_n\rho\big\rangle_{H_n}\\[1mm]
\ns\ds =\int_0^t \big\langle  \big(A_n + A_{1,n}
+B_n \Th_n\big) \hat x_n(r;\G_n\zeta),\tilde
x_n(r;\G_n\rho)
\big\rangle_{H_n} d\tau\\[1mm]
\ns\ds \q + \int_0^t \big\langle \big(C_n +D_n
\Th_n \big)
\hat x_n(r;\G_n\zeta),\tilde x_n(r;\G_n\rho) \big\rangle_{H_n} dW(\tau)\\[1mm]
\ns\ds \q + \int_0^t \big\langle \hat
x_n(r;\G_n\zeta), \big[-A_n\! -\!A_{1,n}\!-\!B_n
\Th_n \!+\!\big(C_n+D_n\Th_n\big)^2\big]^*
\tilde x_n(r;\G_n\rho)
\big\rangle_{H_n} d\tau\\[1mm]
\ns\ds \q - \int_0^t \big\langle \hat
x_n(r;\G_n\zeta), \big(C_n+D_n \Th_n \big)^*
\tilde
x_n(r;\G_n\rho) \big\rangle_{H_n} dW(\tau)\\[1mm]
\ns\ds \q - \int_0^t \big\langle \big(C_n +D_n
\Th_n \big) \hat x_n(r;\G_n\zeta), \big(C_n +D_n
\Th_n \big)^* \tilde x_n(r;\G_n\rho)
\big\rangle_{H_n} d\tau=0.
\end{array}
$$
Hence,
$$
\begin{array}{ll}\ds
\big\langle X_{n,t}\G_n\zeta, \wt
X_{n,t}\G_n\rho\big\rangle_{H_n}\3n&\ds=\big\langle
\hat x_n(t;\G_n\zeta),\tilde
x_n(t;\G_n\rho)\big\rangle_{H_n}\\[1mm]
\ns&\ds = \big\langle
\G_n\zeta,\G_n\rho\big\rangle_{H_n}, \q
\dbP\mbox{-}\as
\end{array}
$$
This implies that $$X_{n,t}\wt X_{n,t}^*=I_n,\q
\dbP\mbox{-a.s.}$$  Namely, $\wt
X_{n,t}^*=X_{n,t}^{-1}$, $\dbP$-a.s.

\vspace{0.1cm}

By \eqref{6.7-eq3} and \eqref{7.24-eq5},  for
any  $\zeta\in H$, it is easy to see that
$ \wt X_{n,t}(\cd)^*\zeta$ converges weakly to
$\wt X(t,\cd)^*\zeta$  in $L^2_{\cF_t}(\Om;H)$
as  $n\to\infty$.
On the other hand, by Lemma \ref{lm8}, we have
that
$$
\ds\lim_{n\to\infty}|X_{n,t}(\cd)-X(t,\cd)|_{L^4_{\cF_t}(\Om;\cL_2(H;V))}=0.
$$
Thus,
$$
\begin{array}{ll}\ds
X_{n_k,t}(\cd)\wt X_{n_k,t}(\cd)^*\zeta \mbox{
converges weakly to }\\[1mm]
\ns\ds X(t,\cd)\wt X(t,\cd)^*\zeta \;\mbox{ in
}\;L^2_{\cF_t}(\Om;V),\q \mbox{ as }k\to\infty.
\end{array}
$$
This yields that $$ X(t,\cd)\wt X(t,\cd)^*\zeta
= \zeta \q\mbox{  in  } V, \q\dbP\mbox{-a.s.}$$
Further, by $X(t,\cd),\wt X(t,\cd)\in \cL(H)$,
$\dbP$-a.s., we deduce that
$$
X(t,\cd)\wt X(t,\cd)^*\zeta = \zeta \q\mbox{ in
} H, \q \dbP\mbox{-a.s.},
$$
which implies that $$X(\cd)\wt X(\cd)^*=I
\q\mbox{ for a.e. } (t,\om)\in [0,T]\times\Om.$$

\vspace{0.1cm}

Put
\begin{equation}\label{6.8-eq13.1}
\left\{
\begin{array}{ll}\ds
P(t,\cd) = Y(t,\cd)\wt X(t,\cd)^*,\qq\forall\;t\in [0,T],\\[1mm]
\ns\ds \Pi(\cd) = Z(\cd)\wt X(\cd)^*,\\[1mm]
\ns\ds
\Lambda(\cd)=\Pi(\cd)-P(\cd)\big(C(\cd)+D(\cd)\Th(\cd)\big).
\end{array}
\right.
\end{equation}
From \eqref{10.10-eq21}, it follows that $$ \wt
X(\cd)^*\in L^4_\dbF(\Om;C([0,T];\cL_2(V';H))).
$$ This, together with \eqref{6.8-eq13.1},
implies that
\begin{equation}\label{10.10-eq23}
\begin{cases}
P(\cd)\in L^2_\dbF(\Om;C([0,T];\cL_2(V';V))),\\[1mm]
\ns\ds \Lambda(\cd) \in
L^2_\dbF(0,T;\cL_2(V';V)).
\end{cases}
\end{equation}

\medskip

{\bf Step 3}. In this step, we construct a
sequence of finite dimensional approximations of
$(P,\Lambda)$, which will serve as approximate
transposition solutions to the equation
\eqref{5.5-eq6}.

\vspace{0.1cm}

Put
\begin{equation}\label{8.25-eq3}
\begin{cases}
P_n= Y_{n} \wt X_{n}^{*},\\[1mm]
\ns\ds\Pi_n= Z_n \wt X_{n}^{*},\\
\ns\ds \Lambda_n=\Pi_n-P_n (C_n+D_n\Th_n).
\end{cases}
\end{equation}
It follows from Lemma \ref{lm8}
that\vspace{0.2mm}
\begin{equation}\label{6.8-12eq4}
\ba{ll}\ds \lim_{n\to \infty}  |\wt X_{n}^* -
\wt X^*|_{L^4_\dbF(\Om;C([0,T];\cL_2(V';H)))}\\[1mm]
\ns\ds =\lim_{n\to \infty}  |\wt X_{n} - \wt
X|_{L^4_\dbF(\Om;C([0,T];\cL_2(H;V)))}=0, \ea
\end{equation}
\begin{equation}\label{6.8-12eq4.1}
\lim_{n\to \infty}  | Y_{n} -
Y|_{L^4_\dbF(\Om;C([0,T];\cL_2(H;V)))}=0,
\end{equation}
and
\begin{equation}\label{6.8-12eq4.2}
\lim_{n\to \infty}  | Z_{n} -
Z|_{L^4_\dbF(\Om;L^2(0,T;\cL_2(H;V)))}=0.
\end{equation}
Now, by \eqref{6.8-eq13.1}, \eqref{8.25-eq3},
\eqref{6.8-12eq4}, \eqref{6.8-12eq4.1} and
\eqref{6.8-12eq4.2}, we deduce that
\begin{equation}\label{6.8-eq4}
\lim_{n\to \infty}  |P_{n} -
P|_{L^2_\dbF(\Om;C([0,T];\cL_2(V';V)))}=0,
\end{equation}
and
\begin{equation}\label{8.25-eq5}
\lim_{n\to \infty}  |\Lambda_{n} -
\Lambda|_{L^2_\dbF(0,T;\cL_2(V';V))}=0.
\end{equation}

By It\^o's formula, and noting
\eqref{7.19-eq1}--\eqref{7.19-eq2}, we obtain
that
$$
\3n\begin{array}{ll} \ds
dP_n\3n&\ds=\Big\{-\big[(A_n+ A_{1,n})^* Y_n+
C^*_n Z_n +  Q_n  X_n\big] X^{-1}_n \\[1mm]
\ns&\ds\q +  Y_n X^{-1}_n \big[(C_n +D_n\Th_n)^2
- A_n-A_{1,n}-B_n\Th_n\big]\\[1mm]
\ns&\ds\q - Z_n X^{-1}_n (C_n
+D_n\Th_n) \Big\}dt\\[1mm]
\ns&\ds\q
+\[ Z_n X^{-1}_n- Y_n X^{-1}_n(C_n+D_n\Th_n)\]dW(t)\\[1mm]
\ns&\ds =\Big\{- (A_n+ A_{1,n})^*P_n- C^*_n
\Pi_n- Q_n\\[1mm]
\ns&\ds\q
+P_n \big[(C_n+D_n\Th_n)^2-A_n-A_{1,n}-B_n\Th_n\big] \\[1mm]
\ns&\ds \q -\Pi_n (C_n +D_n\Th_n)
\Big\}dt+\big[\Pi_n-P_n (C_n+D_n\Th_n)
\big]dW(t).
\end{array}
$$
Hence, by \eqref{8.25-eq3},
$(P_n(\cd),\Lambda_n(\cd))$ solves the following
$\dbR^{n\times n}$-valued backward stochastic
differential equation:
\begin{equation}\label{6.7-eq8}
\left\{
\begin{array}{ll}\ds
dP_n =-\big[P_n(A_n+ A_{1,n}) + (A_n+ A_{1,n})^*
P_n + \Lambda_n C_n + C^*_n \Lambda_n  \\[1mm]
\ns\ds \qq\qq + C^*_n P_n C_n+(P_n B_n + C^*_n
P_nD_n + \Lambda_n D_n) \Th_n + Q_n \big]dt\\[1mm]
\ns\ds\qq\qq + \Lambda_n dW(t)
\qq\qq\qq\qq\qq\qq\qq\mbox{ in }\;[0,T],\cr
\ns\ds P_n(T)=G_n.
\end{array}
\right.
\end{equation}

For $\xi_{k}\in L^4_{\cF_t}(\Om;V')$, $u_{k}\in
L^4_\dbF(\Om;L^2(t,T;V'))$ and $v_{k}\in
L^4_\dbF(\Om;L^2(t,T;V'))$ ($k=1,2$), denote by
$x_1(\cd)$ and $x_2(\cd)$ respectively the mild
solutions to the equations \eqref{op-fsystem1}
and \eqref{op-fsystem2}. For $k=1,2$, let us
introduce the following two (forward) stochastic
differential equations:
\begin{equation}\label{6.8-eq1}
\left\{
\begin{array}{ll}
\ds dx_{k,n} = \big[(A_n+A_{1,n}) x_{k,n}+ u_{k,n}\big]dr \\[1mm]
\ns\ds\qq\qq +\big( C_n(t)x_{k,n} + v_{k,n} \big)dW(r) &\mbox{ in } [t,T],\\[1mm]
\ns\ds x_{k,n}(t)= \xi_{k,n},
\end{array}
\right.
\end{equation}
where
$u_{k,n}=\G_n u_k$, $ v_{k,n}=\G_n v_k$ and $
\xi_{k,n}=\G_n\xi_k$.

Clearly,
\begin{equation}\label{6.8-eq9}
\left\{
\begin{array}{ll}\ds
\lim_{n\to+\infty} \xi_{k,n} = \xi_k &\mbox{ in
} L^4_{\cF_t}(\Om;V'),
\\[1mm]
\ns\ds \lim_{n\to+\infty} u_{k,n} = u_k &\mbox{
in }
L^4_\dbF(\Om;L^2(t,T;V')),\\[1mm]
\ns\ds \lim_{n\to+\infty} v_{k,n} = v_k &\mbox{
in } L^4_\dbF(\Om;L^2(t,T;V')).
\end{array}
\right.
\end{equation}
From Lemma \ref{lm16} and \eqref{6.8-eq9}, for
$k=1,2$, we get that\vspace{1mm}
\begin{equation}\label{6.8-eq10}
\lim_{n\to+\infty} x_{k,n}(\cd) = x_k(\cd)
\mbox{ in } L^4_\dbF(\Om;C([t,T];V')).
\end{equation}
By It\^o's formula, and using
\eqref{6.7-eq8}--\eqref{6.8-eq1}, we arrive at
\begin{eqnarray*}
&&  d \big\langle P_n
x_{1,n},x_{2,n}\big\rangle_{H_n}\\[1mm]
&& = \big\langle d P_n
x_{1,n},x_{2,n}\big\rangle_{H_n} + \big\langle
P_n dx_{1,n},x_{2,n}\big\rangle_{H_n} +
\big\langle P_n
x_{1,n},dx_{2,n}\big\rangle_{H_n}\\[1mm]
&& \q + \big\langle d P_n
dx_{1,n},x_{2,n}\big\rangle_{H_n} + \big\langle
dP_n x_{1,n},dx_{2,n}\big\rangle_{H_n} +
\big\langle
P_n dx_{1,n},dx_{2,n}\big\rangle_{H_n} \\[1mm]
&&  =  \big\langle - \big[P_n(A_n + A_{1,n}) +
(A_n + A_{1,n})^* P_n +  \Lambda_n C_n + C^*_n
\Lambda_n
 + C^*_n P_n C_n   \\[1mm]
&&  \q + (P_n B_n  +  C^*_n P_nD_n  +  \Lambda_n
D_n) \Th_n   +  Q_n \big]x_{1,n}, x_{2,n}
\big\rangle_{H_n} dr  \\[1mm]
&&  \q + \big\langle \Lambda_n x_{1,n},
x_{2,n}\big\rangle_{H_n}dW(r)+ \big\langle P_n
[(A_n + A_{1,n}) x_{1,n} + u_{1,n}]
,x_{2,n}\big\rangle_{H_n}dr\\[1mm]
&&  \q + \big\langle P_n (C_n x_{1,n}\!  +\!
v_{1,n}) ,x_{2,n}\big\rangle_{H_n}dW(r)\!+\!
\big\langle P_nx_{1,n},(A_n \!+\! A_{1,n})
x_{2,n} \!+\!
u_{2,n}\big\rangle_{H_n}dr \\[1mm]
&&  \q+ \big\langle P_nx_{1,n} ,C_n x_{2,n}  +
v_{2,n}\big\rangle_{H_n}dW(r) + \big\langle
\Lambda_n (C_n x_{1,n}  +
v_{1,n}),x_{2,n}\big\rangle_{H_n}dr \\[1mm]
&& \q + \big\langle \Lambda_n x_{1,n},C_n
x_{2,n}  + v_{2,n}\big\rangle_{H_n}dr+
\big\langle P_n (C_n x_{1,n}  + v_{1,n}),C_n
x_{2,n}  +
v_{2,n}\big\rangle_{H_n}dr\\[1mm]
&&  =  \big\langle - \big[ (P_n B_n + C^*_n
P_nD_n + \Lambda_n D_n) \Th_n + Q_n
\big]x_{1,n},x_{2,n} \big\rangle_{H_n}dr
\\[1mm]
&&  \q +\big\langle P_n u_{1,n},
x_{2,n}\big\rangle_{H_n}dr+ \big\langle P_n
x_{1,n}, u_{2,n}\big\rangle_{H_n}dr+ \big\langle
P_n C_n x_{1,n}, v_{2,n}\big\rangle_{H_n}dr \\[1mm]
&& \q + \big\langle P_nv_{1,n}, C_nx_{2,n}+
v_{2,n}\big\rangle_{H_n}dr+ \big\langle
\Lambda_nv_{1,n}, x_{2,n}\big\rangle_{H_n}dr +
\big\langle \Lambda_nx_{1,n}, v_{2,n}
\big\rangle_{H_n}dr \\[1mm]
&&  \q+ \big\langle \Lambda_n x_{1,n},
x_{2,n}\big\rangle_{H_n}dW(r)+ \langle P_n (C_n
x_{1,n}  + v_{1,n}) ,x_{2,n}\rangle_{H_n}dW(r)
\\[1mm]
&& \q + \langle P_nx_{1,n} ,C_n x_{2,n}  +
v_{2,n}\rangle_{H_n}dW(r).
\end{eqnarray*}
This implies that, for any $t\in [0,T]$,
\begin{eqnarray}\label{6.8-eq3}
&&\3n \mE\langle G_n
x_{1,n}(T),x_{2,n}(T)\rangle_{H_n}  + \mE
\int_t^T \big\langle Q_n(\tau) x_{1,n}(\tau),
x_{2,n}(\tau)
\big\rangle_{H_n}d\tau\nonumber\\[1mm]
&&\3n  + \mE \int_t^T \big\langle
\big[P_n(\tau)B_n(\tau) + C_n(\tau)^*
P_n(\tau)D_n(\tau)+ \Lambda_n(\tau)
D_n(\tau)\big]\nonumber\\[1mm]
&&\3n \qq\qq  \times\Th_n(\tau) x_{1,n}(\tau),
x_{2,n}(\tau) \big\rangle_{H_n}d\tau
 \\[1mm]
&&\3n  =\mE\big\langle P_n(t)
\xi_{1,n},\xi_{2,n} \big\rangle_{H_n} + \mE
\int_t^T \big\langle P_n(\tau)u_{1,n}(\tau),
x_{2,n}(\tau)\big\rangle_{H_n}d\tau \nonumber \\[1mm]
&&\3n \q+ \mE\!\! \int_t^T\!\! \big\langle
P_n(\tau)x_{1,n}(\tau),
u_{2,n}(\tau)\big\rangle_{H_n}\!d\tau \! +\!
\mE\! \!\int_t^T\! \!\big\langle
P_n(\tau)C_n(\tau)x_{1,n}(\tau),
v_{2,n}(\tau)\big\rangle_{H_n}\!d\tau \nonumber\\[1mm]
&& \3n \q + \mE
\int_t^T \big\langle  P_n(\tau)v_{1,n}(\tau), C_n(\tau)x_{2,n}(\tau)+ v_{2,n}(\tau)\big\rangle_{H_n}d\tau\nonumber\\
&&\3n \q + \mE \int_t^T \big\langle
\Lambda_n(\tau)v_{1,n}(\tau),
x_{2,n}(\tau)\big\rangle_{H_n}d\tau+ \mE
\int_t^T \big\langle
\Lambda_n(\tau)x_{1,n}(\tau), v_{2,n}(\tau)
\big\rangle_{H_n}d\tau.\nonumber
\end{eqnarray}

\medskip

\ms

{\bf Step 4}. In this step, we derive some
properties of $P$ and $\Lambda$.

\vspace{0.1cm}

Let $t\in [0,T)$ and $\eta\in
L^2_{\cF_t}(\Om;H)$. Consider the following
forward-backward stochastic evolution equation:
\begin{equation}\label{5.14-eq12}
\left\{
\begin{array}{ll}\ds
d x^t(r)=\big(A+A_1+B\Th\big) x^tdr + \big(C +
D \Th \big)x^tdW(r) &\mbox{\rm in } (t,T],\\[1mm]
dy^t(r) = -\big[(A+A_1)^* y^t + C^* z^t + Qx^t
\big] dr +  z^tdW(r) &\mbox{\rm in } [t,T),
\\[1mm] \ns\ds x^t(t)=\eta, \q
y^t(T)= G x^t(T).
\end{array}
\right.
\end{equation}
By Corollary \ref{5.7-prop1}, it is easy to see
that \eqref{5.14-eq12} admits a unique mild
solution\vspace{0.2mm}
$$
\begin{array}{ll}\ds
\big(x^t(\cd),y^t(\cd),z^t(\cd)\big)\big(\equiv\big(x^t(\cd;\eta),y^t(\cd;\eta),z^t(\cd;\eta)\big)\big)\\[1mm]
\ns\ds\in L^2_\dbF(\Om;C([t,T];H))\times
L^2_\dbF(\Om;C([t,T];H))\times L^2_\dbF(t,T;H)
\end{array}
$$
such that
\begin{equation}\label{5.72016-eq3}
R\Th x^t(r)+B^* y^t(r)+D^* z^t(r) =0,\q\mbox{for
a.e. }  (r,\om)\in (t,T)\times\Om.
\end{equation}
For every $r\in [t,T]$, define two families of
operators $X^t_r$ and $Y^t_r$ on
$L^2_{\cF_t}(\Om;H)$ as follows:
$$
X^t_r \eta \triangleq  x^t(r;\eta), \q Y^t_r
\eta  \triangleq y^t(r;\eta).
$$
For a.e. $r\in [t,T]$, define a family of
operators $Z^t(r)$ on $L^2_{\cF_t}(\Om;H)$ by
$$
 Z^t(r)\eta \triangleq z^t(r;\eta).
$$
It follows from Lemmas \ref{lm2}--\ref{lm3}
that\vspace{1mm}
\begin{equation}\label{5.26-eq8.1}
\begin{array}{ll}\ds
|X^t_r\eta|_{L^2_{\cF_r}(\Om;H)}\leq
\cC|\eta|_{L^2_{\cF_t}(\Om;H)},\\[1mm]
\ns\ds |Y^t_r\eta|_{L^2_{\cF_r}(\Om;H)}\leq
\cC|\eta|_{L^2_{\cF_t}(\Om;H)},\\[1mm]
\ns\ds |Z^t(\cd)\eta|_{L^2_{\dbF}(t,T;H)}\leq
\cC|\eta|_{L^2_{\cF_t}(\Om;H)}.
\end{array} \qq\forall\, r\in [t,T],\;\eta\in
{L^2_{\cF_t}(\Om;H)}.
\end{equation}
This indicates that $$X^t_r,  Y^t_r\in
\cL(L^2_{\cF_t}(\Om;H);L^2_{\cF_r}(\Om;H))$$ for
every $r\in [t,T]$ and $$Z^t(\cd)\in
\cL(L^2_{\cF_t}(\Om;H);L^2_{\dbF}(t,T;H)).$$

\vspace{0.1cm}

By \eqref{5.7-eq4.1} and \eqref{5.14-eq12}, it
is easy to see that, for any $\zeta\in H$,
$$
X^t_rX(t)\zeta = x^t(r;X(t)\zeta)=\hat
x(r;\zeta).
$$
Thus,
$$
Y^t_t X(t)\zeta = y^t(t;X(t)\zeta)=Y(t)\zeta
$$
and
$$
Z^t(\tau) X(t)\zeta =
z^t(\tau;X(t)\zeta)=Z(\tau)\zeta,\q \mbox{ for
a.e. }\tau\in [t,T].
$$
This implies that
\begin{equation}\label{5.14-eq13}
Y^t_t =Y(t)\wt X(t)^*\q \mbox{ for all } t\in
[0,T],\;\ \dbP\mbox{-a.s.,}
\end{equation}
and
\begin{equation}\label{5.14-eq13.1}
Z^t(\tau) =Z(\tau)\wt X(t)^*\q \mbox{ for a.e. }
t\in [0,T], \; \tau\in [t,T],\;\dbP\mbox{-a.s.}
\end{equation}

Let $\eta, \xi\in L^2_{\cF_t}(\Om;H)$. Since
$Y^t_r \eta=y^t(r;\eta)$ and $X^t_r
\xi=x^t(r;\xi)$, applying It\^o's formula to
$\langle y^t(\cd;\eta),x^t(\cd;\xi)\rangle_H$
and noting
\eqref{5.14-eq12}--\eqref{5.72016-eq3}, we
obtain that\vspace{0.1cm}
$$
\begin{array}{ll}\ds
\mE\langle GX^t_T \eta,X^t_T\xi \rangle_H  -
\mE\langle Y^t_t\eta, \xi \rangle_H\\[1mm]
\ns\ds =  -\mE \int_t^T \big(\langle
Q(r)X^t_r\eta,X^t_r\xi \rangle_H + \langle
R(r)\Th(r)X^t_r\eta,\Th(r)X^t_r\xi \rangle_H
\big)dr.
\end{array}
$$
Therefore,
$$
\begin{array}{ll}\ds
\mE\langle Y^t_t\eta, \xi \rangle_H\\[1mm]
\ns\ds = \mE\Big\langle (X^{t}_T)^* GX^t_T\eta +
\int_t^T \big((X^{t}_r)^* Q(r)X^t_r\eta  +
(X^{t}_r)^* \Th(r)^* R(r)\Th(r)X^t_r\eta
\big)dr,\xi \Big\rangle_H.
\end{array}
$$
From this, we conclude that, for any $\eta\in
L^2_{\cF_t}(\Om;H)$,
\begin{equation}\label{5.14-eq11}
\begin{array}{ll}\ds
Y^t_t\eta\3n&\ds  = \mE\((X^{t}_T)^* GX^t_T\eta
\\[1mm]
\ns&\ds \qq + \int_t^T \big((X^{t}_r
)^*Q(r)X^t_r\eta + (X^{t}_r)^*\Th(r)^*
R(r)\Th(r)X^t_r\eta \big)dr\;\Big|\;\cF_t\),
\end{array}
\end{equation}
which deduces that $$Y(t)\wt X(t)^*=Y_t^t
\;\mbox{ is symmetric for any }t\in [0,T],\q
\dbP\mbox{-a.s.}$$ Further, \eqref{5.14-eq11}
together with \eqref{5.26-eq8.1} implies that
for any $t\in [0,T]$ and $\eta \in
L^2_{\cF_t}(\Om;$ $H)$,
\begin{equation}\label{5.14-eq15}
\mE|Y^t_t\eta|_H^2 \leq \cC\mE|\eta|_H^2,
\end{equation}
where $\cC$ is  independent of $t\in [0,T]$.
According to \eqref{5.14-eq15}, we find that
\begin{equation}\label{5.14-eq16.1}
|Y(t)\wt
X(t)^*|_{\cL(L^2_{\cF_t}(\Om;H);\;L^2_{\cF_t}(\Om;H))}\leq
\cC.
\end{equation}
Thus, from \eqref{6.8-eq13.1}, \eqref{5.14-eq13}
and \eqref{5.14-eq16.1}, it follows that, for
some positive constant $\cC_0$,
\begin{equation}\label{5.14-eq9}
|P(t)|_{\cL(L^2_{\cF_t}(\Om;H);\;L^2_{\cF_t}(\Om;H))}\le
\cC_0, \qq\forall\;t\in [0,T].
\end{equation}

We claim that,
\begin{equation}\label{12.30-eq3}
|P(t)|_{\cL(H)}\leq \cC_0, \qq\forall\;t\in
[0,T],\q\dbP\mbox{-a.s.}
\end{equation}
Otherwise, there would exist $\e_0>0$ and  $\wt
\Om\in\cF_t$ with $\dbP(\wt\Om)>0$ such that
$$
|P(t,\om)|_{\cL(H)}>\cC_0+\e_0,\qq \mbox{ for
a.e. }\om\in\wt\Om.
$$
Let $\{\eta_k\}_{k=1}^\infty$ be a dense subset
of the unit sphere of $H$. Then, for a.e.
$\om\in\wt\Om$, there is an $\eta_\om\in
\{\eta_k\}_{k=1}^\infty$ such that
$$
|P(t,\om)\eta_\om|_H\geq
|P(t,\om)|_{\cL(H)}-\frac{\e_0}{2}>\cC_0+\frac{\e_0}{2}.
$$
Write
$$
\left\{
\begin{array}{ll}\ds
\Om_1 = \Big\{\om\in\wt\Om\;\Big|\;
|P(t,\om)\eta_1|_H>|P(t,\om)|_{\cL(H)}-\frac{\e_0}{2}\Big\},\\[3mm]
\ns\ds \Om_n =\Big\{\om\in\wt\Om\;\Big|\;
|P(t,\om)\eta_n|_H>|P(t,\om)|_{\cL(H)}-\frac{\e_0}{2}\Big\}\setminus
\(\bigcup_{k=1}^{n-1}\Om_k\), \\
\ns\ds\hspace{7.6cm}\mbox{ for } n=2,3,\cdots.
\end{array}
\right.
$$
Since $$P(t,\om)\eta_n \in L^2_{\cF_t}(\Om;H)
\mbox{ for all } n\in \dbN,$$ we see that
$\{\Om_n\}_{n=1}^\infty\subset \cF_t$ and
$\ds\dbP(\wt\Om)=\sum_{n=1}^\infty\dbP(\Om_n)$.
Hence,
$$
\Big|P(t,\om)\sum_{n=1}^\infty
\chi_{\Om_n}\eta_n\Big|_H \geq
\cC_0+\frac{\e_0}{2},\qq \mbox{ for a.e.
}\om\in\wt\Om.
$$
Therefore,
\begin{equation*}\label{12.30-eq1}
\mE\Big|P(t,\om)\sum_{n=1}^\infty
\chi_{\Om_n}\eta_n\Big|_H^2 \geq
\(\cC_0+\frac{\e_0}{2}\)^2\dbP(\wt\Om).
\end{equation*}
On the other hand, it follows from
\eqref{5.14-eq9} that
\begin{equation*}\label{12.30-eq2}
\begin{array}{ll}\ds
\mE\Big|P(t,\om)\sum_{n=1}^\infty
\chi_{\Om_n}\eta_n\Big|_H^2 \3n&\ds\leq
\cC_0^2\mE\Big|\sum_{n=1}^\infty
\chi_{\Om_n}\eta_n\Big|_H^2
=\cC_0^2\sum_{n=1}^\infty\dbP(\Om_n) =
\cC_0^2\dbP(\wt\Om).
\end{array}
\end{equation*}
These lead to a contradiction. Hence,
\eqref{12.30-eq3} holds. Since the constant
$\cC_0$ (in \eqref{5.14-eq9}) is independent of
$t\in [0,T]$, it follows that \vspace{1mm}
\begin{equation}\label{12.30-eq4}
|P(t,\om)|_{\cL(H)}\leq \cC_0, \qq \mbox{ for
a.e. }(t,\om)\in [0,T]\times\Om.
\end{equation}

Similar to the proof of \eqref{5.14-eq11}, we
can show that for any $\eta\in
L^2_{\cF_t}(\Om;V')$,
$$
\begin{array}{ll}\ds
Y^t_t\eta \3n&\ds = \mE\((X^{t}_T)^* GX^t_T\eta
+ \int_t^T \big((X^{t}_r )^*Q(r)X^t_r\eta +
(X^{t}_r)^*\Th(r)^* R(r)\Th(r)X^t_r\eta
\big)dr\;\Big|\;\cF_t\).
\end{array}
$$
This, together with ({\bf AS3}) and ({\bf AS4}),
implies that $$
|P(t)|_{\cL(L^2_{\cF_t}(\Om;V');\;L^2_{\cF_t}(\Om;V'))}\le
\cC. $$ Then, similar to the proof of
\eqref{12.30-eq4}, we obtain that
\begin{equation}\label{12.30-eq4.1}
|P(t,\om)|_{\cL(V')}\leq \cC, \qq \mbox{ for
a.e. }(t,\om)\in [0,T]\times\Om.
\end{equation}

Next, we prove that (Recall that $\Lambda(\cd)$
is a pointwise defined operator)
 \bel{20170101e1}
 \Lambda(t,\om)=\Lambda(t,\om)^*,\qq
\mbox{for a.e. }(t,\om)\in (0,T)\times \Om.
 \ee
For this purpose, let $\zeta,\kappa\in V'$,
$\zeta_n=\G_n\zeta$ and $\kappa_n=\G_n\kappa$
for $n=1,2,\cds$ From \eqref{6.7-eq8}, we get
that£¬ for each $t\in [0,T]$,
\begin{equation}\label{7.3-eq1}
\begin{array}{ll}\ds
\q\big\langle
G_n\zeta_n,\kappa_n\big\rangle_{\dbR^n}-\big\langle
P_n(t)\zeta_n,\kappa_n\big\rangle_{\dbR^n}\\[1mm]
\ns\ds =-\int_t^T\big\langle\big[P_n(A_n
+A_{1,n})+ (A_n +A_{1,n})^* P_n + \Lambda_n C_n
+ C^*_n \Lambda_n  \\[1mm]
\ns\ds \qq\qq + C^*_n P_n C_n+(P_n B_n + C^*_n
P_nD_n + \Lambda_n D_n) \Th_n + Q_n
\big]\zeta_n,\kappa_n\big\rangle_{\dbR^n}d\tau\\[1mm]
\ns\ds\q  +
\int_t^T\big\langle\Lambda_n\zeta_n,\kappa_n
\big\rangle_{\dbR^n}dW(\tau).
\end{array}
\end{equation}
By \eqref{8.25-eq5}, we know that
\begin{equation}\label{10.10-eq24}
\begin{array}{ll}\ds
\lim_{n\to\infty}\mE\Big|\int_t^T\big\langle\Lambda
\zeta,\kappa
\big\rangle_{V,V'}dW(\tau)-\int_t^T\big\langle\Lambda_n\zeta_n,\kappa_n
\big\rangle_{\dbR^n}dW(\tau)\Big|^2\\[1mm]
\ns\ds =
\lim_{n\to\infty}\mE\Big|\int_t^T\big\langle\Lambda
\zeta,\kappa
\big\rangle_{V,V'}dW(\tau)-\int_t^T\big\langle\Lambda_n\zeta_n,\kappa_n
\big\rangle_{V,V'}dW(\tau)\Big|^2\\[1mm]
\ns\ds =
\lim_{n\to\infty}\mE\Big|\int_t^T\big\langle\Lambda
\zeta - \Lambda_n\zeta_n ,\kappa
\big\rangle_{V,V'}dW(\tau) + \int_t^T\big\langle
\Lambda_n\zeta_n,\kappa -
\kappa_n \big\rangle_{V,V'}dW(\tau)\Big|^2\\[1mm]
\ns\ds \leq 2 \lim_{n\to\infty}\mE
\int_0^T\big|\big\langle\Lambda \zeta -
\Lambda_n\zeta_n ,\kappa
\big\rangle_{V,V'}\big|^2d\tau \!+
2\lim_{n\to\infty}\mE\int_0^T\!\big|\big\langle
\Lambda_n\zeta_n,\kappa \!-\! \kappa_n
\big\rangle_{V,V'}\big|^2d\tau\\[1mm]
\ns\ds \leq 2 \lim_{n\to\infty}\mE
\int_0^T\big|\Lambda \zeta -
\Lambda_n\zeta_n\big|_V^2
\big|\kappa\big|_{V'}^2d\tau +
2\lim_{n\to\infty}\mE\int_t^T\big|
\Lambda_n\zeta_n\big|_{V}^2\big|\kappa -
\kappa_n
\big|_{V'}^2d\tau \\[1mm]
\ns\ds \leq 2 \lim_{n\to\infty}\mE
\int_0^T\big|\Lambda \zeta -
\Lambda_n\zeta_n\big|_V^2
\big|\kappa\big|_{V'}^2d\tau +
2\lim_{n\to\infty}\mE\int_0^T\big|
\Lambda\zeta\big|_{V}^2\big|\kappa - \kappa_n
\big|_{V'}^2d\tau = 0.
\end{array}
\end{equation}
\vspace{1mm} Thus, there is a subsequence
$\{n_k\}_{k=1}^\infty\subset \{n\}_{n=1}^\infty$
such that
\begin{equation}\label{10.10-eq25}
\lim_{k\to\infty}\int_t^T\big\langle\Lambda_{n_k}\zeta_{n_k},\kappa_{n_k}
\big\rangle_{\dbR^{n_k}}dW(\tau)=\int_t^T\big\langle\Lambda
\zeta,\kappa \big\rangle_{V,V'}dW(\tau), \q
\dbP\mbox{-a.s.}
\end{equation}
Combining \eqref{6.8-eq4}, \eqref{8.25-eq5},
\eqref{7.3-eq1} and \eqref{10.10-eq25}, we see
that\vspace{1mm}
\begin{equation}\label{8.25-eq1}
\begin{array}{ll}\ds
\big\langle G
\zeta,\kappa\big\rangle_{H}-\big\langle
P(t)\zeta,\kappa\big\rangle_{H}\\[1mm]
\ns\ds =-\int_t^T\big\langle\big[P(A +A_{1}) +
(A +A_{1})^* P + \Lambda C+ C \Lambda +
C^*  P C  \\[1mm]
\ns\ds \qq\qq +(P B + C^* PD + \Lambda D)\Th + Q
\big]\zeta,\kappa\big\rangle_{V,V'}d\tau\\[1mm]
\ns\ds\q   + \int_t^T\big\langle\Lambda
\zeta,\kappa \big\rangle_{V,V'}dW(\tau).
\end{array}
\end{equation}

Since the operator $\Lambda(\cd)$ is pointwise
defined, it holds that
$$
\big\langle\Lambda \zeta,\kappa
\big\rangle_{V,V'}=\big\langle\Lambda
\zeta,\kappa
\big\rangle_{H}=\big\langle\zeta,\Lambda
^*\kappa \big\rangle_{H},\qq\forall
\zeta,\kappa\in V'.
$$
Therefore, by changing the positions of $\kappa$
and $\zeta$ in \eqref{8.25-eq1}, we obtain that,
for any $\zeta,\kappa\in V'$,
\begin{equation}\label{2.10-eq3}
\begin{array}{ll}\ds
\big\langle
\zeta,G\kappa\big\rangle_{H}-\big\langle
\zeta,P(t)\kappa\big\rangle_{H}\\[1mm]
\ns\ds =-\int_t^T\big\langle\big[P(A +A_{1}) +
(A +A_{1})^* P + \Lambda C + C \Lambda +
C^*  P C  \\[1mm]
\ns\ds \qq\qq +(P B + C^* PD + \Lambda D)\Th + Q
\big]^*\kappa,\zeta\big\rangle_{V,V'}d\tau\\[1mm]
\ns\ds\q  + \int_t^T\big\langle
\zeta,\Lambda^*\kappa \big\rangle_{H}dW(\tau).
\end{array}
\end{equation}
From \eqref{8.25-eq1} and \eqref{2.10-eq3}, we
find that
\begin{eqnarray}\label{8.25-eq1.1}
&& 0 = - \int_t^T\!\big\langle\big\{\big[
\Lambda C + C \Lambda + (P B  + C^* PD + \Lambda
D)\Th\big]\\[1mm]
&&\qq - \big[ \Lambda C  + C \Lambda + (P B +
C^* PD + \Lambda D)\Th\big]^*
\big\}\zeta,\kappa\big\rangle_{H}d\tau \nonumber \\[1mm]
&&  \qq + \int_t^T\big\langle(\Lambda-\Lambda^*)
\zeta,\kappa \big\rangle_{H}dW(\tau).\nonumber
\end{eqnarray}
By \eqref{8.25-eq1.1} and the uniqueness of the
decomposition of semimartingales, we conclude
that for any $\zeta,\kappa\in V'$,\vspace{1mm}
$$
\big\langle\big(\Lambda(t,\om)-\Lambda(t,\om)^*\big)\zeta,
\kappa\big\rangle_{H} =0,\q \mbox{ for a.e. }
(t,\om)\in (0,T)\times\Om,
$$
which gives \eqref{20170101e1}.

\vspace{0.1cm}

From \eqref{5.7-eq5.1}, it follows that
\begin{equation}\label{5.7-eq6.1}
R\Th+B^*P+D^*\Pi =0, \q \mbox{ a.e. }
(t,\om)\in[0,T]\times\Om.
\end{equation}
This, together with \eqref{5.14-eq9},
\eqref{12.30-eq4.1} and ({\bf AS4}), implies
that
\begin{equation}\label{6.7-eq10}
D^*\Pi\in \Upsilon_2(H;U)\cap \Upsilon_2(V';\wt
U).
\end{equation}
According to \eqref{6.8-eq13.1},
\eqref{6.7-eq10}, ({\bf AS3}) and  ({\bf AS4}),
it holds that
\begin{equation}\label{6.8-eq14}
D^*\Lambda=D^*\Pi-D^*P(C+D\Th)\in
\Upsilon_2(H;U)\cap \Upsilon_2(V';\wt U).
\end{equation}
From \eqref{6.8-eq14}, we see that
\begin{equation}\label{6.8-eq14.1}
\Lambda D  \in \Upsilon_2(U;H)\cap
\Upsilon_2(\wt U;V').
\end{equation}

\ms

\ms

{\bf Step 5}. In this step,  we prove that a
variant of \eqref{6.18-eq1} holds. We shall do
this by taking $n\to\infty$ in \eqref{6.8-eq3}.

\vspace{0.1cm}

Denote by $\wt U'$ the dual space of $\wt U$
with respect to the pivot space $U$. From ({\bf
AS4}), \eqref{6.8-eq4}, \eqref{8.25-eq5} and
\eqref{6.8-eq10}, we obtain that
\begin{eqnarray}\label{10.10-eq3.1}
&&\3n\2n \lim_{n\to+\infty}\mE \int_t^T
\big\langle \big[P_n(\tau)B_n(\tau) +
C_n(\tau)^* P_n(\tau)D_n(\tau)+ \Lambda_n(\tau)
D_n(\tau)\big]\nonumber\\[1mm]
&& \qq\qq\q \times\Th_n(\tau) x_{1,n}(\tau),
x_{2,n}(\tau) \big\rangle_{H_n}d\tau \nonumber \\[1mm]
&& \3n\2n =  \lim_{n\to+\infty} \mE \int_t^T
\big\langle \Th_n(\tau) x_{1,n}(\tau),\nonumber\\[1mm]
&&\qq\qq \big[B_n(\tau)^*P_n(\tau)^*\! \! +\!
D_n(\tau)^*P_n(\tau)^*C_n(\tau) \!+\!
D_n(\tau)^*\Lambda_n(\tau)^* \big] x_{2,n}(\tau)
\big\rangle_{H_n}\! d\tau\\
&& \3n\2n = \mE \int_t^T \big\langle \Th(\tau)
x_{1}(\tau),\nonumber \\[1mm]
&&\qq\qq\big[B(\tau)^*P(\tau) + D(\tau)^*
P(\tau)C(\tau)+ D(\tau)^*\Lambda(\tau) \big]
x_{2}(\tau)\big\rangle_{\wt U, \wt
U'}d\tau\nonumber
\\
&& \3n\2n = \mE \int_t^T \big\langle
\big[P(\tau)B(\tau) + C^*(\tau) P(\tau)D(\tau)+
\Lambda(\tau) D(\tau)\big]\Th(\tau) x_{1}(\tau),
x_{2}(\tau)\big\rangle_{V, V'}d\tau.\nonumber
\end{eqnarray}
By \eqref{8.25-eq5} and \eqref{6.8-eq10}, we see
that, for $k=1,2$,
\begin{equation}\label{10.10-eq33}
\begin{array}{ll}\ds
\lim_{n\to\infty}\big|\Lambda_{n}(\cd)x_{k,n}(\cd)-\Lambda(\cd)x_{k}(\cd)\big|_{L^{\frac{4}{3}}_\dbF(\Om;L^2(0,T;V))}\\[1mm]
\ns\ds \leq \lim_{n\to\infty}
\big|\Lambda_{n}(\cd)\big(x_{k,n}(\cd)
-x_{k}(\cd)\big)\big|_{L^{\frac{4}{3}}_\dbF(\Om;L^2(0,T;V))}\\[1mm]
\ns\ds\q +
\lim_{n\to\infty}\big|\big(\Lambda_{n}(\cd)
-\Lambda(\cd)\big)x_{k}(\cd)\big|_{L^{\frac{4}{3}}_\dbF(\Om;L^2(0,T;V))} \\[1mm]
\ns\ds \leq
\lim_{n\to\infty}\big|\Lambda(\cd)\big|_{L^2_\dbF(0,T;\cL_2(V';V))}\big|x_{k,n}(\cd)-x_{i}(\cd)\big|_{L^{4}_\dbF(\Om;L^\infty(0,T;V'))}
\\
\ns\ds \q+
\lim_{n\to\infty}\big|\Lambda_{n}(\cd)-\Lambda(\cd)\big|_{L^2_\dbF(0,T;\cL_2(V';V))}
\big|x_{k}(\cd)\big|_{L^{4}_\dbF(\Om;L^\infty(0,T;V'))}=0.
\end{array}
\end{equation}

By \eqref{10.10-eq3.1}--\eqref{10.10-eq33},
using a similar argument for other terms in
\eqref{6.8-eq3}, and noting \eqref{20170101e1},
we can take $n\to\infty$ on both sides of this
equality to get that
\begin{eqnarray}\label{10.10-eq5}
&&\3n \mE\langle G x_{1}(T),x_{2}(T)\rangle_{H}
+ \mE \int_t^T \big\langle Q(\tau) x_{1}(\tau),
x_{2}(\tau)
\big\rangle_{H}d\tau\nonumber\\
&&\3n \q + \mE\! \int_t^T\! \big\langle
\big[P(\tau)B(\tau)\! +\! C(\tau)^*
P(\tau)D(\tau)\!+\! \Lambda(\tau)
D(\tau)\big]\Th(\tau) x_{1}(\tau), x_{2}(\tau)
\big\rangle_{V, V'}d\tau\nonumber
  \\
&&\3n  =\mE\big\langle P(t)\xi_{1},\xi_{2}
\big\rangle_{V, V'}+ \mE \int_t^T \big\langle
u_{1}(\tau),
P(\tau)^*x_{2}(\tau)\big\rangle_{V', V}d\tau   \nonumber\\
&&\3n \q + \mE \int_t^T \big\langle
P(\tau)x_{1}(\tau), u_{2}(\tau)\big\rangle_{V,
V'}d\tau+ \mE \int_t^T \big\langle
P(\tau)C(\tau)x_{1}(\tau),
v_{2}(\tau)\big\rangle_{ V, V'}d\tau \\
&&\3n\q + \mE\int_t^T \big\langle
P(\tau)v_{1}(\tau), C(\tau)x_{2}(\tau)+
v_{2}(\tau)\big\rangle_{V, V'}d\tau
\nonumber\\
&&\3n \q + \mE \int_t^T \big\langle v_{1}(\tau),
\Lambda(\tau)x_{2}(\tau)\big\rangle_{V',
V}d\tau+ \mE \int_t^T \big\langle
\Lambda(\tau)x_{1}(\tau), v_{2}(\tau)
\big\rangle_{V, V'}d\tau.\nonumber
\end{eqnarray}
Noting that $P(t)\in
\cL(L^2_{\cF_t}(\Om;H);L^2_{\cF_t}(\Om;H))$, we
have
\begin{equation}\label{10.10-eq34}
\begin{array}{ll}\ds
\mE\big\langle P(t)\xi_{1},\xi_{2}
\big\rangle_{V, V'}+ \mE \int_t^T \big\langle
u_{1}(\tau),
P(\tau)^*x_{2}(\tau)\big\rangle_{V', V}d\tau \\
\ns\ds\q + \mE \int_t^T \big\langle
P(\tau)x_{1}(\tau), u_{2}(\tau)\big\rangle_{V,
V'}d\tau  + \mE \int_t^T \big\langle
P(\tau)C(\tau)x_{1}(\tau),
v_{2}(\tau)\big\rangle_{ V, V'}d\tau \\
\ns\ds\q  + \mE\int_t^T \big\langle
P(\tau)v_{1}(\tau), C(\tau)x_{2}(\tau)+
v_{2}(\tau)\big\rangle_{V, V'}d\tau \\
\ns\ds = \mE\big\langle P(t)\xi_{1},\xi_{2}
\big\rangle_{H}+ \mE \int_t^T \big\langle
u_{1}(\tau),
P(\tau)^*x_{2}(\tau)\big\rangle_{H}d\tau  \\
\ns\ds \q + \mE \int_t^T \big\langle
P(\tau)x_{1}(\tau), u_{2}(\tau)\big\rangle_{V,
V'}d\tau + \mE \int_t^T \big\langle
P(\tau)C(\tau)x_{1}(\tau),
v_{2}(\tau)\big\rangle_{H}d\tau \\
\ns\ds\q + \mE\int_t^T \big\langle
P(\tau)v_{1}(\tau), C(\tau)x_{2}(\tau)+
v_{2}(\tau)\big\rangle_{H}d\tau.
\end{array}
\end{equation}
From \eqref{6.8-eq14.1} and noting again
$P(t)\in
\cL(L^2_{\cF_t}(\Om;H);L^2_{\cF_t}(\Om;H))$, we
obtain that
\begin{equation}\label{10.10-eq35}
\begin{array}{ll}\ds
\mE \int_t^T \big\langle \big[P(\tau)B(\tau) +
C(\tau)^* P(\tau)D(\tau)+ \Lambda(\tau)
D(\tau)\big]\Th(\tau) x_{1}(\tau), x_{2}(\tau)
\big\rangle_{V, V'}d\tau\\[2mm]
\ns\ds = \mE \int_t^T \big\langle
\big[P(\tau)B(\tau) + C(\tau)^* P(\tau)D(\tau)+
\Lambda(\tau) D(\tau)\big]\Th(\tau) x_{1}(\tau),
x_{2}(\tau) \big\rangle_{H}d\tau.
\end{array}
\end{equation}
Combing \eqref{10.10-eq5}, \eqref{10.10-eq34}
and \eqref{10.10-eq35}, we conclude that
\begin{eqnarray}\label{vari6.18-eq1}
&& \mE\langle Gx_{1}(T),x_{2}(T)\rangle_{H} +\mE
\int_t^T \big\langle Q(\tau) x_{1}(\tau),
x_{2}(\tau)
\big\rangle_{H}d\tau\nonumber\\
&&\q  + \mE\! \int_t^T\! \big\langle
\big[P(\tau)B(\tau)\! +\! C(\tau)^*
P(\tau)D(\tau)+ \Lambda(\tau)
D(\tau)\big]\Th(\tau) x_{1}(\tau), x_{2}(\tau)
\big\rangle_{H}d\tau
\nonumber  \\
&&   = \mE\big\langle P(t) \xi_{1},\xi_{2}
\big\rangle_{H} + \mE \int_t^T \big\langle
P(\tau)u_{1}(\tau),
x_{2}(\tau)\big\rangle_{H}d\tau  \\
&&  \q  + \mE \int_t^T \big\langle
P(\tau)x_{1}(\tau),
u_{2}(\tau)\big\rangle_{H}d\tau+ \mE
\int_t^T\big\langle P(\tau)C(\tau)x_{1}(\tau),
v_{2}(\tau)\big\rangle_{H}d\tau\nonumber\\ &&
\q+ \mE
\int_t^T \big\langle  P(\tau)v_{1}(\tau), C(\tau)x_{2}(\tau)+v_{2}(\tau)\big\rangle_{H}d\tau\nonumber\\
&&  \q + \mE \int_t^T \big\langle v_{1}(\tau),
\Lambda(\tau)x_2(\tau)\big\rangle_{ V', V}d\tau+
\mE \int_t^T \big\langle \Lambda(\tau)x_1(\tau),
v_{2}(\tau) \big\rangle_{ V,V'}d\tau.\nonumber
\end{eqnarray}

\vspace{0.2cm}

{\bf Step 6}. In this step,  we prove that a
variant of \eqref{10.10-eq10} holds.

\vspace{2mm}

For any $\xi_k\in L^2_{\cF_t}(\Om;H)$, $u_k\in
L^2_\dbF(t,T;H)$ and $v_k\in L^2_\dbF(t,T;U)$
($k=1,2$), denote by $x_1(\cd)$ and $x_2(\cd)$
respectively the mild solutions to the equations
\eqref{op-fsystem1} and \eqref{op-fsystem2}. We
can find six sequences
$$\{\xi_{k}^j\}_{j=1}^\infty\in
L^4_{\cF_t}(\Om;V'),$$
$$\{u_{k}^j\}_{j=1}^\infty \subset
L^{4}_{\dbF}(\Om;L^2(t,T;V'))$$ and
$$\{v_{k}^j\}_{j=1}^\infty  \subset
L^{4}_{\dbF}(\Om;L^2(t,T;\wt U)),$$ such that
\begin{equation}\label{10.10-eq11}
\begin{cases}\ds
\lim_{j\to\infty}\xi_{k}^j = \xi_k \q  \mbox{ in
}\;L^2_{\cF_t}(\Om;H),\\[1mm]
\ns\ds\lim_{j\to\infty}u_{k}^j = u_k \q  \mbox{
in
}\;L^2_{\dbF}(t,T;H),\\[1mm]
\ns\ds \lim_{j\to\infty}v_{k}^j = v_k \q  \mbox{
in }\;L^2_{\dbF}(t,T;U).
\end{cases}
\end{equation}
Denote by $x_1^j(\cd)$ (\resp $x_2^j(\cd)$) the
mild solution to the equation
\eqref{op-fsystem1} (\resp \eqref{op-fsystem2})
with $\xi_{1}$, $u_{1}$ and $v_{1}$ (\resp
$\xi_{2}$, $u_{2}$ and $v_{2}$) replaced
respectively by $\xi_{1}^j$, $u_{1}^j$ and
$v_{1}^j$ (\resp $\xi_{2}^j$, $u_{2}^j$ and
$v_{2}^j$),   and by $x^j_{k,n}$  the solution
to \eqref{6.8-eq1} with  $\xi_{k,n}$, $u_{k,n}$
and $v_{k,n}$ replaced respectively by $\G_n
\xi^j_{k}$, $\G_n u^j_{k}$ and $D_n v^j_{k}$. It
follows from \eqref{6.8-eq3} that\vspace{1mm}
\begin{eqnarray}\label{9.7-eq7}
&&\3n \mE\langle G_n
x_{1,n}^j(T),x_{2,n}^j(T)\rangle_{H_n}  + \mE
\int_t^T \big\langle Q_n(\tau) x_{1,n}^j(\tau),
x_{2,n}^j(\tau)
\big\rangle_{H_n}d\tau\nonumber\\
&&\3n \q + \mE \int_t^T \big\langle
\big[P_n(\tau)B_n(\tau) + C_n(\tau)^*
P_n(\tau)D_n(\tau)+ \Lambda_n(\tau)
D_n(\tau)\big]\nonumber\\
&& \qq\qq\times\Th_n(\tau) x_{1,n}^j(\tau),
x_{2,n}^j(\tau) \big\rangle_{H_n}d\tau
\nonumber  \\
&&\3n  =\mE\big\langle P_n(t)
\xi_{1}^j,\xi_{2}^j \big\rangle_{H_n}+ \mE
\int_t^T \big\langle P_n(\tau)u_{1}^j(\tau),
x_{2,n}^j(\tau)\big\rangle_{H_n}d\tau \nonumber\\
&&\3n \q  + \mE \int_t^T \big\langle
P_n(\tau)x_{1,n}^j(\tau),
u_{2}^j(\tau)\big\rangle_{H_n}d\tau \\
&&\3n \q + \mE \int_t^T \big\langle
P_n(\tau)C_n(\tau)x_{1,n}^j(\tau),
D_n(\tau)v_{2}^j(\tau)\big\rangle_{H_n}d\tau \nonumber\\
&& \3n \q + \mE\int_t^T \big\langle
P_n(\tau)D_n(\tau)v_{1}^j(\tau),
C_n(\tau)x_{2,n}^j(\tau)+
D_n(\tau)v_{2}^j(\tau)\big\rangle_{H_n}d\tau\nonumber
\\
&&\3n  \q + \mE \int_t^T \big\langle
\Lambda_n(\tau)D_n(\tau)v_{1}^j(\tau),
x_{2,n}^j(\tau)\big\rangle_{H_n}d\tau \nonumber\\
&&\3n \q + \mE \int_t^T \big\langle
\Lambda_n(\tau)x_{1,n}^j(\tau),
D_n(\tau)v_{2}^j(\tau)
\big\rangle_{H_n}d\tau.\nonumber
\end{eqnarray}

From Assumption ({\bf AS4}) and
\eqref{10.10-eq33}, for $k=1,2$ and $j\in\dbN$,
we have that
\begin{equation}\label{8.25-eq7}
\begin{array}{ll}\ds
\lim_{n\to+\infty} D(\cd,\cd)^*
\Lambda_{n}(\cd,\cd) x_{k,n}^j(\cd,\cd)
=D(\cd,\cd)^*\Lambda(\cd,\cd)
x_{k}^j(\cd,\cd)\q\mbox{ in }\;
L^{\frac{4}{3}}_{\dbF}(\Om;L^2(t,T;\wt U')).
\end{array}
\end{equation}
Therefore, we get that
\begin{equation}\label{9.7-eq12}
\begin{array}{ll}\ds
\ds \lim_{n\to\infty}\mE \int_t^T \big\langle
D_n(\tau)v_{1}^j(\tau),
\Lambda_n(\tau)x_{2,n}^j(\tau)\big\rangle_{H_n}d\tau
\\[1mm]
\ns\ds= \mE \int_t^T \big\langle v_{1}^j(\tau),
D(\tau)^*\Lambda(\tau)x_{2}^j(\tau)\big\rangle_{\wt
U, \wt U'}d\tau.
\end{array}
\end{equation}
By \eqref{6.8-eq14}, we see that that
$D^*\Lambda x_{2}^j\in L^2_\dbF(t,T;U)$. Hence,
\begin{equation}\label{8.25-eq9}
\begin{array}{ll}\ds
\ds \lim_{j\to\infty}\mE \int_t^T \big\langle
v_{1}^j(\tau),
D(\tau)^*\Lambda(\tau)x_{2}^j(\tau)\big\rangle_{\wt
U, \wt U'}d\tau \\[1mm] \ns\ds= \mE \int_t^T
\big\langle v_{1}(\tau),
D(\tau)^*\Lambda(\tau)x_{2}(\tau)\big\rangle_Ud\tau.
\end{array}
\end{equation}
Similarly,
\begin{equation}\label{8.25-eq10}
\begin{array}{ll}\ds
\ds \lim_{j\to\infty}\lim_{n\to\infty}\mE
\int_t^T \big\langle
\Lambda_n(\tau)x_{1,n}^j(\tau),
D_n(\tau)v_{2,n}^j(\tau)
\big\rangle_{H_n}d\tau\\
\ns\ds= \mE \int_t^T \big\langle
D(\tau)^*\Lambda(\tau)x_{1}(\tau),
v_{2}(\tau)\big\rangle_Ud\tau.
\end{array}
\end{equation}

Noting that $P(t)\in
\cL(L^2_{\cF_t}(\Om;H);L^2_{\cF_t}(\Om;H))$, we
have that
$$
\lim_{j\to\infty}P(t)\xi_{1}^j= P(t)\xi_{1}\q
\mbox{ in }\;\, L^2_{\cF_t}(\Om;H).
$$
Thus, we get that
\begin{equation*}\label{10.10-eq6}
\lim_{j\to\infty}\mE\big\langle
P(t)\xi_{1}^j,\xi_{2}^j \big\rangle_{V,
V'}=\lim_{j\to\infty}\mE\big\langle
P(t)\xi_{1}^j,\xi_{2}^j \big\rangle_{H} =
\mE\big\langle P(t)\xi_{1},\xi_{2}
\big\rangle_{H}.
\end{equation*}
Since $P(\cd)\in \Upsilon_2(H)$, we see that
$$
\lim_{j\to\infty}P(\cd)^*x_{2}^j=
P(\cd)^*x_{2}\q \mbox{ in } \;L^2_{\dbF}(t,T;H).
$$
Thus,
\begin{equation*}\label{10.10-eq7}
\begin{array}{ll}\ds
\lim_{j\to\infty}\mE \int_t^T \big\langle
u_{1}^j(\tau),
P(\tau)^*x_{2}^j(\tau)\big\rangle_{V', V}d\tau\\
\ns\ds =\lim_{j\to\infty}\mE \int_t^T
\big\langle u_{1}^j(\tau),
P(\tau)^*x_{2}^j(\tau)\big\rangle_{H}d\tau\\
\ns\ds = \mE \int_t^T \big\langle u_{1}(\tau),
P(\tau)^*x_{2}(\tau)\big\rangle_{H}d\tau\\
\ns\ds = \mE \int_t^T \big\langle P(\tau)
u_{1}(\tau), x_{2}(\tau)\big\rangle_{H}d\tau.
\end{array}
\end{equation*}

By \eqref{8.25-eq9}, \eqref{8.25-eq10} and a
similar proof of \eqref{6.18-eq1}, we can get
that
\begin{eqnarray}\label{8.25-eq10.1}
&& \3n\3n\mE\langle
Gx_{1}(T),x_{2}(T)\rangle_{H} +\mE \int_t^T
\big\langle Q(\tau) x_{1}(\tau), x_{2}(\tau)
\big\rangle_{H}d\tau\nonumber\\
&&\3n\3n \q + \mE \int_t^T \big\langle
\big[P(\tau)B(\tau) + C(\tau)^* P(\tau)D(\tau)+
\Lambda(\tau) D(\tau)\big]\Th(\tau) x_{1}(\tau),
x_{2}(\tau) \big\rangle_{H}d\tau
\nonumber  \\
&& \3n\3n = \mE\big\langle P(t) \xi_{1},\xi_{2}
\big\rangle_{H} + \mE \int_t^T \big\langle
P(\tau)u_{1}(\tau),
x_{2}(\tau)\big\rangle_{H}d\tau \nonumber \\
&& \3n\3n\q + \mE \int_t^T \big\langle
P(\tau)x_{1}(\tau),
u_{2}(\tau)\big\rangle_{H}d\tau + \mE
\int_t^T\!\big\langle P(\tau)C(\tau)x_{1}(\tau),
D(\tau)v_{2}(\tau)\big\rangle_Ud\tau  \\
&& \3n\3n\q  +  \mE
\int_t^T\! \big\langle  P(\tau)D(\tau)v_{1}(\tau), C(\tau)x_{2}(\tau)\!+\!D(\tau)v_{2}(\tau)\big\rangle_{H}d\tau\nonumber\\
&& \3n\3n \q + \mE \int_t^T \big\langle
v_{1}(\tau),
D(\tau)^*\Lambda(\tau)x_2(\tau)\big\rangle_Ud\tau+
\mE \int_t^T \big\langle
D(\tau)^*\Lambda(\tau)x_1(\tau), v_{2}(\tau)
\big\rangle_Ud\tau.\nonumber
\end{eqnarray}

\ms

\ms

{\bf Step 7}. In this step,  we prove that the
assertion 1) in Definition \ref{4.8-def2} holds.
We first show that
$$K\geq 0,\q \mbox{a.e. } (t,\om)\in[0,T]\times\Om.$$
Let us replace the $v_1$ in \eqref{op-fsystem1}
and $v_2$  in \eqref{op-fsystem2} by $Dv_1$ and
$Dv_2$, respectively.

\vspace{0.1cm}

From  \eqref{5.7-eq6.1}, we see that
\begin{equation}\label{5.7-eq9}
\begin{array}{ll}\ds
0\3n&\ds=B^* P + D^*\big[\Lambda +
P(C+D\Th)\big] + R\Th\\[1mm]
\ns&\ds = B^* P + D^* PC + D^*\Lambda +K\Th.
\end{array}
\end{equation}
Thus,
\begin{equation}\label{6.8-eq15}
PB + C^*PD + \Lambda^*D  = -\Th^*K^*.
\end{equation}
Thanks to \eqref{8.25-eq10.1} and
\eqref{6.8-eq15}, and noting that
$\Lambda(\cd)^*=\Lambda(\cd)$  and
$K(\cd)^*=K(\cd)$, we obtain that
\begin{eqnarray}\label{5.7-eq7}
&&\3n\3n \mE\langle G x_1(T),x_2(T)\rangle_H +
\mE \int_t^T \big\langle Q(\tau) x_1(\tau),
x_2(\tau) \big\rangle_{H}d\tau \nonumber  \\[1mm]
&&\3n\3n \q -\mE \int_t^T \big\langle
\Th(\tau)^* K(\tau) \Th(\tau) x_1(\tau),
x_2(\tau) \big\rangle_{H}d\tau
\\[1mm]
&&\3n\3n  =\mE\big\langle P(t) \xi_1,\xi_2
\big\rangle_{H} + \mE \int_t^T \big\langle
P(\tau)u_1(\tau), x_2(\tau)\big\rangle_{H}d\tau
\nonumber \\[1mm]
&&\3n\3n  \q + \mE \int_t^T \big\langle
P(\tau)x_1(\tau), u_2(\tau)\big\rangle_{H}d\tau
+ \mE \int_t^T \big\langle P(\tau)C(\tau)x_1
(\tau), D(\tau)v_2 (\tau)\big\rangle_{H}d\tau \nonumber  \\[1mm]
&&\3n\3n\q + \mE
\int_t^T \big\langle  P(\tau)D(\tau)v_1 (\tau), C (\tau)x_2 (\tau)\!+\! D(\tau)v_2(\tau)\big\rangle_{H}d\tau\nonumber\\
&& \3n\3n \q  + \mE \int_t^T \big\langle
v_1(\tau),
D(\tau)^*\Lambda(\tau)x_2(\tau)\big\rangle_Ud\tau+
\mE \int_t^T \big\langle
D(\tau)^*\Lambda(\tau)x_1(\tau), v_2(\tau)
\big\rangle_Ud\tau.\nonumber
\end{eqnarray}

For any $(s,\eta)\in[0,T)\times
L_{\cF_s}^2(\Omega;H)$,  similar to the proof of
\eqref{5.31-eq1}, thanks to \eqref{5.7-eq7}, and
noting \eqref{9.7-eq10} and \eqref{6.8-eq15}, we
can show that
\begin{equation}\label{7.3-eq11}
\ba{ll}\ds
\cJ(s,\eta;u(\cd))\3n&\ds=\frac{1}{2}\dbE\Big(\big\langle
P(s)\eta,\eta\big\rangle_H+\int_s^T\big\langle K(u-\Th x),u-\Th x\big\rangle_Udr\Big)\\[1mm]
\ns&\ds =\cJ\big(s,\eta;\Th(\cd)\bar x(\cd)\big)
+\frac{1}{2}\dbE\int_s^T\big\langle K(u-\Th
x),u-\Th x\big\rangle_U d\tau. \ea
\end{equation}
Hence,
\begin{equation}\label{8.22-eq1}
\begin{array}{ll}\ds
\frac{1}{2}\mE\langle P(s)\eta,\eta
\rangle_H=\cJ(s,\eta;\Th(\cd)\bar x(\cd))\leq
\cJ(s,\eta;u(\cd)),\qq\forall\; u(\cd)\in
L^2_\dbF(s,T;U),
\end{array}
\end{equation}
if and only if
$$K\geq 0,\q\mbox{a.e. } (t,\om)\in[0,T]\times\Om.$$

Recall that $\dbF$ stands for the progressive
$\si$-field (in $[0,T]\times\Omega$) with
respect to $\mathbf{F}$. Clearly, a process
$\f:[0,T]\times\Om\to U$ is
$\mathbf{F}$-progressively measurable if and
only if it is $\dbF$-measurable. Note that for
any $\mathbf{F}$-adapted process $\f(\cd)$,
there is an $\mathbf{F}$-progressively
measurable process $\tilde \f(\cd)$ which is
stochastically equivalent to $\f(\cd)$ (see
\cite[pp. 68]{Meyer1} for example). Thus,  a
process $\f:[0,T]\times\Om\to U$ is
$\mathbf{F}$-adapted if and only if it is
$\dbF$-measurable.

\vspace{0.1cm}

Put
$$
\Xi_1\triangleq\big\{(t,\om)\in
(0,T)\times\Om\;\big|\; K(t,\om)h=0 \mbox{ for
some nonzero }h\in U\big\}
$$
and
$$
\Xi_2\triangleq\big\{(t,\om)\in
(0,T)\times\Om\;\big|\;|K(t,\om)h|_U>0 \mbox{
for all }h\in U_1\big\},
$$
where $$U_1\triangleq\big\{h\in
U\;\big|\;|h|_U=1\big\}.$$

\vspace{0.1cm}

Clearly, $$\Xi_1\cap \Xi_2=\emptyset$$ and
$$\Xi_1\cup \Xi_2=(0,T)\times\Om.$$ By the
definition of $\Xi_2$, we see that
$$
\Xi_2= \bigcup_{m=1}^\infty\Big\{(t,\om)\in
(0,T)\times\Om\;\Big|\;|K(t,\om)h|_U>\frac{1}{m}
\mbox{ for all }h\in U_1 \Big\}.
$$
Let $U_0$ be a countable dense subset of $U_1$.
Then
\begin{equation}\label{7.3-eq5}
\begin{array}{ll}\ds
\Xi_2\3n&\ds=
\bigcup_{m=1}^\infty\Big\{(t,\om)\in
(0,T)\times\Om\;\Big|\;|K(t,\om)h|_U>\frac{1}{m}
\mbox{
for all }h\in U_0 \Big\}\\[1mm]
\ns&\ds = \bigcup_{m=1}^\infty\bigcap_{h\in
U_0}\Big\{(t,\om)\in
(0,T)\times\Om\;\Big|\;|K(t,\om)h|_U>\frac{1}{m}\Big\}.
\end{array}
\end{equation}
Since $K(\cd,\cd)h\in L^2_\dbF(0,T;U)$, we get
that, for any $h\in U$,
$$
\Big\{(t,\om)\in
(0,T)\times\Om\;\Big|\;|K(t,\om)h|_U>\frac{1}{m}\Big\}\in\dbF.
$$
This, together with \eqref{7.3-eq5}, implies
that $\Xi_2\in \dbF$. So does $\Xi_1$.

\vspace{0.1cm}

We now show that
$$
K>0 \mbox{ for a.e. } (t,\omega)\in
[0,T]\times\Omega.
$$
Let us use the contradiction argument and assume
that this was untrue. Then the measure (given by
the product measure of the Lebesgue measure on
$[0,T]$ and the probability measure $\dbP$) of
$\Xi_1$ would be positive.

\vspace{0.1cm}

For a.e. $(t,\omega)\in\Xi_1$, put
$$
\Upsilon(t,\om)\triangleq\big\{h\in U_1\;\big|\;
K(t,\om)h=0\big\}.
$$
Clearly, $\Upsilon(t,\om)$ is closed in $U$.
Define a map $F:(0,T)\times\Om\to 2^{U}$ as
follows:
$$
F(t,\om)=\left\{
\begin{array}{ll}\ds
\Upsilon(t,\om) , &\mbox{ if }(t,\om)\in \Xi_1
\\[1mm]
\ns\ds 0, &\mbox{ if }(t,\om)\in  \Xi_2.
\end{array}
\right.
$$
Then, $F(t,\om)$ is closed for a.e. $(t,\om)\in
(0,T)\times\Om$.

\vspace{0.1cm}

We now prove that $F$ is $\dbF$-measurable. Let
$O$ be a closed subset of $U$ and $O_1=O\cap
U_1$. Put
\bel{201615e1} \Si_1 \triangleq\big\{(t,\om)\in
(0,T)\times\Om\;\big|\;F(t,\om)\cap
O\neq\emptyset\big\} \ee
and
\bel{201615e1.1} \Si_2\triangleq\big\{(t,\om)\in
(0,T)\times\Om\;\big|\;F(t,\om)\cap
O_1\neq\emptyset\big\}. \ee
Clearly, $\Si_1\supset\Si_2$. Moreover,
 $$
 \Si_1=\left\{
 \ba{ll}
 \Si_2\cup \Xi_2,\ &\mbox{if }0\in O,\\[2mm]
 \Si_2,\ &\mbox{if }0\notin O.
 \ea\right.
 $$
Write
$$
\Si_3 \triangleq\big\{(t,\om)\in
(0,T)\times\Om\;\big|\;|K(t,\om)h|_U>0 \mbox{
for all }h\in O_1 \big\}.\vspace{1mm}
$$
Clearly, $$\Si_2\cap\Si_3=\emptyset$$ and
$$(0,T)\times\Om = \Si_2\cup\Si_3.$$ Similar to
the above (for the proof of $\Xi_2\in \dbF$), we
can show that $\Si_3\in \dbF$. Hence, $\Si_2\in
\dbF$ and therefore so does $\Si_1$.

\vspace{0.1cm}

Now we apply Lemma \ref{lm6} to $F(\cd,\cd)$
with $(\wt\Om,\wt\cF)=((0,T)\times\Om,\dbF)$ to
find an $\mathbf{F}$-adapted process $f$ such
that
$$
Kf=0 \mbox{ for a.e. }(t,\om)\in (0,T)\times\Om.
$$
Noting that
$$
|f(t,\om)|_U\leq 1 \mbox{ for a.e. } (t,\om)\in
(0,T)\times\Om, \vspace{1mm}$$
we find that $f\in L^2_\dbF(0,T;U)$.
Furthermore, we have
$$
|f(t,\om)|_U= 1 \mbox{ for a.e. }(t,\om)\in
\Xi_1,
$$
which concludes that $|f|_{L^2_\dbF(0,T;U)}>0$.

\vspace{0.1cm}

By \eqref{7.3-eq11}, we see that $\Th \bar x +
f$ is also an optimal control. This contradicts
the uniqueness of the optimal control. Hence,
$K(t,\om)$ is invertible (but $K(t,\om)^{-1}$
does not need to be bounded) for a.e.
$(t,\om)\in [0,T]\times\Om$.

\vspace{0.1cm}

Further, we show that the domain of
$K(t,\om)^{-1}$ is dense in $U$ for a.e.
$(t,\om)\in [0,T]\times\Om$.

Denote by $\cR(K(t,\om))$ the range of
$K(t,\om)$. Clearly, $\cR(K(t,\om))\subset U$.

\vspace{0.1cm}

Put
$$
\wt\Xi_1\triangleq\big\{(t,\om)\in
(0,T)\times\Om\;\big|\;
\cR(K(t,\om))^{\perp}\neq \{0\}\big\}
$$
and
$$
\wt \Xi_2\triangleq\big\{(t,\om)\in
(0,T)\times\Om\;\big|\;\cR(K(t,\om))^{\perp}=
\{0\}\big\}.
$$
Clearly, $\wt\Xi_1\cup\wt\Xi_2=(0,T)\times\Om$.
By the definition of $\wt\Xi_2$,  we see that
$$
\begin{array}{ll}\ds
\wt\Xi_2= \bigcup_{m=1}^\infty\Big\{(t,\om)\in
(0,T)\times\Om\;\Big|\;\forall\,\tilde h\in
U_0,\mbox{ there is an }h\in U_0 \mbox{ such
that }\\
\ns\ds\hspace{5.4cm}\big|\langle K(t,\om)h,
\tilde h\rangle_U\big|> \frac{1}{m} \Big\}.
\end{array}
$$
Then
\begin{equation}\label{7.3-eq6}
\wt\Xi_2= \bigcap_{\tilde h\in U_0}\bigcup_{h\in
U_0} \bigcup_{m=1}^\infty\Big\{(t,\om)\in
(0,T)\times\Om\;\Big|\;\big|\langle K(t,\om)h,
\tilde h\rangle_U\big|> \frac{1}{m} \Big\}.
\end{equation}
Since $K(\cd,\cd)h\in L^2_\dbF(0,T;U)$, it
follows that, for any $h,\tilde h\in U$,
\begin{equation}\label{7.3-eq7}
\Big\{(t,\om)\in
(0,T)\times\Om\;\Big|\;\big|\langle K(t,\om)h,
\tilde h\rangle_U\big|>
\frac{1}{m}\Big\}\in\dbF.
\end{equation}
From \eqref{7.3-eq6} and \eqref{7.3-eq7}, we see
that $\wt\Xi_2\in \dbF$. Hence, $\wt\Xi_1\in
\dbF$.

\vspace{0.1cm}

It suffices to prove that $\cR(K(t,\om))$ is
dense in $U$ for a.e. $(t,\om)\in
(0,T)\times\Om$. To show this, we use the
contradiction argument. If $\cR(K(t,\om))$ was
not dense in $U$ for a.e. $(t,\om)\in
(0,T)\times\Om$., then the measure of $\wt\Xi_1$
would be positive.

\vspace{0.1cm}

For a.e. $(t,\omega)\in\wt\Xi_1$,
put\vspace{1mm}
$$
\wt\Upsilon(t,\om)\triangleq\big\{\tilde h\in
U_1\;\big|\; \langle K(t,\om)h, \tilde
h\rangle_U=0, \; \forall\;h\in U\big\}.
$$
Clearly, $\wt\Upsilon(t,\om)$ is closed in $U$.

\vspace{0.1cm}

Define a map $\wt F:(0,T)\times\Om\to 2^{U}$ as
follows:
$$
\wt F(t,\om)=\left\{
\begin{array}{ll}\ds
\wt\Upsilon(t,\om) , &\mbox{ if }(t,\om)\in
\wt\Xi_1
\\[1mm]
\ns\ds 0, &\mbox{ if }(t,\om)\in \wt \Xi_2.
\end{array}
\right.
$$
Then, $\wt F(t,\om)$ is closed for a.e.
$(t,\om)\in (0,T)\times\Om$.

\vspace{0.1cm}

We now prove that $\wt F$ is $\dbF$-measurable.
Similar to \eqref{201615e1} and
\eqref{201615e1.1}, put
$$
\wt\Si_1 \triangleq\big\{(t,\om)\in
(0,T)\times\Om\;\big|\;\wt F(t,\om)\cap
O\neq\emptyset\big\},
$$
and
$$
\wt\Si_2\triangleq\big\{(t,\om)\in
(0,T)\times\Om\;\big|\;\wt F(t,\om)\cap
O_1\neq\emptyset\big\}.\vspace{1mm}
$$
If $0\in O$, then
$\wt\Si_1=\wt\Si_2\cup\wt\Xi_2$. If $0\notin O$,
then $\wt\Si_1=\wt\Si_2$. Hence, we only need to
show that $\wt\Si_2\in \dbF$. Write
$$
\begin{array}{ll}\ds
\wt\Si_3 \triangleq \Big\{(t,\om) \in
(0,T)\times\Om\;\Big|\; \forall\,\tilde h \in
O_1,\mbox{ there is an }h\in  U_1 \mbox{ so that
}\\ \ns\ds \hspace{4.4cm}\langle K(t,\om)h,
\tilde h\rangle_U
>0 \Big\}.
\end{array}
$$
Then, $$(0,T)\times\Om = \wt\Si_2\cup\wt\Si_3$$
and $$\wt\Si_2\cap\wt\Si_3=\emptyset.$$ Hence,
it suffices to show that $\wt\Si_3\in \dbF$. Let
$O_0$ be a countable dense subset of $O_1$.
Clearly,
\begin{eqnarray}\label{7.3-eq9}
&&\wt\Si_3 =
\bigcup_{m=1}^\infty\Big\{(t,\om)\in
(0,T)\times\Om\;\Big|\;\forall\,\tilde h\in
O_1,\mbox{ there is an }h\in  U_1 \mbox{ such
that }\nonumber\\
&& \hspace{5.4cm}\langle K(t,\om)h, \tilde
h\rangle_U> \frac{1}{m} \Big\}\nonumber\\
&&\q \;\, = \bigcup_{m=1}^\infty\Big\{(t,\om)\in
(0,T)\times\Om\;\Big|\;\forall\,\tilde h\in
O_0,\mbox{ there is an }h\in U_0 \mbox{ such
that }\\
&& \hspace{5.4cm}\langle K(t,\om)h, \tilde
h\rangle_U> \frac{1}{m} \Big\}\nonumber\\
&&\q \;\, = \bigcup_{m=1}^\infty\bigcap_{\tilde
h\in O_0}\bigcup_{h\in U_0} \Big\{(t,\om)\in
(0,T)\times\Om\;\Big|\;\langle K(t,\om)h, \tilde
h\rangle_U> \frac{1}{m}\Big\}.\nonumber
\end{eqnarray}
For any $m\in \dbN$, $\tilde h\in O_0$ and $h\in
U_0$, noting that $K(\cd,\cd)h\in
L^2_\dbF(0,T;U)$, we deduce that
\begin{equation}\label{7.3-eq10}
\Big\{(t,\om)\in (0,T)\times\Om\;\Big|\;\langle
K(t,\om)h, \tilde h\rangle_U>
\frac{1}{m}\Big\}\in\dbF.
\end{equation}
From \eqref{7.3-eq9} and \eqref{7.3-eq10}, it
follows that $\wt\Si_3\in \dbF$. Hence,
$\wt\Si_2\in \dbF$.

\ms

Before continuing the proof, we recall the
following known measurable selection result
(e.g. \cite{Wagner}).
\begin{lemma}\label{lm6}
Let $F :\; (\Om,\cF) \to 2^H$ be a closed-valued
set
 mapping, $F (\om)\neq \emptyset$ for every
$\om\in\Om$, and for each open set $O\subset H$,
$$
F^{-1}(O)\triangleq\big\{
\om\in\Om\;\big|\;F(\om)\cap O\neq\emptyset
\big\}\in\cF.
$$
Then $F$ has a measurable selection $f:\Om\to
H$, i.e., there is an $H$-valued,
$\cF$-measurable function $f$ such that $f(\om)
\in F(\om)$ for every $\om\in\Om$.
\end{lemma}

\vspace{0.1cm}

Now let us return the proof of Theorem
\ref{5.7-th1}. We apply Lemma \ref{lm6} to $\wt
F(\cd,\cd)$ with
$(\wt\Om,\wt\cF)=((0,T)\times\Om,\dbF)$ to find
an $\mathbf{F}$-adapted process $\tilde f$ such
that
$$
\langle K(t,\om)h,\tilde f(t,\om)\rangle_U=0, \q
\forall\;h\in U, \mbox{ a.e. } (t,\om)\in
(0,T)\times\Om.
$$
Since
$$|\tilde f(t,\om)|_U\leq 1 \mbox{ for
a.e. }(t,\om)\in (0,T)\times\Om,$$
it holds that $\tilde f\in L^2_\dbF(0,T;U)$.
Furthermore, we have
$$
|\tilde f(t,\om)|_U= 1 \mbox{ for a.e.
}(t,\om)\in \wt\Xi_1,
$$
which implies that $|\tilde
f|_{L^2_\dbF(0,T;U)}>0$.

\vspace{0.1cm}

We claim that $\Th \bar x + \tilde f$ is also an
optimal control. Indeed, by the choice of
$\tilde f$, it holds that
$$
\dbE\int_s^T\big\langle K(u-\Th x-\tilde
f),\tilde f \big\rangle_U dr  =0
$$
and
$$
\begin{array}{ll}\ds
\dbE\int_s^T\big\langle K\tilde f,u-\Th x-\tilde
f\big\rangle_U dr = \dbE\int_s^T\big\langle
\tilde f,K(u-\Th x-\tilde f)\big\rangle_U dr=0.
\end{array}
$$
Therefore, for any $u(\cd)\in L^2_\dbF(s,T;U)$,
\begin{equation}\label{7.3-eq12}
\begin{array}{ll}\ds
\dbE\int_s^T\big\langle K(u-\Th x-\tilde
f),u-\Th x-\tilde f\big\rangle_U dr \\
\ns\ds = \dbE\int_s^T\big\langle K(u-\Th
x),u-\Th x\big\rangle_U dr.
\end{array}
\end{equation}
According to \eqref{7.3-eq11} and
\eqref{7.3-eq12}, we obtain that
$$\cJ(s,\eta;\Th(\cd)\bar x(\cd)+\tilde f)\leq \cJ(s,\eta;u),\q\forall u(\cd)\in L^2_\dbF(s,T;U),$$
which indicates that $\Th \bar x + \tilde f$ is
also an optimal control. This leads to a
contradiction to the uniqueness of the optimal
controls. Hence, $K(t,\om)^{-1}$ is densely
defined in $U$ for a.e. $(t,\om)\in
(0,T)\times\Om$.

\vspace{0.1cm}

Now, let us prove that $K(t,\om)^{-1}$ is a
closed operator for a.e. $(t,\om)\in
(0,T)\times\Om$. Let
$\{h_j\}_{j=1}^\infty\subset
\cD(K(t,\om)^{-1})$, $h\in U$ and $\hat h\in U$
satisfy
\begin{equation}\label{7.3-eq13}
\lim_{j\to\infty}h_j = h \mbox{ in }U
\end{equation}
and
\begin{equation}\label{7.3-eq14}
\lim_{j\to\infty}K(t,\om)^{-1}h_j = \hat h
\mbox{ in }U.
\end{equation}
From \eqref{7.3-eq14}, we obtain that $$\ds
\lim_{j\to\infty} h_j = K(t,\om)\hat h \mbox{ in
} U.$$
This, together with \eqref{7.3-eq13}, implies
that $ h = K(t,\om)\hat h$.
Hence, $h\in \cD(K(t,\om)^{-1})$ and
$K(t,\om)^{-1}h = \hat h$. This indicates that
the operator $K(t,\om)^{-1}$ is closed.
Therefore, the assertion 1) in Definition
\ref{4.8-def2} holds.

\vspace{0.1cm}

By \eqref{5.7-eq9},  we find that
\begin{equation}\label{5.31-eq2}
-K^{-1}(B^* P+ D^*\Lambda+D^* PC) =\Th.
\end{equation}
This implies \eqref{5.7-eq5}. Moreover, from
\eqref{vari6.18-eq1}, \eqref{5.7-eq7} and
\eqref{5.31-eq2}, we see that the assertions 2)
and 3) in Definition \ref{4.8-def2} holds.

\vspace{0.2cm}

{\bf Step 8}. In this step, we prove the
uniqueness of the transposition solution to
\eqref{5.5-eq6}.

Assume that $$(P_1(\cd),\Lambda(\cd)),
(P_2(\cd),\Lambda_{2}(\cd))\in C_{\dbF,w}([0,T];
L^{\infty}(\Om;\cL(H))) \times
L^2_{\dbF,w}(0,T;\cL(H))$$ are two transposition
solutions to \eqref{5.5-eq6}.

\ms

From \eqref{8.22-eq1}, we have that for any
$s\in[0,T)$ and $\eta\in L^2_{\cF_s}(\Om;H)$,
\begin{equation}\label{8.22-eq2}
\frac{1}{2}\mE\langle P_1(s)\eta,\eta
\rangle_H=\cJ(s,\eta;\Th(\cd)x(\cd))=\frac{1}{2}\mE\langle
P_2(s)\eta,\eta \rangle_H.
\end{equation}
Thus, for any $\xi,\eta\in L^2_{\cF_s}(\Om;H)$,
we have that\vspace{0.4mm}
$$
\mE\langle P_1(s)(\eta+\xi),\eta+\xi \rangle_H
=\mE\langle P_2(s)(\eta+\xi), \eta+\xi
\rangle_H,
$$
and
$$
\mE\langle P_1(s)(\eta -\xi),\eta -\xi \rangle_H
=\mE\langle P_2(s)(\eta -\xi),\eta -\xi
\rangle_H.\vspace{3mm}
$$
These, together with $P_1(\cd)=P_1(\cd)^*$ and
$P_2(\cd)=P_2(\cd)^*$, imply that\vspace{1mm}
\begin{equation}\label{8.22-eq3}
\mE\langle P_1(s)\eta,\xi \rangle_H  =\mE\langle
P_2(s)\eta,\xi \rangle_H,\qq
\forall\;\xi,\eta\in
L^2_{\cF_s}(\Om;H).\vspace{2mm}
\end{equation}
Hence, $$P_1(s)=P_2(s) \mbox{ for any }  s\in
[0,T], \; \dbP\mbox{-a.s.}$$

\vspace{0.1cm}

Let $v_2=0$ in \eqref{op-fsystem2}. By
\eqref{6.18-eq1} and noting\vspace{1mm}
$$
\big\langle K(\cd)^{-1} L(\cd) x_{1}(\cd),
L(\cd)x_{2}(\cd) \big\rangle_{H}=-\big\langle
\Th(\cd) x_{1}(\cd), K(\cd)\Th(\cd)x_{2}(\cd)
\big\rangle_{H},
$$
we see that for any $\xi_1\in
L^2_{\cF_0}(\Om;H)$, $u_1\in
L^4_\dbF(\Om;L^2(0,T;H))$, $v_1 \in
L^4_\dbF(\Om;$ $L^2(0,T;$ $V'))$,
\begin{equation}\label{8.22-eq9}
0=\mE \int_0^T \big\langle v_1(s),
\big(\Lambda_1(s) -
\Lambda_2(s)\big)x_2(s)\big\rangle_{V',V}ds.
\end{equation}
Consequently,
\begin{equation}\label{9.7-eq13}
\big(\Lambda_1 - \Lambda_2 \big)x_2 = 0 \;\mbox{
in }\; L^{\frac{4}{3}}_\dbF(\Om;L^2(0,T;V)).
\end{equation}
\vspace{0.1cm} This, together with Lemma
\ref{lm11}, implies that\vspace{1mm}
$$\Lambda_1=\Lambda_2 \mbox{ in
}L^2_\dbF(0,T;\cL_2(H,V)).$$  Hence, the desired
uniqueness follows. This completes the proof of
Theorem \ref{5.7-th1}.


\section{Existence of transposition solutions to some operator-valued BSREs}\label{sect6}


This chapter is devoted to proving the existence
of transposition solutions to \eqref{5.5-eq6}
under suitable assumptions on $A$, $A_1$, $B$,
$C$, $D$, $Q$, $G$ and $R$. For simplicity, we
assume that $U=H$.

\vspace{0.1cm}

Let us introduce the following assumptions:

\vspace{0.25cm}

({\bf AS5}) {\it $A$, $A_1$, $B$, $C$, $Q$, $G$
and $R$ are all infinite dimensional block
diagonal matrices and $D$ is invertible. }

\vspace{0.25cm}

For $F =A, A_1, B, C, Q, G, R$, denote by $F_k$
the $k$-th matrix of the main diagonal block (In
particular, $A_{1,k}$ is understood in this
way).

\vspace{0.25cm}

({\bf AS6}) {\it For each $k\in\dbN$, $F_k$ is
an $m_k\times m_k$ matrix and $m_k\leq M$ for
some $M\in\dbN$.}

\vspace{0.25cm}

({\bf AS7}) {\it For each $k\in\dbN$,
$$
R_k(t) = R_k(0) + \int_0^t R_{1,k}(s)ds +
\int_0^t R_{2,k}(s)dW(s),
$$
where $R_{1,k}(\cd), R_{2,k}(\cd) \in
L^\infty_\dbF(0,T;\cS(\dbR^{m_k}))$. }

\vspace{0.25cm}

({\bf AS8}) {\it For each $k\in\dbN$,
$$
R_kB_k+C_k^\top R_k+R_{2,k}=0,
$$
and
$$
\begin{array}{ll}\ds
\wt Q_k\3n&\ds\triangleq Q_k-R_{1,k} + C_k^\top
R_kC_k
+ R_k (B_kC_k-A_k-A_{1,k})\\
\ns&\ds\q + (B_kC_k-A_k-A_{1,k})^\top R_k \geq
0.
\end{array}
$$
}

\vspace{0.25cm}

({\bf AS9}) {\it For each $k\in\dbN$, the
Malliavinian derivative of $F_k$ is uniformly
bounded, i.e., there is a constant $\cC>0$ such
that for a.e. $\tau,t\in [0,T]$ and $k\in\dbN$,
$|\cD_\tau F_k(t)|\leq \cC$, $\dbP$-a.s., where
$\cD_\tau F_k(t)$ is the Malliavinian derivative
of $F_k(t)$ at $\tau$.}

\vspace{0.25cm}

We have the following result:
\begin{theorem}\label{10.25-th2}
Let Assumptions ({\bf AS5})--({\bf AS9}) hold.
Then the equation \eqref{5.5-eq6} admits a
unique transposition solution $$(P,\Lambda)\in
C_\dbF([0,T];L^\infty(\Om;\cL(H)))\times
L^\infty_\dbF(0,T;\cL(H)).$$
\end{theorem}
\begin{remark}
One may expect that \eqref{5.5-eq6} would admit
a transposition solution $(P,\Lambda)\in
C_\dbF([0,T];$ $L^\infty(\Om;\cL(H)))\times
L^\infty_\dbF(\Om;L^2(0,T;\cL(H)))$ without
further assumptions. Unfortunately, this is
incorrect even in finite dimensions, i.e.,
$H=\dbR^n$ (e.g. \cite[Example 6.2]{LWZ}).
\end{remark}

The uniqueness result in Theorem \ref{10.25-th2}
is  obvious. We only need to show the existence
result. Since $D$ is invertible, without loss of
generality, we may assume that $D=I$. Otherwise,
one simply takes $v=D^{-1}u$ as the control. Let
us present some preliminaries as follows.

\vspace{0.1cm}

Under Assumption ({\bf AS5}), the equation
\eqref{5.5-eq6} can be written as infinitely
many matrix equations as follows:
\begin{equation}\label{10.25-eq1}
\left\{
\begin{array}{ll}\ds
dP_k =-\big[ P_k(A_k+A_{1,k}) +
(A_k+A_{1,k})^\top P_k + \Lambda_k C_k +
C_k^\top
\Lambda_k  \\
\ns\ds \qq\qq  + C_k^\top P_kC_k+ Q_k- L_k^\top
K_k^{-1} L_k \big]dt+ \Lambda_k dW(t) &\mbox{ in
}[0,T],\cr
\ns\ds P_k(T)=G_k,
\end{array}
\right.
\end{equation}
where $ L_k= B_k^\top P_k+D_k^\top
(P_kC_k+\Lambda_k)$. By  Theorem 2.2 in
\cite{QZ}, the equation \eqref{10.25-eq1} admits
a unique adapted solution $(P_k,\Lambda_k)\in
L^\infty_\dbF(\Om;C([0,T];\cS(\dbC^{m_k})))\times
L^2_\dbF(0,T;\cS(\dbC^{m_k}))$ with $P_k\geq 0$.

\vspace{0.15cm}

Put
$$\Phi_k=(R_k+P_k)^{-1}, \qq \Psi_k =
-\Phi_k(\Lambda_k+R_{2,k})\Phi_k^{-1}.
$$
Then it is easy to see that $(\Phi_k,\Psi_k)$
solves the following equation (e.g. \cite{QZ}):
\begin{equation}\label{10.25-eq2}
\left\{
\begin{array}{ll}\ds
d\Phi_k =
\big[(A_k+A_{1,k}-B_kC_k)\Phi_k+\Phi_k(A_k+A_{1,k}-B_kC_k)^\top
\\
\ns\ds \qq\q\; - B_k\Phi_kB_k^\top+ B_k \Psi_k +
\Psi_kB_k^\top  + \Phi_k \wt Q_k \Phi_k \big]dt+
\Psi_k dW(t) &\mbox{ in }[0,T],\cr
\ns\ds \Phi_k(T)=(R_k(T)+G_k)^{-1}.
\end{array}
\right.
\end{equation}
Since $K_k$ is a bounded matrix-valued process
with a positive lower bound, we deduce that
$X\in
L^\infty_\dbF(\Om;C([0,T];\cS(\dbR^{m_k})))$.

\vspace{0.1cm}

Next, we prove the following result:
\begin{proposition}\label{10.25-prop1}
Let ({\bf AS7})--({\bf AS9}) hold. Then,
$(\Phi_k,\Psi_k)\in
\dbL_{1,2}^a(\dbC^{m_k})\times
\dbL_{1,2}^a(\dbC^{m_k})$. Furthermore, there
exists a version of $(\mathbf{D}_\tau
\Phi_k,\mathbf{D}_\tau \Psi_k)$ which solves the
following equation:
\begin{equation}\label{10.25-prop1-eq1}
\begin{cases}\ds
d\mathbf{D}_\tau \Phi_k= \mathbf{h}_kdt +  h_kdt
+ \mathbf{D}_\tau \Psi_k dW(t) &\mbox{
in }[\tau,T],\\
\ns\ds d\mathbf{D}_\tau
\Phi_k(T)=\mathbf{D}_\tau (R_k(T)+G_k)^{-1},
\end{cases}
\end{equation}
where
$$
\begin{array}{ll}\ds
\mathbf{h}_k=(A_k+A_{1,k}-B_kC_k)\mathbf{D}_\tau
\Phi_k+\mathbf{D}_\tau
\Phi_k(A_k+A_{1,k}-B_kC_k)^\top
 \\
\ns\ds \qq - B_k\mathbf{D}_\tau \Phi_k
B_k^{\top}+ B_k \mathbf{D}_\tau \Psi_k +
\mathbf{D}_\tau Z_kB_k^{\top}  + \mathbf{D}_\tau
\Phi_k \wt Q_k \Phi_k + \Phi_k \wt
Q_k\mathbf{D}_\tau \Phi_k
\end{array}
$$
and
$$
\begin{array}{ll}\ds
h_k=\mathbf{D}_\tau(A_{1,k}-B_kC_k)
\Phi_k+\Phi_k\mathbf{D}_\tau
(A_{1,k}-B_kC_k)^\top
- \mathbf{D}_\tau B_k X_kB_k^{\top} \\
\ns\ds \qq - B_k \Phi_k\mathbf{D}_\tau
B_k^{\top}+ \mathbf{D}_\tau B_k \Psi_k +
\Psi_k\mathbf{D}_\tau B_k^{\top}  + \Phi_k
\mathbf{D}_\tau \wt Q_k \Phi_k.
\end{array}
$$
Further,
\begin{equation}\label{10.25-prop1-eq2}
\mathbf{D}_\tau \Phi_k = 0,\q \mathbf{D}_\tau
\Psi_k=0 \qq\mbox{ in }[0,\tau].
\end{equation}
\end{proposition}

{\it Proof of Proposition \ref{10.25-prop1}}\,:
Let $(\Phi_{k,j},\Psi_{k,j})$ be a sequence
defined recursively by $\Phi_{k,0}=\Psi_{k,0}=0$
and
\begin{equation}\label{10.25-eq3}
\left\{
\begin{array}{ll}\ds
d\Phi_{k,n+1} =
\big[(A_k+A_{1,k}-B_kC_k)\Phi_{k,n+1}+\Phi_{k,n+1}(A_k+A_{1,k}-B_kC_k)^\top
\\[1mm]
\ns\ds \q\qq\qq - B^k\Phi_{k,n+1}B_k^{\top}+ B_k
\Psi_{k,n+1} + \Psi_{k,n+1}B_k^{\top} +
\Phi_{k,n+1} \wt Q_k \Phi_{k,n} \\[1mm]
\ns\ds \q\qq\qq+ \Phi_{k,n}\wt Q_k \Phi_{k,n+1}
\big]dt+ \Psi_{k,n+1} dW(t) \qq\mbox{ in
}[0,T],\\[1mm]
\ns\ds \Phi_{k,n+1}(T)=(R_k(T)+G_k)^{-1}.
\end{array}
\right.
\end{equation}
By Lemma 4.2 in \cite{QZ}, we conclude that
\begin{equation}\label{10.25-eq12}
\lim_{n\to\infty}(\Phi_{k,n},\Psi_{k,n})=
(\Phi_k,\Psi_k)  \mbox{ in
}\;L^\infty_\dbF(\Om;C([0,T];\cS(\dbC^{m_k})))\times
L^2_\dbF(0,T;\cS(\dbC^{m_k})).
\end{equation}
Clearly, $$(\Phi_{k,0},\Psi_{k,0})\in
\dbL_{1,2}^a(\dbC^{m_k})\times
\dbL_{1,2}^a(\dbC^{m_k}).$$  Suppose that
$$(\Phi_{k,n},\Psi_{k,n})\in
\dbL_{1,2}^a(\dbC^{m_k})\times
\dbL_{1,2}^a(\dbC^{m_k}).$$ Then, by Proposition
\ref{prop2}, we see that
$$(\Phi_{k,n+1},\Psi_{k,n+1})\in
\dbL_{1,2}^a(\dbC^{m_k})\times
\dbL_{1,2}^a(\dbC^{m_k}).$$ Thus,
$$(\Phi_{k,n},\Psi_{k,n})\in
\dbL_{1,2}^a(\dbC^{m_k})\times
\dbL_{1,2}^a(\dbC^{m_k}),\qq \forall n\in\dbN.$$

Next, let us show that
\begin{equation}\label{10.25-eq4}
\lim_{n\to\infty}(\mathbf{D}_\tau
\Phi_{k,n},\mathbf{D}_\tau \Psi_{k,n}) =
(\widehat \Phi_{k,\tau}, \widehat \Psi_{k,\tau})
\q\mbox{ in }\;\dbL_{1,2}^a(\dbC^{m_k})\times
\dbL_{1,2}^a(\dbC^{m_k}).
\end{equation}
Here $(\widehat \Phi_{k,\tau}, \widehat
\Psi_{k,\tau})$ solves the following equation:
\begin{equation}\label{10.25-eq5}
\begin{cases}\ds
d\widehat \Phi_{k,\tau}=\mathbf{f_k}dt  + f_k
dt+ \widehat \Psi_{k,\tau} dW(t) &\mbox{
in }\;[\tau,T],\\[1mm]
\ns\ds d\widehat
\Phi_{k,\tau}(T)=\mathbf{D}_\tau
(R_k(T)+G_k)^{-1}.
\end{cases}
\end{equation}
where
\begin{equation}\label{10.25-eq8}
\begin{array}{ll}\ds
\mathbf{f_k}= (A_k+A_{1,k}-B_kC_k)\widehat
\Phi_{k,\tau}+\widehat
\Phi_{k,\tau}(A_k+A_{1,k}-B_kC_k)^\top
\\[1mm]
\ns\ds \qq - B_k\widehat
\Phi_{k,\tau}B_k^{\top}+ B_k \widehat
\Psi_{k,\tau} + \widehat \Psi_{k,\tau}B_k^{\top}
+ \widehat \Phi_{k,\tau} \wt Q_k \Phi_k + \Phi_k
\wt Q_k \widehat \Phi_{k,\tau}
\end{array}
\end{equation}
and
\begin{equation}\label{10.25-eq9}
\begin{array}{ll}\ds
f_k=\mathbf{D}_\tau(A_{1,k}-B_kC_k)
\Phi_k+\Phi_k\mathbf{D}_\tau
(A_{1,k}-B_kC_k)^\top
- \mathbf{D}_\tau B_k \Phi_kB_k^{\top} \\[1mm]
\ns\ds \qq  - B_k \Phi_k\mathbf{D}_\tau
B_k^{\top}+ \mathbf{D}_\tau B_k \Psi_k +
\Psi_k\mathbf{D}_\tau B_k^{\top}  + \Phi_k
\mathbf{D}_\tau \wt Q_k \Phi_k.
\end{array}
\end{equation}
From \eqref{10.25-eq3}, it is easy to see that
$(\mathbf{D}_\tau \Phi_{k,n},\mathbf{D}_\tau
\Psi_{k,n})$ solves
\begin{equation}\label{10.25-eq6}
\begin{cases}\ds
d\mathbf{D}_\tau \Phi_{k,n}= \mathbf{g}_{k,n}dt
+ g_{k,n} dt+ \mathbf{D}_\tau \Psi_{k,n} dW(t)
&\mbox{
in }[\tau,T],\\[1mm]
\ns\ds d\widehat
\Phi_{k,\tau}(T)=\mathbf{D}_\tau
(R_k(T)+G_k)^{-1}.
\end{cases}
\end{equation}
where
\begin{equation}\label{10.25-eq10}
\begin{array}{ll}\ds
\mathbf{g}_{k,n}=(A_k+A_{1,k}-B_kC_k)\mathbf{D}_\tau
\Phi_{k,n}+\mathbf{D}_\tau
\Phi_{k,n}(A_k+A_{1,k}-B_kC_k)^\top
 \\[1mm]
\ns\ds \qq\q - B_k\mathbf{D}_\tau
\Phi_{k,n}B_k^{\top}+ B_k \mathbf{D}_\tau
\Psi_{k,n} + \mathbf{D}_\tau \Psi_{k,n}
B_k^{\top} + \mathbf{D}_\tau \Phi_{k,n} \wt Q_k
\Phi_{k,n-1}\\[1mm]
\ns\ds \qq\q + \Phi_{k,n-1} \wt Q_k
\mathbf{D}_\tau \Phi_{k,n}
\end{array}
\end{equation}
and
\begin{equation}\label{10.25-eq11}
\begin{array}{ll}\ds
g_{k,n}=\mathbf{D}_\tau(A_{1,k}-B_kC_k)
\Phi_{k,n}+\Phi_{k,n}\mathbf{D}_\tau
(A_{1,k}-B_kC_k)^\top\\[1mm]
\ns\ds \qq\q
- \mathbf{D}_\tau B_k \Phi_{k,n}B_k^{\top} - B_k \Phi_k\mathbf{D}_\tau B_k^{\top}+ \mathbf{D}_\tau B_k \Psi_{k,n}\\[1mm]
\ns\ds \qq\q   + \Psi_{k,n}\mathbf{D}_\tau
B_k^{\top}  + \Phi_{k,n} \mathbf{D}_\tau (\wt
Q_k \Phi_{k,n-1})+ \mathbf{D}_\tau(\Phi_{k,n-1}
\wt Q_k) \Phi_{k,n}.
\end{array}
\end{equation}
By the equations \eqref{10.25-eq5} and
\eqref{10.25-eq6}, we deduce that
\begin{equation}\label{10.25-eq16}
\begin{cases}\ds
d\big(\mathbf{D}_\tau \Phi_{k,n}  -  \widehat
\Phi_{k,\tau}\big) = (\mathbf{g}_{k,n} -
\mathbf{f}_{k,n})dt + (g_{k,n}-f_{k,n}) dt\\[1mm]
\ns\ds\hspace{3.32cm} + \big(\mathbf{D}_\tau
\Psi_{k,n}-\widehat \Psi_{k,\tau}\big) dW(t)
&\mbox{
in }[\tau,T],\\[1mm]
\ns\ds d\widehat
\Phi_{k,\tau}(T)=\mathbf{D}_\tau
(R_k(T)+G_k)^{-1}.
\end{cases}
\end{equation}
Then, for all $\tau\in [0,T]$,
\begin{eqnarray}\label{10.25-eq14}
&&\3n |\mathbf{D}_\tau \Phi_{k,n}-\widehat
\Phi_{k,\tau}|^2_{L^2_\dbF(0,T;\dbR^{m_k\times
m_k})}+ |\mathbf{D}_\tau \Psi_{k,n}-\widehat
\Psi_{k,\tau}|^2_{L^2_\dbF(0,T;\dbR^{m_k\times
m_k})}\nonumber
\\[1mm]
&&\3n \leq \cC \mE\[ \int_\tau^T
\big(|\mathbf{D}_\tau \Phi_{k,n} \wt Q_k
(\Phi_{k,n-1} -  \Phi_k) + (\Phi_{k,n-1}-
\Phi_k) \wt Q_k \mathbf{D}_\tau
\Phi_{k,n}|_{\dbR^{m_k\times m_k}}\nonumber\\[1mm]
&&\3n \qq +|g_{k,n}-g_{k,n}|_{\dbR^{m_k\times
m_k}} \big)ds\]^2.
\end{eqnarray}
It follows from \eqref{10.25-eq12} that
$$
\ba{ll}
\ds\lim_{n\to\infty}\mE\[\int_\tau^T\big(|\mathbf{D}_\tau
\Phi_{k,n} \wt Q_k (\Phi_{k,n-1}- \Phi_k)+
(\Phi_{k,n-1}- \Phi_k) \wt Q_k \mathbf{D}_\tau
\Phi_{k,n}|_{\dbR^{m_k\times m_k}}\\
\ns\ds\qq\q +|g_{k,n}-g_{k,n}|_{\dbR^{m_k\times
m_k}} \big)ds\]^2 =0. \ea
$$
This, together with \eqref{10.25-eq14}, implies
that
\begin{equation}\label{10.25-eq12xx}
\lim_{n\to\infty}\big(|\mathbf{D}_\tau
\Phi_{k,n}-\widehat
\Phi_{k,\tau}|^2_{L^2_\dbF(0,T;\dbR^{m_k\times
m_k})}+ |\mathbf{D}_\tau \Psi_{k,n}-\widehat
\Psi_{k,\tau}|^2_{L^2_\dbF(0,T;\dbR^{m_k\times
m_k})}\big) =0.
\end{equation}
Since $\dbL_{1,2}^a(\dbC^{m_k})\times
\dbL_{1,2}^a(\dbC^{m_k})$ is closed, the limit
$(\Phi_k,\Psi_k)$ of
$\{(\Phi_{k,n},\Psi_{k,n})\}_{n=1}^\infty$
belongs to $\dbL_{1,2}^a(\dbC^{m_k})\times
\dbL_{1,2}^a(\dbC^{m_k})$ and $(\widehat
\Phi_{k,\tau},\widehat \Psi_{k,\tau})$ is a
version of $(\mathbf{D}_\tau
\Phi_k,\mathbf{D}_\tau \Psi_k)$.

\vspace{0.1cm}

Now we prove that for the above version of the
Malliavinian derivative of $\Phi$,
$\mathbf{D}_\tau \Phi(\tau)=\Psi(\tau)$ for a.e.
$\tau\in [0,T]$.

\vspace{0.1cm}

For $t<\tau$, it follows from \eqref{10.25-eq2}
that
$$
\begin{array}{ll}\ds
\Phi_k(\tau) =\Phi_k(t) + \int_t^\tau
\big[(A_k+A_{1,k}-B_kC_k)\Phi_k+\Phi_k(A_k+A_{1,k}-B_kC_k)^\top
 \\
\ns\ds \hspace{3.82cm} - B_k\Phi_kB_k^{\top}+
B_k \Psi_k + \Psi_kB_k^{\top} + \Phi_k \wt Q_k
\Phi_k \big]ds\\
\ns\ds \hspace{1.82cm}+ \int_t^\tau \Psi_k
dW(s).
\end{array}
$$
This, together with Lemma \ref{10.25-lm1},
implies that for $t< \si\leq\tau$,
\begin{equation}\label{10.25-eq20}
\begin{array}{ll}\ds
\mathbf{D}_\si \Phi_k(\tau)\\[1mm]
\ns\ds =\Psi_k(\si) + \int_\si^\tau
\big[(A_k+A_{1,k}-B_kC_k)\mathbf{D}_\si
\Phi_{k}+\mathbf{D}_\si
\Phi_{k}(A_k+A_{1,k}-B_kC_k)^\top
\\[1mm]
\ns\ds \qq\qq\q\qq - B_k\mathbf{D}_\si
\Phi_{k}B_k^{\top} + B_k \mathbf{D}_\si \Psi_{k}
+ \mathbf{D}_\si \Psi_{k}B_k^{\top}  +
\mathbf{D}_\si \Phi_k\wt Q_k \Phi_k
\\[1mm]
\ns\ds  \qq\qq\q\qq + \Phi_k \wt Q_k
\mathbf{D}_\si
\Phi_k+\mathbf{D}_\si(A_{1,k}-B_kC_k)
\Phi_k+\Phi_k\mathbf{D}_\si
(A_{1,k}-B_kC_k)^\top
 \\[1mm]
\ns\ds \qq\qq\q\qq - \mathbf{D}_\si B_k
\Phi_kB_k^{\top}- B_k \Phi_k\mathbf{D}_\si
B_k^{\top}+ \mathbf{D}_\si B_k \Psi_k +
\Psi_k\mathbf{D}_\si B_k^{\top} \\[1mm]
\ns\ds \qq\qq\q\qq  + \Phi_k \mathbf{D}_\si \wt
Q_k \Phi_k\big]ds + \int_\si^\tau \mathbf{D}_\si
\Psi_k dW(s).
\end{array}
\end{equation}
By taking $\si=\tau$, we find that
$\mathbf{D}_\tau \Phi_k(\tau)=\Psi_k(\tau)$,
$\dbP$-a.s., for a.e. $\tau\in [0,T]$.

\ms

With the aid of Proposition \ref{10.25-prop1},
we can prove the following result:
\begin{proposition}\label{10.25-prop2}
Let ({\bf AS7})--({\bf AS9}) hold. Then the
solution $$(P_k,\Lambda_k)\in
L^\infty_\dbF(0,T;\cS(\dbC^{m_k})) \times
L^\infty_\dbF(0,T; \cS(\dbC^{m_k})).$$
\end{proposition}
{\it Proof}\,: We only need to show that
$$\Lambda_k \in
L^\infty_\dbF(0,T;\cS(\dbC^{m_k})).$$ Since
$\Psi_k =
-\Phi_k(\Lambda_k+R_{2,k})(\Phi_k)^{-1}$, it
suffices to prove that $$\Psi_k \in
L^\infty_\dbF(0,T;\dbC^{m_k\times m_k}).$$

For any  $\a>0$, by It\^o's formula, we have
that
\begin{eqnarray}\label{10.25-eq17}
&&\3n\3n e^{\a t}|\mathbf{D}_\tau
\Phi_k(t)|_{\dbC^{m_k\times m_k}}^2 +\int_t^T
e^{\a s}\big(\a|\mathbf{D}_\tau \Phi_k
|_{\dbC^{m_k\times m_k}}^2+|\mathbf{D}_\tau
\Psi_k |_{\dbC^{m_k\times m_k}}^2\big)ds\nonumber\\
&&\3n\3n = e^{\a T}|\mathbf{D}_\tau
(R_k(T)+G_k)^{-1}|_{\dbC^{m_k\times m_k}}^2 -
2\int_t^T e^{\a s}\big\langle \mathbf{D}_\tau
\Phi_k, \mathbf{D}_\tau
\Psi_k\big\rangle_{\dbC^{m_k\times m_k}}dW(s)\nonumber\\
&&\3n\3n\q - 2\int_t^T e^{\a s}\big\langle
\mathbf{D}_\tau \Phi_k, \mathbf{h}_k
\big\rangle_{\dbC^{m_k\times m_k}}ds - 2\int_t^T
e^{\a s}\big\langle \mathbf{D}_\tau \Phi_k, h_k
\big\rangle_{\dbC^{m_k\times m_k}}ds \\
&&\3n\3n \leq e^{\a T} |\mathbf{D}_\tau
(R_k(T)+G_k)^{-1}|_{\dbC^{m_k\times m_k}}^2 -
2\int_t^T e^{\a s}\big\langle \mathbf{D}_\tau
\Phi_k, \mathbf{D}_\tau
\Psi_k\big\rangle_{\dbR^{m_k\times
m_k}}dW(s) \nonumber\\
&&\3n\3n\q  + \cC\int_t^T e^{\a s}
\big(|\mathbf{D}_\tau \Phi_k|_{\dbR^{m_k\times
m_k}}^2+|\mathbf{D}_\tau
\Phi_k|_{\dbC^{m_k\times m_k}}|\mathbf{D}_\tau
\Psi_k|_{\dbC^{m_k\times m_k}} +
|h_k|^2_{\dbC^{m_k\times m_k}}\big)ds\nonumber
\\
&&\3n\3n \leq e^{\a T} |\mathbf{D}_\tau
(R_k(T)+G_k)^{-1}|_{\dbC^{m_k\times m_k}}^2 -
2\int_t^T e^{\a s}\big\langle \mathbf{D}_\tau
\Phi_k, \mathbf{D}_\tau
\Psi_k\big\rangle_{\dbR^{m_k\times
m_k}}dW(s) \nonumber\\
&&\3n\3n\q  + \cC\int_t^T e^{\a s}
|\mathbf{D}_\tau \Phi_k|_{\dbR^{m_k\times
m_k}}^2 ds + \frac{1}{2}\int_t^T e^{\a s}
|\mathbf{D}_\tau \Psi_k|_{\dbR^{m_k\times
m_k}}^2 ds \nonumber\\
&&\3n\3n\q + \cC \int_t^T e^{\a s}
|h_k|^2_{\dbC^{m_k\times m_k}}ds.\nonumber
\end{eqnarray}
By taking $\a=2\cC$, we get from
\eqref{10.25-eq17} that
$$
\!\!\!\begin{array}{ll}\ds e^{\a
t}|\mathbf{D}_\tau \Phi_k(t)|_{\dbR^{m_k\times
m_k}}^2 +\int_t^T e^{\a
s}\big(\a|\mathbf{D}_\tau
\Phi_k|_{\dbR^{m_k\times
m_k}}^2+|\mathbf{D}_\tau
\Psi_k|_{\dbC^{m_k\times m_k}}^2\big)ds\\
\ns\ds   \leq  e^{\a T} |\mathbf{D}_\tau (R_k(T)
+ G_k)^{-1}|_{\dbC^{m_k\times m_k}}^2 - 2
\int_t^T e^{\a s}\big\langle \mathbf{D}_\tau
\Phi_k, \mathbf{D}_\tau
\Psi_k\big\rangle_{\dbC^{m_k\times m_k}}dW(s)\\
\ns\ds\q + \cC \int_t^T e^{\a s}
|h_k|_{\dbC^{m_k\times m_k}}^2ds,
\end{array}
$$
which implies that
\begin{equation*}\label{10.25-eq18}
\begin{array}{ll}\ds
e^{\a t}|\mathbf{D}_\tau
\Phi(t)|_{\dbC^{m_k\times m_k}}^2\\
\ns\ds \leq e^{\a T} \mE\big(|\mathbf{D}_\tau
(R_k(T)+G_k)^{-1}|_{\dbR^{m_k\times
m_k}}^2|\cF_t\big) + \cC\mE\(\int_t^T e^{\a s}
|h_k|_{\dbC^{m_k\times m_k}}^2ds\Big|\cF_t\).
\end{array}
\end{equation*}
Since  $|\mathbf{D}_\tau
(R_k(T)+G_k)^{-1}|_{\dbC^{m_k\times m_k}}^2$ and
$|h_k|_{\dbC^{m_k\times m_k}}^2$ are bounded, we
complete the proof.

\ms

Now we are in a position to prove Theorem
\ref{10.25-th2}.

\ms

{\it Proof of Theorem \ref{10.25-th2}}:\, Let
$$P={\rm diag}(P_1,P_2,\cds,P_k,\cds),\qq\Lambda={\rm diag}(\Lambda_1,\Lambda_2,\cds,\Lambda_k,\cds).
$$
Clearly,
\begin{equation}\label{10.25-eq23}
(P,\Lambda)\in
C_\dbF([0,T];L^\infty(\Om;\cL(H)))\times
L^\infty_\dbF(0,T;\cL(H)).
\end{equation}
Since  for all $k\in\dbN$, $$P_k(t,\om)\geq 0\
\mbox{ for a.e. }(t,\om)\in [0,T]\times\Om,$$ we
get that $$P\geq 0\ \mbox{ for a.e. }(t,\om)\in
[0,T]\times\Om.$$ Thus,
$$K(t,\om)=R(t,\om)+D(t,\om)^*P(t,\om)D(t,\om)>0
\mbox{ for a.e. }(t,\om)\in [0,T]\times\Om.$$
Similar to the argument of Step 8 in the proof
of Theorem \ref{5.7-th1}, one can easily obtain
that the left inverse $K(t,\om)^{-1}$ is a
densely defined closed operator for a.e.
$(t,\om)\in [0,T]\times\Om$.

Write $\ds n_j\triangleq\sum_{k=1}^{j} m_k$,
where $j\in\dbN$. For any $k\in\dbN$,  $t\in
[0,T]$, $\xi_1,\xi_2\in L^4_{\cF_t}(\Om;H)$,
$u_1(\cd), u_2(\cd)\in L^4_\dbF(\Om;L^2(t,T;H))$
and $v_1(\cd), v_2(\cd)$ $\in
L^4_\dbF(\Om;L^2(t,T;V'))$, and for $\ell=1,2$,
let
$$
\begin{cases}\ds
\xi_{\ell,k}=(\G_{n_{k+1}}-\G_{n_{k}})\xi_\ell,
\\
\ns\ds
u_{\ell,k}=(\G_{n_{k+1}}-\G_{n_{k}})u_\ell, \\
\ns\ds
v_{\ell,k}=(\G_{n_{k+1}}-\G_{n_{k}})v_\ell,\\
\ns\ds
x_{\ell,k}=(\G_{n_{k+1}}-\G_{n_{k}})x_\ell.
\end{cases}
$$
By It\^o's formula, we obtain that
\begin{eqnarray*}\label{10.25-eq21}
&&\3n\3n\mE\langle G^kx_{1,k}(T),
x_{2,k}(T)\rangle_{H} +\mE \int_t^T \big\langle
Q_k(\tau) x_{1,k}(\tau),
x_{2,k}(\tau) \big\rangle_{H}d\tau\nonumber\\
&&\3n\3n\q - \mE \int_t^T \big\langle
(K_k(\tau))^{-1} L_k(\tau) x_{1,k}(\tau),
L_k(\tau)x_{2,k}(\tau) \big\rangle_{H}d\tau
\\
&&\3n\3n = \mE\big\langle P_k(t)
\xi_{1,k},\xi_{2,k} \big\rangle_{H} + \mE
\int_t^T \big\langle P_k(\tau)u_{1,k}(\tau),
x_{2,k}(\tau)\big\rangle_{H}d\tau  \nonumber\\
&&\3n\3n \q + \mE \int_t^T\! \big\langle
P_k(\tau)x_{1,k}(\tau),
u_{2,k}(\tau)\big\rangle_{H}d\tau + \mE\!
\int_t^T\!\big\langle
P_k(\tau)C_k(\tau)x_{1,k}(\tau),
v_{2,k}(\tau)\big\rangle_{H}d\tau \nonumber\\
&&\3n\3n \q + \mE
\int_t^T \big\langle  P_k(\tau)v_{1,k}(\tau), C_k(\tau)x_{2,k}(\tau)+v_{2,k}(\tau)\big\rangle_{H}d\tau \nonumber\\
&&\3n\3n \q + \mE \int_t^T \big\langle
v_{1,k}(\tau),
\Lambda_k(\tau)x_{2,k}(\tau)\big\rangle_{
H}d\tau+ \mE \int_t^T \big\langle
\Lambda_k(\tau)x_{1,k}(\tau), v_{2,k}(\tau)
\big\rangle_{H}d\tau.\nonumber
\end{eqnarray*}
This, together with definitions of $P$ and
$\Lambda$, implies that $(P,\Lambda)$ satisfies
\eqref{6.18-eq1}. Similarly, we can show that
they fulfills \eqref{10.10-eq10}. Thus,
$(P,\Lambda)$ is a transposition solution to
\eqref{5.5-eq6}.


\section{Some examples of controlled SPDEs}\label{sect7}


In this chapter, we shall give some illuminating
examples for SLQs for stochastic wave, parabolic
and Schr\"odinger equations. One can also
consider SLQs for other controlled stochastic
partial differential equations, such as
stochastic KdV equations, stochastic beam
equations, etc. We omit it here.

\ms In this chapter, we let $n\in\dbN$ and
$\cO\subset\dbR^n$ be a bounded domain with a
$C^\infty$ boundary $\pa\cO$.

\subsection{SLQs for stochastic wave equations}

Stochastic wave equations are widely used to
describe vibrations of string/ membrane under
the perturbations of random noises and
propagation of waves in random media (e.g.
\cite{Funaki,Kotelenez}). In this section, we
consider the SLQs for the following controlled
stochastic wave equations:
\begin{equation}\label{8.25-eq14}
\begin{cases}
dy_t - \D ydt = (a_1y+ b_1u)dt + (a_2y + b_2
u)dW(t)
&\mbox{ in } \cO\times (0,T),\\[1mm]
\ns\ds y=0 &\mbox{ on }\pa\cO\times (0,T),\\[1mm]
\ns\ds y(0)=y_0,\q y_t(0)=y_1 &\mbox{ in } \cO,
\end{cases}
\end{equation}
with the following cost functional
$$
\cJ(0,y_0,y_1;u)\triangleq \mE\int_0^T\int_\cO
(q|y|^2 + r|u|^2)dxdt + \mE\int_\cO g|y(T)|^2dx.
$$
Here $(y_0,y_1)\in H_0^1(\cO)\cap L^2(\cO)$,
\begin{equation}\label{8.25-eq12}
\begin{cases}
\ds  a_1,a_2,b_1,b_2\in
L^\infty_{\dbF}(0,T;C^{2n}(\overline\cO)),\\[1mm]
\ns\ds
u\in L^2_\dbF(0,T;L^2(\cO)),  \\[1mm]
\ns\ds q,r \in
L^\infty_{\dbF}(0,T;C^{2n}(\overline\cO)), \\[1mm]
\ns\ds g\in
L^\infty_{\cF_T}(\Om;C^{2n}(\overline\cO)),\\[1mm]
\ns\ds q\geq 0,\;\;r\geq 1,\;\;g\geq 0.
\end{cases}
\end{equation}

\ms

Our optimal control problem is as follows:

\ms

\no\bf Problem (wSLQ). \rm For each
$(y_0,y_1)\in H_0^1(\cO)\times L^2(\cO)$, find a
$\bar u(\cd)\in L^2_\dbF(0,T;$ $L^2(\cO))$ such
that
\begin{equation}\label{5.2-eq3xxx}
\cJ\big(0,y_0,y_1;\bar
u(\cd)\big)=\inf_{u(\cd)\in
L^2_\dbF(0,T;L^2(\cO))}\cJ\big(0,y_0,y_1;u(\cd)\big).
\end{equation}

Problem (wSLQ) is a concrete example of Problem
(SLQ) with the following setting:

\begin{itemize}
  \item $H=H_0^1(\cO)\times L^2(\cO)$ and $U=L^2(\cO)$;

\ms

  \item The operator $A$ is
defined as follows:
$$
\begin{cases}
D(A)=[H^2(\cO)\cap H_0^1(\cO)]\times H_0^1(\cO),\\[1mm]
\ns\ds A\left(
          \begin{array}{c}
            \f_1 \\
            \f_2 \\
          \end{array}
        \right)
  = \left(
          \begin{array}{c}
            \f_2 \\
            \D \f_1 \\
          \end{array}
        \right),\q \forall \left(
          \begin{array}{c}
            \f_1 \\
            \f_2 \\
          \end{array}
        \right)\in D(A);
\end{cases}
$$

\ms

\item $A_1y=(0,a_1y)^\top$, $Bu=(0,b_1u)^\top$, $Cy=(0,a_2y)^\top$ and $Du=(0,b_2u)^\top$;

\ms

\item The operators $Q$, $R$ and $G$ are given by
$$
\left\{ \ba{ll} \ds\langle Q y,y\rangle_H =
\int_\cO q(|\nabla y|^2 + |y_t|^2)dx,\\[2mm]
\ns\ds\langle R u,u\rangle_H = \int_\cO
r|u|^2dx,\\[2mm] \ns\ds \langle G y(T),y(T)\rangle_H
= \int_\cO g(|\nabla y(T)|^2+|y_t(T)|^2)dx. \ea
\right.
$$
\end{itemize}

Since $q\geq 0$, $r\geq 1$ and $g\geq 0$, it is
easy to check that ({\bf AS1}) holds.

\vspace{0.1cm}

Define an operator $\widehat A$ as follows:
$$
\begin{cases}
D(\widehat A)=H^2(\cO)\cap H_0^1(\cO),\\[1mm]
\ns\ds \widehat A\f =  -\D\f,\q \forall \f\in
D(\widehat A).
\end{cases}
$$
Denote by  $\{\hat \lambda_j\}_{j=1}^\infty$ the
eigenvalues of $\widehat A$ and $\{ \hat
e_j\}_{j=1}^\infty$ the corresponding
eigenvectors with $|\hat e_j|_{L^2(\cO)}=1$ for
$j\in\dbN$.  Clearly, $\big\{\pm
i\sqrt{\hat\lambda_j}\big\}_{j=1}^\infty$ are
the eigenvalues of $A$ and $\ds\Big\{\Big(\pm
\frac{1}{i\sqrt{\hat \lambda_j}} \hat e_j,\hat
e_j\Big)\Big\}_{j=1}^\infty$ are the
corresponding eigenvectors.  It is well known
that $\ds\Big\{\Big(\pm \frac{1}{i\sqrt{\hat
\lambda_j}} \hat e_j,\hat
e_j\Big)\Big\}_{j=1}^\infty$ constitutes an
orthonormal basis of $H^1_0(\cO)\times
L^2(\cO)$. Hence, ({\bf AS2}) holds.

\vspace{0.1cm}

Denote by $V$ the completion of the Hilbert
space $H^1_0(\cO)\times L^2(\cO)$ with the norm
$$
\begin{array}{ll}\ds
|(f_1,f_2)|_V^2 = \sum_{j=1}^\infty |\hat
\lambda_j|^{-n}\big( |f_{1,j}|^2  +
|f_{2,j}|^2\big), \\
\ns\ds \qq\qq \mbox{ for all
}f_1=\sum_{j=1}^\infty f_{1,j} \hat
\lambda_j^{-\frac{1}{2}} e_j \in H_0^1(\cO)\;
\mbox{ and }\; f_2=\sum_{j=1}^\infty f_{2,j} e_j
\in L^2(\cO).
\end{array}
$$
By the asymptotic distribution of the
eigenvalues of $\widehat A$ (e.g., \cite[Chapter
1, Theorem 1.2.1]{Safarov}), we see that
$$
\begin{array}{ll}\ds
|I|^2_{\cL_2(H;V)}\3n&\ds= \sum_{j=1}^\infty
\left(\Big(\pm \frac{1}{i\sqrt{\hat \lambda_j}}
e_j,e_j\Big)^\top, \Big(\pm \frac{1}{i\sqrt{\hat
\lambda_j}} e_j,e_j\Big)^\top\right)_{V}\\
\ns&\ds =4\sum_{j=1}^\infty |\hat
\lambda_j|^{-n}<\infty.
\end{array}
$$
Hence, the embedding from $H_0^1(\cO)\times
L^2(\cO)$ to $V$ is Hilbert-Schmidt. From the
definition of $V$, it follows that
$V'=D(A^{n})$. By the classical theory of
operator semigroup (e.g., \cite[Chapter II,
Section 3]{EN}), one can show that $A$ generates
a $C_0$-semigroup on $V'$ .

\vspace{0.1cm}

Let $\a\in C^{2n}(\overline\cO)$. For any $f\in
V'$, one has
$$
\begin{array}{ll}\ds
|\a f|_{V'}\3n&\ds=|\a f|_{D(A^{n})} = |A^{n}
(\a f)|_{H^1_0(\cO)\times L^2(\cO)}\\[1mm]
\ns&\ds\leq
|\a|_{C^{2n}(\overline\cO)}|f|_{D(A^{n})} =
|\a|_{C^{2n}(\overline\cO)}|f|_{V'}.
\end{array}
$$
From this, we conclude that $C\in
L^\infty_\dbF(0,T;\cL(V'))$, $G\in
L^\infty_{\cF_T}(\Om;\cL(V'))$ and $Q\in
L^\infty_\dbF(0,T;\cL(V'))$. Thus, ({\bf AS3})
holds.

\vspace{0.1cm}

Let $\wt U  =D(\widehat A^n)$. Clearly, $\wt U$
is dense in $L^2(\cO)$.  From the definition of
$B$, $D$ and $R$, we find that $R\in
L^\infty_\dbF(0,T;\cL(\wt U))$ and $B,D\in
L^2_\dbF(0,T;\cL(\wt U;V'))$. Therefore, ({\bf
AS4}) holds.

\vspace{0.2cm}

\subsection{SLQ for stochastic parabolic
equations}

Stochastic parabolic equations are widely used
to describe diffusion processes under the
perturbations of random noises(e.g.
\cite{Arnold,Greenwood}). In this section, we
consider the SLQs for the following controlled
stochastic parabolic equations:
\begin{equation}\label{5.27-eq1}
\begin{cases}
dy -   \D ydt = (a_1y+ b_1u)dt + (a_2y + b_2
u)dW(t)
&\mbox{ in } \cO\times (0,T),\\
\ns\ds y=0 &\mbox{ on }\pa\cO\times (0,T),\\
\ns\ds y(0)=y_0 &\mbox{ in } \cO,
\end{cases}
\end{equation}
with the following cost functional
$$
\cJ(0,y_0;u)\triangleq \mE\int_0^T\int_\cO
(q|y|^2 + r|u|^2)dxdt + \mE\int_\cO g|y(T)|^2dx.
$$
Here $y_0\in L^2(\cO)$,
\begin{equation}\label{5.27-eq2}
\begin{cases} \ds
a_1,\, b_1, \, b_2,\,\frac{1}{b_2},\,q \in
L^\infty_{\dbF}(0,T), \;
a_2=-b_1,\\[2mm]
\ns\ds
u\in L^2_\dbF(0,T;L^2(\cO)),\\
\ns\ds r(t)=1+\int_0^t r_0(s)ds, \qq r_0 \in L^\infty_{\dbF}(0,T),\\
\ns\ds g\in L^\infty_{\cF_T}(\Om),\q g\geq 0,\\[2mm]
\ns\ds q-r_0-b_1^2r-2a_1 r \geq 0,
\end{cases}
\end{equation}
and  the Malliavinian derivatives of $a_1$,
$b_1$, $b_2$, $r_0$, $q$ and $g$ are   uniformly
bounded with respect to $(t,\om)\in
[0,T]\times\Om$.

\ms

\begin{remark}
$\ds\frac{1}{b_2} \in L^\infty_{\dbF}(0,T)$
means that the control in the diffusion term is
acted on the whole domain. It looks like a
little restrictive. Nevertheless, in some
important control systems, such kind of control
can be achieved. For example, in the propagation
of propagation of an electric potential in a
neuron, one can impose such a control by putting
a electronic field on the neuron.
\end{remark}

\ms

Consider the following optimal control problem:

\ms

\no\bf Problem (pSLQ). \rm For each $y_0\in
L^2(\cO)$, find a $\bar u(\cd)\in
L^2_\dbF(0,T;L^2(\cO))$ such that
\begin{equation}\label{5.27-eq3}
\cJ\big(0,y_0;\bar u(\cd)\big)=\inf_{u(\cd)\in
L^2_\dbF(0,T;L^2(\cO))}\cJ\big(0,y_0;u(\cd)\big).
\end{equation}

Problem (pSLQ) is a concrete example of Problem
(SLQ) with the following setting:

\begin{itemize}
  \item $H=U=L^2(\cO)$;

\ms

  \item The operator $A$ is
defined as follows:
$$
\begin{cases}
D(A)=H^2(\cO)\cap H_0^1(\cO),\\[1mm]
\ns\ds A\f =  \D\f,\q \forall \f\in D(A);
\end{cases}
$$
\item $A_1y=a_1y$, $Bu=b_1u$, $Cy=-b_1y$ and $Du=b_2u$;

\ms

\item The operators $Q$, $R$ and $G$ are given by
$$
\begin{cases}\ds
\langle Q y,y\rangle_H = \int_\cO q|y|^2dx, \\[2mm]
\ns\ds \langle R u,u\rangle_H = \int_\cO
r|u|^2dx,\\[2mm]
\ns\ds\langle G y(T),y(T)\rangle_H = \int_\cO
g|y(T)|^2dx.
\end{cases}
$$
\end{itemize}

Clearly, $A$, $A_1$, $B$, $C$, $Q$, $G$ and $R$
are all infinite dimensional  diagonal matrices
and $D$ is invertible. Hence ({\bf AS5}) and
({\bf AS6}) hold. From the third line of
\eqref{5.27-eq2}, we find that ({\bf AS7}) is
satisfied.

\vspace{0.1cm}

From the first line of \eqref{5.27-eq2}, it is
clear that $B=-C$. Hence, $$RB+C^*
R=rb_1I-b_1rI=0.$$ Denote by
$\{\lambda_j\}_{j=1}^\infty$ the eigenvalues of
$A$ and $\{e_j\}_{j=1}^\infty$ the corresponding
eigenvectors with $| e_j|_{L^2(\cO)}=1$ for
$j\in\dbN$.   By \eqref{5.27-eq2} and noting
that $r$ is independent of the spatial variable,
we find that for any $j\in \dbN$,
$$
\begin{array}{ll}\ds
\lan [Q -r_0I + C^*R C  + R (BC-A-A_1) + (BC-A-A_1)^* R] e_j,  e_j \ran{}_H\\[1mm]
\ns\ds = \lan
[(q-r_0+b_1^2r)I+rI(-b_1^2I-\D-a_1)+
(-b_1^2I-\D-a_1)rI] e_j,  e_j \ran{}_H\\[1mm]
\ns\ds =  (q-r_0+b_1^2r) +2r (-b_1^2-a_1
+\lambda_j)=q-r_0-b_1^2r-2a_1 r+\lambda_j > 0.
\end{array}
$$
Therefore, ({\bf AS8}) holds.

\vspace{0.1cm}

Since the Malliavinian derivatives of $a_1$,
$b_1$, $b_2$, $r_0$, $q$ and $g$ are   uniformly
bounded with respect to $(t,\om)\in
[0,T]\times\Om$, we see that ({\bf AS9}) is
satisfied.

\vspace{0.2cm}

\subsection{SLQs for stochastic Schr\"odinger
equations}

Stochastic Schr\"odinger equations are useful
tools to describe open quantum systems(e.g.
\cite{Kol}). In this section, we consider the
SLQs for the following controlled stochastic
parabolic equations:
\begin{equation}\label{10.25-eq25}
\begin{cases}
dy -  i\D ydt = (a_1y+ b_1u)dt + (a_2y + b_2
u)dW(t)
&\mbox{ in } \cO\times (0,T),\\[1mm]
\ns\ds y=0 &\mbox{ on }\pa\cO\times (0,T),\\[1mm]
\ns\ds y(0)=y_0 &\mbox{ in } \cO,
\end{cases}
\end{equation}
with the following cost functional
$$
\cJ(0,y_0;u)\triangleq \mE\int_0^T\int_\cO
(q|y|^2 + r|u|^2)dxdt + \mE\int_\cO g|y(T)|^2dx.
$$
Here $y_0\in L^2(\cO)$. The conditions on the
coefficients will be given below.

\ms

Our optimal control problem is as follows:

\ms

\no\bf Problem (sSLQ). \rm For each $y_0\in
L^2(\cO)$, find a $\bar u(\cd)\in
L^2_\dbF(0,T;L^2(\cO))$ such that\vspace{0.1cm}
\begin{equation}\label{10.25-eq26}
\cJ\big(0,y_0;\bar u(\cd)\big)=\inf_{u(\cd)\in
L^2_\dbF(0,T;L^2(\cO))}\cJ\big(0,y_0;u(\cd)\big).
\end{equation}

 Problem (sSLQ) is a concrete example of Problem
(SLQ) with the following setting:

\begin{itemize}
  \item $H=U=L^2(\cO)$;\vspace{0.1cm}
  \item The operator $A$ is
defined as follows:
$$
\begin{cases}
D(A)=H^2(\cO)\cap H_0^1(\cO),\\[1mm]
\ns\ds A\f = i \D\f,\q \forall \f\in D(A);
\end{cases}
$$
\item $A_1=0$, $Bu=au$, $Cy=-ay$ and $Du=u$;\vspace{0.1cm}
\item The operators $Q$, $R$ and $G$ are given by
$$
\begin{cases}\ds
\langle Q y,y\rangle_H = \int_\cO q|y|^2dx, \\[1mm]
\ns\ds \langle R u,u\rangle_H = \int_\cO
r|u|^2dx,\\[1mm]
\ns\ds\langle G y(T),y(T)\rangle_H = \int_\cO
g|y(T)|^2dx.
\end{cases}
$$
\end{itemize}

To guarantee that ({\bf AS1})--({\bf AS4}) hold,
we assume the coefficients fulfill the following
conditions:
\begin{equation}\label{5.27-eq4}
\begin{cases}
\ds  a_1,a_2,b_1,b_2\in
L^\infty_{\dbF}(0,T;C^{2n}(\overline\cO)),\\[1mm]
\ns\ds
u\in L^2_\dbF(0,T;L^2(\cO)),  \\[1mm]
\ns\ds q,r \in
L^\infty_{\dbF}(0,T;C^{2n}(\overline\cO)), \\[1mm]
\ns\ds g\in
L^\infty_{\cF_T}(\Om;C^{2n}(\overline\cO)),\\[1mm]
\ns\ds q\geq 0,\;\;r\geq 1,\;\;g\geq 0.
\end{cases}
\end{equation}

\ms

Since $q\geq 0$, $r\geq 1$ and $g\geq 0$, it is
easy to check that ({\bf AS1}) holds.

\vspace{0.1cm}

Write  $\{ \lambda_j\}_{j=1}^\infty$ for the
eigenvalues of $A$ and $\{e_j\}_{j=1}^\infty$
the corresponding eigenvectors such that
$|e_j|_{L^2(\cO)}=1$ for $j\in\dbN$.  It is well
known that $\{e_j\}_{j=1}^\infty$ constitutes an
orthonormal basis of $ L^2(\cO)$. Hence, ({\bf
AS2}) holds.

\vspace{0.1cm}

Denote by $V$ the completion of the Hilbert
space $L^2(\cO)$ with the norm
$$
\begin{array}{ll}\ds
|f|_V^2 = \sum_{j=1}^\infty |\lambda_j|^{-n}
|f_{j}|^2,  \qq \mbox{ for all } f
=\sum_{j=1}^\infty f_{j} e_j \in L^2(\cO).
\end{array}
$$
By the asymptotic distribution of the
eigenvalues of $A$ (e.g., \cite[Chapter 1,
Theorem 1.2.1]{Safarov}), we see that
$$
\begin{array}{ll}\ds
|I|^2_{\cL_2(H;V)} = \sum_{j=1}^\infty
\left(e_j,e_j\right)_{V} =4\sum_{j=1}^\infty
|\lambda_j|^{-n}<\infty.
\end{array}
$$
Hence, the embedding from $ L^2(\cO)$ to $V$ is
Hilbert-Schmidt. From the definition of $V$, it
follows that $V'=D((iA)^{n/2})$. By the
classical theory of operator semigroup (e.g.,
\cite[Chapter II, Section 3]{EN}), one can show
that $A$ generates a $C_0$-semigroup on $V'$ .

\vspace{0.1cm}

Let $\a\in C^{2n}(\overline\cO)$. For any $f\in
V'$, one has
$$
\begin{array}{ll}\ds
|\a f|_{V'}\3n&\ds=|\a f|_{D((iA)^{n/2})} =
|(iA)^{n/2}
(\a f)|_{L^2(\cO)}\\[1mm]
\ns&\ds\leq
|\a|_{C^{2n}(\overline\cO)}|f|_{D((iA)^{n/2})} =
|\a|_{C^{2n}(\overline\cO)}|f|_{V'}.
\end{array}
$$
From this, we conclude that $$C\in
L^\infty_\dbF(0,T;\cL(V')), \q  G\in
L^\infty_{\cF_T}(\Om;\cL(V')), \q Q\in
L^\infty_\dbF(0,T;\cL(V')).$$ Thus, ({\bf AS3})
holds.

\vspace{0.1cm}

Let $\wt U  =D((iA)^n)$. Clearly, $\wt U$ is
dense in $L^2(\cO)$.  From the definitions of
$B$, $D$ and $R$, we find that $$R\in
L^\infty_\dbF(0,T;\cL(\wt U)), \q B,D\in
L^2_\dbF(0,T;\cL(\wt U;V')).$$ Therefore, ({\bf
AS4}) holds.

\ms

Next, we give conditions on the coefficients
such that ({\bf AS5})--({\bf AS9}) hold. We
assume that\vspace{0.1cm}
\begin{equation}\label{10.25-eq24}
\begin{cases} \ds
a_1,\, b_1, \, b_2,\,\frac{1}{b_2},\,q \in
L^\infty_{\dbF}(0,T;\dbR), \;
a_2=-b_1,\\[2mm]
\ns\ds
u\in L^2_\dbF(0,T;L^2(\cO)),\\
\ns\ds r(t)=1+\int_0^t r_0(s)ds, \qq r_0 \in L^\infty_{\dbF}(0,T;\dbR),\\
\ns\ds g\in L^\infty_{\cF_T}(\Om;\dbR),\q g\geq 0,\\[2mm]
\ns\ds q-r_0-b_1^2r-2a_1 r \geq 0,
\end{cases}
\end{equation}
and  the Malliavinian derivative of $a_1$,
$b_1$, $b_2$, $r_0$, $q$ and $g$ are   uniformly
bounded with respect to $(t,\om)\in
[0,T]\times\Om$.

\ms

Clearly, $A$, $A_1$, $B$, $C$, $Q$, $G$ and $R$
are all infinite dimensional  diagonal matrices
and $D$ is invertible. Hence ({\bf AS5}) and
({\bf AS6}) hold. From \eqref{10.25-eq24}, we
find that ({\bf AS7}) is satisfied.

\vspace{0.1cm}

From the first line of \eqref{10.25-eq24}, we
see that $B=-C$. Consequently,\vspace{0.1cm}
$$RB+C^* R=raI-arI=0.$$
By \eqref{10.25-eq24}, and noting that $r$ is
independent of the spatial variable and
$A^*=-A$, we find that for any $j\in
\dbN$,\vspace{0.1cm}
$$
\begin{array}{ll}\ds
\lan [Q -r_0I + C^*R C  + R (BC-A) + (BC-A)^* R] e_j, e_j \ran{}_H\\[1mm]
\ns\ds = \lan [(q-r_0+a^2r)I+rI(-a^2I-A)+
(-a^2I-A)^*rI] e_j, e_j \ran{}_H \\[1mm]
\ns\ds=q-r_0-a^2r > 0.
\end{array}
$$
Therefore, ({\bf AS8}) holds.

\vspace{0.1cm}

Since the Malliavinian derivative of $a_1$,
$b_1$, $b_2$, $r_0$, $q$ and $g$ are   uniformly
bounded with respect to $(t,\om)\in
[0,T]\times\Om$, we see that ({\bf AS9}) is
satisfied.


\section*{Acknowledgement}


This work is supported by the NSF of China under
grants 11971334 and 11931011, the NSFC-CNRS
Joint Research Project under grant 11711530142
and the PCSIRT under grant IRT$\_$16R53 from the
Chinese Education Ministry, and the Fundamental
Research Funds for the Central Universities in
China under grant 2015SCU04A02.

\end{document}